\documentclass[sn-mathphys,Numbered]{sn-jnl}
\usepackage{soul}
\usepackage{graphicx}
\usepackage{multirow}
\usepackage{amsmath,amssymb,amsfonts}%
\usepackage{mathtools} 
\usepackage{amsthm}%
\usepackage{mathrsfs}%
\usepackage[title]{appendix}%
\usepackage{xcolor}%
\usepackage{textcomp}%
\usepackage{manyfoot}%
\usepackage{subcaption}
 \usepackage{array}
\usepackage{booktabs}%
\usepackage{algorithm}%
\usepackage{algorithmicx}%
\usepackage{algpseudocode}%
\usepackage{listings}%
\usepackage{bm}
\usepackage{geometry}
\usepackage{todonotes}
\usepackage{environ}
\usepackage{csquotes}
\usepackage{circuitikz}
\usepackage{xspace}
\newcommand{\MATLAB}{\textsc{Matlab}\xspace}
\makeatletter
\newcommand*\dashline{\rotatebox[origin=c]{90}{$\dabar@\dabar@\dabar@$}}
\makeatother

\renewcommand{\b}[1]{\mathbf #1}
\newcommand{\byref}{\b{y}_{\textrm{ref}}}

\def\code#1{\mbox{\texttt{#1}}}
\newcommand{\ie}{i.\,e.\ }
\newcommand{\eg}{e.\,g.\ }
\newcommand{\dt}{\Delta t}
\newcommand{\dx}{\Delta x}
\newcommand{\ii}{\mathrm{i}}

\newcommand{\tend}{t_\mathrm{end}}
\newcommand{\dd}{\mathrm{d}}

\DeclareMathOperator{\re}{Re}
\DeclareMathOperator{\im}{Im}

\newcommand{\Fj}{\mathbf{F}^{[j]}}

\newcommand{\byk}{\b y^{(k)}}

\newcommand{\Rmax}{R_{\max}}
\newcommand{\Smax}{S_{\max}}
\newcommand{\tol}{\text{tol}}

\newcommand{\step}{\text{Step}}
\newcommand{\Tset}{\text{T}}
\newcommand{\test}{\text{test}}

\newcommand{\bspi}{\boldsymbol{\pi}}
\newcommand{\bssigma}{\boldsymbol{\sigma}}

\renewcommand{\O}{\mathcal O}
\newcommand{\R}{\mathbb R}
\newcommand{\C}{\mathbb C}

\newcommand{\norm}[1]{\lVert #1\rVert}
\DeclareMathOperator{\const}{const}
\DeclareMathOperator{\diag}{diag}
\newcommand{\qta}{\quad\text{ and }\quad}

\renewcommand{\Vec}[1]{\begin{pmatrix*}[r]#1\end{pmatrix*}}

\theoremstyle{thmstyleone}%
\newtheorem{thm}{Theorem}
\newtheorem{prop}[thm]{Proposition}%

\theoremstyle{thmstyletwo}%
\newtheorem{rem}{Remark}%

\theoremstyle{thmstylethree}%
\newtheorem{defn}{Definition}%

\usepackage{tikz}
\usepackage{tkz-euclide}
\usetikzlibrary{patterns}
\definecolor{colorA}{rgb}{0,0.447,0.741}
\definecolor{colorB}{rgb}{0.85,0.325,0.098}
\definecolor{colorE}{rgb}{0.929,0.694,0.125}
\definecolor{colorF}{rgb}{0.494,0.184,0.556}
\definecolor{colorD}{rgb}{0.466,0.674,0.188}
\definecolor{colorC}{rgb}{0.301,0.745,0.933}
\definecolor{colorG}{rgb}{0.635,0.078,0.184}

\newlength{\tickl}    
\tikzset{axes/.style={thick,-latex}}
\tikzset{lineplot/.style={thick}}
\tikzset{arrow/.style={thick,-latex}} 
\tikzset{thick arrow/.style={ultra thick,-latex}} 
\tikzset{grid lines/.style={very thin,gray!30}}	
\tikzset{point/.style={radius=2pt}}
\tikzset{help line/.style={black,thin,dashed}} 

\usetikzlibrary{shapes.geometric}

\makeatletter
\newsavebox{\measure@tikzpicture}
\NewEnviron{scaletikzpicturetowidth}[1]{%
	\def\tikz@width{#1}%
	\def\tikzscale{1}\begin{lrbox}{\measure@tikzpicture}%
		\BODY
	\end{lrbox}%
	\pgfmathparse{#1/\wd\measure@tikzpicture}%
	\edef\tikzscale{\pgfmathresult}%
	\BODY
}
\makeatother

\usetikzlibrary{decorations.markings}
\usepackage{pgfplots}

\tikzset{mylabel/.style  args={at #1 #2  with #3}{
		postaction={decorate,
			decoration={
				markings,
				mark= at position #1
				with  \node [#2] {#3};
} } } }

\begin{document}

	\title[Efficient MPRK Schemes]{Using Machine Learning to Design Time Step Size Controllers for Stable Time Integrators}

\author*[1]{\fnm{Thomas} \sur{Izgin}}\email{izgin@mathematik.uni-kassel.de}

\author[2]{\fnm{Hendrik} \sur{Ranocha}}\email{hendrik.ranocha@uni-mainz.de}


\affil*[1]{\orgdiv{Department of Mathematics}, \orgname{University of Kassel}, \orgaddress{\street{Heinrich--Plett--Straße~40}, \city{Kassel}, \postcode{34132}, \country{Germany}}. ORCID:~0000-0003-3235-210X}

\affil[2]{\orgdiv{Institute of Mathematics}, \orgname{Johannes Gutenberg University Mainz}, \orgaddress{\street{Staudingerweg 9}, \city{Mainz}, \postcode{55128}, \country{Germany}}. ORCID:~0000-0002-3456-2277}

\abstract{%

We present a new method for developing time step controllers based on a technique from the field of machine learning. This method is applicable to stable time integrators that have an embedded scheme, i.e., that have local error estimation similar to Runge-Kutta pairs.
	To design good time step size controllers using these
	error estimates, we propose to use Bayesian optimization. In particular,
	we design a novel objective function that captures important properties
	such as tolerance convergence and computational stability.
	We apply our new approach to several modified Patankar--Runge--Kutta (MPRK) schemes and a Rosenbrock-type scheme, equipping them with controllers based
	on digital signal processing which extend classical PI and PID controllers.
	We demonstrate that the optimization process yields controllers that are
	at least as good as the best controllers chosen from a wide range of
	suggestions available for classical explicit and implicit time integration
	methods by providing work-precision diagrams on a variety of ordinary and partial differential equations.
}




	\maketitle

	\section{Introduction}
	Numerically solving systems of ordinary differential equations (ODEs)
	\begin{equation}\label{eq:init}
		\b y'(t)=\b f(\b y(t), t),\quad \b y(0)=\b y^0\in \R^N,
	\end{equation}
	gives rise to many challenges which are addressed theoretically by topics such  as order of accuracy and convergence, stability, efficiency, and the preservation of linear invariants or further properties of the solution such as positivity.
	Runge--Kutta (RK) methods are commonly used examples of numerical methods solving the initial value problem \eqref{eq:init}, and all the above mentioned topics have been excessively studied in the past, see for instance \cite{B16, HNW1993, HNWII}.

	Traditionally, the first step size control mechanism for numerical approximations
	of ODEs uses an estimate of the local error and multiplies the current time step
	size by a factor derived from the estimate to maximize the time step size while
	keeping the error below a given tolerance \cite[Section~II.4]{HNW1993}. This approach is basically a deadbeat
	I controller and can be improved, in particular when the method is not operating
	in the asymptotic regime of small time step sizes. By including some history
	of the control process, improved PI controllers were introduced in
	\cite{gustafsson1988pi,gustafsson1991control} and extended to more general
	digital signal processing (DSP) controllers in
	\cite{soderlind2006time,soderlind2006adaptive}. For a deeper insight to this topic we refer also to \cite{S2002,soderlind2003digital,Z1964, G1994,AS2017}.
	In those publications general, time control theoretical arguments are used to deduce stable controllers. However, given a specific scheme, there is hope to further improve the time step control. One example is given in \cite{optimizedRK}, where the authors equip RK methods with \enquote{optimized} DSP time step controllers.
	To that end, the authors sampled the domain of the hyperparameters from the controller and deduced the \enquote{optimal} parameters for the given test problems by means of minimizing the maximum, the median, or the $95\%$ percentile of the right-hand side evaluations followed by a final choice based on human interaction.

	It is the purpose of this work to exclude the human interaction by incorporating ideas from Bayesian optimization \cite{Bull2011,GSA2014, SLA2012}.
	The main benefit of this approach, besides human resources and failure, is that no expensive grid search is needed anymore, since the Bayesian optimization detects regions in the space of the hyperparameters resulting in an expensive, yet \emph{bad} time step controller.
	Additionally, the tool will recognize which of the hyperparameters have larger impact on the performance of the controller --- and which have not --- so that the search is somewhat efficient. Of course,  the term \emph{bad} must be declared to the optimization tool. Indeed, the major obstacle in this approach is to come up with an appropriate cost function.
	Our approach of constructing such a cost function is to take into account many aspects a human would consider when analyzing work precision diagrams. To validate that the proposed cost function is a reasonable choice, we will compare the search result with standard parameters from the literature and built-in solvers from \MATLAB for a wide range of problems.

	For a good comparison, we apply this novel technique to a second-order Rosenbrock method with approximate Jacobian, which is also used in the built-in \MATLAB solver \code{ode23s}. Additionally, we equip so-called modified Patankar--Runge--Kutta (MPRK) schemes \cite{KM18} for the first time with a customized time step controller. These methods are motivated by a major drawback of traditional RK methods as a subclass of general linear methods \cite{HNWII,J09}, namely the order constraint for \emph{unconditional positivity}, meaning that there exist only first-order accurate RK methods capable of producing positive approximations for $\b y^0>\b 0$ and arbitrarily large time step sizes \cite{BC78,Sandu02}. In contrast to that, MPRK schemes are nonlinear, yet linearly implicit numerical time integrators which are unconditionally positivity-preserving -- a key property of a numerical method crucial for a wide range of problems, where the method otherwise becomes inefficient or fails, see for instance \cite{STKB,KMpos,SSPMPRK2,sandu2001positive,nuesslein2021positivitypreserving} and the literature mentioned therein.

	Modified Patankar--Runge--Kutta methods stem from the Patankar trick \cite[Section 7.2-2]{Patankar1980}, which is a nonlinear weight multiplied to the RK coefficients \cite{BDM03,KM18,KM18Order3}.

	Recently, a new approach generalizing the theory of NB-series \cite{Ntrees} was presented in \cite{NSARK} providing a technique for systematically deriving conditions for arbitrary high order nonstandard additive Runge--Kutta (NSARK)  methods to which MPRK schemes belong.
	Indeed, NSARK schemes contain all Patankar-type schemes, e.g., MPDeC methods based on Deferred Correction (DeC) methods, which are conservative, unconditionally positive, and of arbitrary high order \cite{MPDeC}. Notably, these schemes were already used to preserve a positive water height when solving the shallow water equations \cite{CMOT21}. Further positivity-preserving Patankar-type schemes are discussed in \cite{gBBKS,martiradonna2020geco,SSPMPRK2,SSPMPRK3}.

	A key achievement of \cite{NSARK} was a proof that Patankar-type schemes always posses an embedded method, which is necessary for our controller. Additionally, the algorithm we present works fine with stable schemes as we do not incorporate knowledge on the stability region into the cost function. While the second-order Rosenbrock method we will use is A-stable, the stability analysis for the nonlinear Patankar-type schemes relies on a Lyapunov stability analysis recently developed \cite{IKM2122, izgin2022stability,IKMnonlin22, IOE22StabMP,gecostab}.


	Altogether, up to now the theory for order of consistency and convergence as well as stability of RK schemes was adapted to MPRK methods.
	Moreover, the schemes are also proven to be unconditionally positive and conservative. However, MPRK methods have not yet been improved with respect to their efficiency. Up to now, only standard step size controllers were used for a particular PRK method \cite{KMP21adap}. However, the construction of a tailored time step controller is still missing and constitutes one objective of the present work.

	Here, we focus on DSP controllers (reviewed in Section~\ref{sec:controller}) to investigate several MPRK schemes which we recall in Section~\ref{sec:schemes}. We emphasise that our methods rely on reinterpreting ODEs as production-destruction systems (PDSs), which is always possible. However, we demonstrate that there is no unique way to write a given ODE as a PDS. As we will encounter non-autonomous systems which are non-conservative, we present the definition of MPRK schemes in this context and discuss how this affects their order of accuracy as well as positivity. We then introduce several MPRK schemes of different orders of accuracy, along with their respective embedded methods, and review their stability properties.
	We explain our methodology in Section~\ref{sec:methodology}, where we also introduce our cost function.
	For the search of an optimal controller, we consider multiple test problems described in Appendix~\ref{sec:training},
	each of which challenges the DSP controller in a different way.
	We derive tailored controller parameters and validate the resulting methods with further test problems introduced
	in Appendix~\ref{sec:validation}, leading to the improved controllers described
	in Section~\ref{sec:PIDpara}.
	Finally, we summarize our findings in Section~\ref{sec:summary} and come to a conclusion.

\section{Numerical Schemes}
In this work we consider to classes of methods; modified Patankar--Runge--Kutta methods as they are not equipped with tailored controllers, and the \code{ode23s}-Rosenbrock method as a proof of concept that our methodology is not limited to MPRK schemes, and for a direct comparison with the built-in \code{ode23s} matlab solver. For clarity, we recall also the Rosebrock method used and highlight the differences our implementation has in contrast to \code{ode23s}.
\subsection{Rosenbrock-Type method}
Given the parameter $\gamma = (2+\sqrt2)^{-1}$, the Jacobian $\b J=\partial_{\b y}\b f(t_n,\b y^n)$, the time derivative $\b T=\partial_t \b f(t_n,\b y^n)$, and the matrix $\b W = (\b I - \gamma \dt \b J)$, the iterates of the \code{ode23s} solver \cite{shampine1997matlab} are computed by
\begin{equation*}
\begin{aligned}
\b W	\b k_1 & = \b f(t_n, \b y^n) + \gamma \dt \b T\\
\b W	\b k_2 &= \b f(t_n +0.5\dt, \b y^n + 0.5\dt \b k_1) - \b k_1 + \b W \b k_1\\
\b y^{n+1} &= \b y^n + \dt \b k_2,
\end{aligned}
\end{equation*}
which possesses the first-order embedded method given by \[\bssigma = \b y^n +\dt \b k_1.\]
Note that \code{ode23s} uses a different local error estimate using yet another evaluation of the right-hand side, see \cite{shampine1997matlab} for the details. Also, we use a central finite difference to compute the Jacobian and time derivative, \ie we use
\begin{equation*}
	\begin{aligned}
		\b T &\coloneqq \frac{\b f(t_n+ h,\b y^n) - \b f(t_n - h,\b y^n)}{2h},\quad h \coloneqq 10^{-8}\\
		 \b J \b e_i &\coloneqq \frac{\b f(t_n,\b y^n + h\b e_i) - \b f(t_n,\b y^n-h\b e_i)}{2h}, \quad i=1,\dotsc,N,
	\end{aligned}
\end{equation*}
where $\b e_i$ denoted the $i$th standard unit vector in $\R^N$.
	\subsection{Modified Patankar--Runge--Kutta Methods}\label{sec:schemes}
	We consider three modified Patankar--Runge--Kutta (MPRK) schemes, two of which are families of third-order methods and one is a family of second-order schemes.
	MPRK methods from \cite[Definition 2.1]{KM18} were first defined for the time integration of autonomous positive and conservative production-destruction systems (PDS), i.e., equations of the form
	\begin{equation} \label{eq:PDS}
		y_i'(t)=\widetilde f_i(\b y(t))=\sum_{j=1}^N (p_{ij}(\b y(t))-d_{ij}(\b y(t))),\quad \b y(0)=\b y^0>\b 0,
	\end{equation}
	where $p_{ij}(\b y(t)), d_{ij}(\b y(t)) \geq 0$ for all $\b y(t)>\b 0$ and $p_{ij}=d_{ji}$ as well as $p_{ii}=0$ for $i,j=1,\dotsc,N$. Here and in the following, vector inequalities are to be understood pointwise. Note that every right-hand side $\widetilde f_i\colon \R^N\to \R$ can be split into production and destruction terms setting \[p_{i1}(\b y)=\max\{0,\widetilde f_i(\b y)\}, \quad d_{i1}(\b y)=-\min\{0,\widetilde f_i(\b y)\},\quad p_{ij}=d_{ij}=0 \text{ for } i\neq j,\]
	so that the conservativity constrain $p_{ij}=d_{ji}$ is the only property that is not fulfilled in general. It is also worth mentioning that the additive splitting into production and destruction terms is not uniquely determined. For instance, considering
	\begin{equation*}
		\b y'=\Vec{y_2 + y_4 - y_1\\
			y_1-y_2\\
			y_1 - y_3\\
			y_3 -y_1-y_4},
	\end{equation*}
	the terms $p_{14}(\b y)=y_4$, $p_{12}(\b y)=y_2$ and $p_{43}(\b y)=y_3$ are a straightforward choice. However, both
	\[p_{21}(\b y)=y_1,\quad p_{34}(\b y)=y_1  \]
	and
	\[p_{31}(\b y)=y_1,\quad p_{24}(\b y)=y_1\]
	complete the splitting into a PDS, whereby we set $p_{mn}=0$ for the remaining production terms and $p_{ij}=d_{ji}$.

	Nevertheless, in this work, we apply MPRK schemes to general positive systems which are non-autonomous, \ie we consider a non-conservative PDS, a so-called \emph{production-destruction-rest system} (PDRS) of the form
	\begin{equation}
		y_i'(t)= f_i(\b y(t), t)=r_i(\b y(t), t) + \sum_{j=1}^N (p_{ij}(\b y(t), t)-d_{ij}(\b y(t), t)),\quad \b y(0)=\b y^0\in \R^N_{>0}, \label{eq:PDS_ODE}
	\end{equation}
	where $p_{ij}=d_{ji}$ and the rest term is also split according to
	\begin{equation}
		r_i(\b y(t), t) = r_i^p(\b y(t), t) - r_i^d(\b y(t), t)\label{eq:rp_rd}
	\end{equation}
	with $r_i^p,r_i^d\geq 0$ for $t\geq 0$ and  $\b y(t)\geq \b 0$.
	Note that $r_i^p$ and $r_i^d$ can always be constructed, for example by using the functions $\max$ and $\min$ as described above.
	The autonomous version of the PDRS \eqref{eq:PDS_ODE} was already considered in \cite{IssuesMPRK}, however, there the rest term was not split according to \eqref{eq:rp_rd}. As a consequence, the numerical solution is not guaranteed to stay positive using the Patankar-trick \cite{Patankar1980} unless $r_i>0$. This forced the authors of \cite{IssuesMPRK} to reformulate the HIRES problem (see \eqref{eq:hires} in the Appendix~\ref{sec:training}) such that $r_i>0$. In contrast, we will be using the splitting \eqref{eq:rp_rd} and can avoid the transformation.
	
	The existence, uniqueness and positivity of the solution of \eqref{eq:PDS_ODE} was discussed in \cite{FSpos}. In what follows, we are assuming that such a positive solution exists.

	To guarantee the positivity of the numerical approximation, we start with an explicit RK method and use the modified Patankar trick \cite{BDM03} on the PDS part while the rest term is only treated with the Patankar trick. Thus, $r_i^p$ will not be weighted and $r_i^d$ will be treated like a destruction term, resulting in the following new definition.
	\begin{defn}
	Given an explicit $s$-stage RK method described by a non-negative Butcher array, \ie $\b A,\b b,\b c \geq\b 0$, we define the corresponding MPRK scheme applied to the PDRS \eqref{eq:PDS_ODE}, \eqref{eq:rp_rd} by
	\begin{equation}\label{eq:MPRK_PDRS}
		\begin{aligned}
			y_i^{(k)}=& y_i^n + \dt\sum_{\nu=1}^{k-1}a_{k\nu}\left( r_i^p(\b y^{(\nu)}, t_n + c_\nu\dt) + \sum_{j=1}^N  p_{ij}(\b y^{(\nu)}, t_n + c_\nu\dt)\frac{y_j^{(k)}}{\pi_j^{(k)}}\right. \\
			& \left. - \left( r_i^d(\b y^{(\nu)}, t_n + c_\nu\dt)+\sum_{j=1}^Nd_{ij}(\b y^{(\nu)}, t_n + c_\nu\dt) \right)\frac{y_i^{(k)}}{\pi_i^{(k)}}\right),\quad k=1,\dotsc,s,\\
			y_i^{n+1}=& y_i^n + \dt\sum_{k=1}^{s}b_{k}\left( r_i^p(\b y^{(k)}, t_n + c_k\dt) + \sum_{j=1}^N  p_{ij}(\b y^{(k)}, t_n + c_k\dt)\frac{y_j^{n+1}}{\sigma_j}\right. \\
			& \left. - \left(r_i^d(\b y^{(k)}, t_n + c_k\dt) + \sum_{j=1}^Nd_{ij}(\b y^{(k)}, t_n + c_k\dt)\right)\frac{y_i^{n+1}}{\sigma_i}\right)
		\end{aligned}
	\end{equation}
	for $i=1,\dotsc,N$, where $\pi_i^{(k)},\sigma_i$ are the so-called \emph{Patankar-weight denominators} (PWDs), which are required to be positive for any $\dt\geq 0$, and independent of the corresponding numerators $y_i^{(k)}$ and $y_i^{n+1}$, respectively.
	\end{defn}
As described in \cite{KM18}, MPRK schemes can be written in matrix notation, which in the case of PDRS is given below.
	\begin{rem}
		In matrix notation, \eqref{eq:MPRK_PDRS} can be rewritten as
		\begin{equation}\label{eq:MPRK_PDRS_matrix}
			\begin{aligned}
				\b M^{(k)}\byk&= \b y^n+\dt \sum_{\nu=1}^{k-1}a_{k\nu} \b r^p(\b y^{(\nu)}, t_n + c_\nu\dt),\quad k=1,\dotsc,s, \\
				\b M\b y^{n+1}&= \b y^n + \dt\sum_{k=1}^{s}b_{k} \b r^p(\b y^{(k)}, t_n + c_k\dt),
			\end{aligned}
		\end{equation}
		where $\b r^p=(r_1^p,\dotsc, r_N^p)^T$, $\b M^{(k)}=(m^{(k)}_{ij})_{1\leq i,j\leq N}$, and $\b M=(m_{ij})_{1\leq i,j\leq N}$ with
		\begin{equation*}
			\begin{aligned}
				m^{(k)}_{ii}&= 1+ \dt \sum_{\nu=1}^{k-1}a_{k\nu}\left(r_i^d(\b y^{(\nu)}, t_n + c_\nu\dt) + \sum_{j=1}^Nd_{ij}(\b y^{(\nu)}, t_n + c_\nu\dt) \right)\frac{1}{\pi_i^{(k)}}, \\
				m^{(k)}_{ij}&= -\dt \sum_{\nu=1}^{k-1}a_{k\nu}p_{ij}(\b y^{(\nu)}, t_n + c_\nu\dt)\frac{1}{\pi_j^{(k)}}, \quad i\neq j
			\end{aligned}
		\end{equation*}
		as well as
		\begin{equation*}
			\begin{aligned}
				m_{ii}&= 1+ \dt \sum_{k=1}^{s}b_k\left(r_i^d(\byk, t_n + c_k\dt)+\sum_{j=1}^N d_{ij}(\byk, t_n + c_k\dt) \right)\frac{1}{\sigma_i}, \\
				m_{ij}&= -\dt \sum_{k=1}^{s}b_kp_{ij}(\byk, t_n + c_k\dt)\frac{1}{\sigma_j}, \quad i\neq j.
			\end{aligned}
		\end{equation*}
		We want to note that this scheme always produces positive approximations if $\b y^0>\b 0$ as $\b M^T$ is still an M-matrix following the same lines as in \cite[Lemma 2.8]{KM18} and exploiting $\b r^d\geq \b 0$. However, if it is known that the analytic solution is not positive due to the existence of the rest term $\b r$, then one may consider choosing the splitting $\b r^d=\b 0$ and $\b r^p = \b r$ in the MPRK scheme \eqref{eq:MPRK_PDRS}. This means that we drop the non-negativity constraint on $\b r^p$ and do not perform the Patankar trick on the rest term. As a result, the right-hand sides in \eqref{eq:MPRK_PDRS_matrix} are allowed to be negative, and thus, the stage vectors and iterates of the MPRK scheme are not forced to stay positive anymore.
	\end{rem}

	Next, we want to explain in what sense the given definition of MPRK schemes generalizes the existing ones from \cite{KM18,IssuesMPRK}.
	To that end, we assume that $\pi^{(k)}_i, \sigma_i$ only depend on the $i$th component of the stages, \ie \begin{equation} \label{eq:pik_sigma}
		\pi_i^{(k)}=\pi_i^{(k)}(y_i^n,y_i^{(1)},\dotsc, y_i^{(k-1)}) \qta \sigma_i=\sigma_i(y_i^n, y_i^{(1)},\dotsc, y_i^{(s)}),
	\end{equation}
	which was already assumed in \cite{AGKM_Oliver} and includes the PWDs developed so far, see \cite{KM18,KM18Order3, SSPMPRK2, SSPMPRK3, IssuesMPRK}.
	\begin{prop}
		Assume that the PWDs $\bssigma$ and $\bspi^{(k)}$ satisfy the assumption \eqref{eq:pik_sigma} for $k=1,\dotsc,N$, and that $\sum_{\nu=1}^sa_{k\nu}=c_k$ as well as $\sum_{k=1}^sb_k=1$.
		Then the MPRK scheme \eqref{eq:MPRK_PDRS} applied to \eqref{eq:PDS_ODE}, \eqref{eq:rp_rd} produces the same approximations as when applied to the corresponding autonomous system $\b Y'(t)=\b F(\b Y(t))$ using
		\[\b Y(t)=(\b y(t),t),\quad \b F(\b Y(t))=\begin{pmatrix}
			\b f(\b Y(t))\\ 1
		\end{pmatrix} \]
		together with the natural choice of writing the right-hand side as a PDRS, \ie
		\begin{equation}\label{eq:Ftrans_to_autonom}
			\b F(\b Y(t)) = \begin{pmatrix}
				r_1(\b Y(t)) + \sum_{j=1}^N (p_{1j}(\b Y(t))-d_{1j}(\b Y(t)))\\
				\vdots\\
				r_N(\b Y(t)) + \sum_{j=1}^N (p_{Nj}(\b Y(t))-d_{Nj}(\b Y(t)))\\
				r_{N+1}\end{pmatrix}, \quad r_{N+1}=r_{N+1}^p=1,
		\end{equation}
		which means that $p_{j,N+1}=d_{j,N+1}=0$ and $p_{N+1,j}=d_{N+1,j}=0$ for $j=1,\dotsc, N+1$.
	\end{prop}
	\begin{proof}
		The MPRK scheme \eqref{eq:MPRK_PDRS} applied to $\b Y'(t)=\b F(\b Y(t))$ with $\b F$ from \eqref{eq:Ftrans_to_autonom} reads
		\begin{equation}\label{eq:MPRK_PDRS_autonom}
			\begin{aligned}
				Y_i^{(k)}=& Y_i^n + \dt\sum_{\nu=1}^{k-1}a_{k\nu}\left( r_i^p(\b Y^{(\nu)}) + \sum_{j=1}^N  p_{ij}(\b Y^{(\nu)})\frac{Y_j^{(k)}}{\pi_j^{(k)}(Y_j^n,Y_j^{(1)},\dotsc, Y_j^{(k-1)})}\right. \\
				& \left. - \left(r_i^d(\b Y^{(\nu)}) + \sum_{j=1}^Nd_{ij}(\b Y^{(\nu)}) \right)\frac{Y_i^{(k)}}{\pi_i^{(k)}(Y_i^n,Y_i^{(1)},\dotsc, Y_i^{(k-1)})}\right),\quad k=1,\dotsc,s,\\
				Y_i^{n+1}=& Y_i^n + \dt\sum_{k=1}^{s}b_{k}\left( r_i^p(\b Y^{(k)}) + \sum_{j=1}^N  p_{ij}(\b Y^{(k)})\frac{Y_j^{n+1}}{\sigma_j(Y_j^n,Y_j^{(1)},\dotsc, Y_j^{(s)})}\right. \\
				& \left. - \left(r_i^d(\b Y^{(k)}) + \sum_{j=1}^Nd_{ij}(\b Y^{(k)})\right)\frac{Y_i^{n+1}}{\sigma_i(Y_i^n,Y_i^{(1)},\dotsc, Y_i^{(s)})}\right).
			\end{aligned}
		\end{equation}
		For $i=N+1$, this reduces to
		\begin{equation} \label{eq:tk_tn+1}
			\begin{aligned}
				t_{(k)}&=t_n + \dt \sum_{\nu=1}^{k-1}a_{k\nu}=t_n+c_k\dt,\\
				t_{n+1} &= t_n +\dt \sum_{k=1}^sb_k=t_n+\dt.
			\end{aligned}
		\end{equation}
		Furthermore, for $i\leq N$, we know from the assumption \eqref{eq:pik_sigma} that $\pi_i^{(k)}(Y_i^n,Y_i^{(1)},\dotsc, Y_i^{(k-1)})=\pi_i^{(k)}(y_i^n,y_i^{(1)},\dotsc, y_i^{(k-1)})$ and similarly for $\sigma_i$. Thus, we end up with
		\begin{equation*}
			\begin{aligned}
				y_i^{(k)}=& y_i^n + \dt\sum_{\nu=1}^{k-1}a_{k\nu}\left( r_i^p(\b y^{(\nu)},t_{(\nu)}) + \sum_{j=1}^N  p_{ij}(\b y^{(\nu)}, t_{(\nu)})\frac{y_j^{(k)}}{\pi_j^{(k)}}\right. \\
				& \left. - \left(r_i^d(\b y^{(\nu)}, t_{(\nu)})+\sum_{j=1}^Nd_{ij}(\b y^{(\nu)}, t_{(\nu)})  \right)\frac{y_i^{(k)}}{\pi_i^{(k)}}\right),\quad k=1,\dotsc,s,\\
				y_i^{n+1}=& y_i^n + \dt\sum_{k=1}^{s}b_{k}\left( r_i^p(\b y^{(k)}, t_{(k)}) + \sum_{j=1}^N  p_{ij}(\b y^{(k)}, t_{(k)})\frac{y_j^{n+1}}{\sigma_j}\right. \\
				& \left. - \left(r_i^d(\b y^{(k)}, t_{(k)}) + \sum_{j=1}^Nd_{ij}(\b y^{(k)}, t_{(k)})\right)\frac{y_i^{n+1}}{\sigma_i}\right).
			\end{aligned}
		\end{equation*}
		Substituting \eqref{eq:tk_tn+1} into these equations, the proof is finished by noting that the resulting linear systems always possess a unique solution.
	\end{proof}
	\begin{rem}
		To see that the order conditions derived in \cite{NSARK, KM18, KM18Order3} remain valid also for the non-autonomous PDRS case, we rewrite the MPRK scheme \eqref{eq:MPRK_PDRS} applied to the transformed system \eqref{eq:Ftrans_to_autonom} as a nonstandard additive Runge--Kutta (NSARK) method, for which the order conditions are already known, see \cite{NSARK}. To that end, we split \eqref{eq:Ftrans_to_autonom} according to
		\begin{align*}
			\b F(\b Y(t))= \sum_{j=1}^{N+1}  \Fj(\b Y(t)),
		\end{align*}
		using  $\b F^{[N+1]}(\b Y(t))=(r_1^p(\b Y(t)),\dotsc,r_N^p(\b Y(t)),r_{N+1})^T\in \R^{N+1}$ and \[F^{[j]}_i(\b Y(t))=\begin{cases}p_{ij}(\b Y(t)), &i\neq j, \\-\left(r_i^d(\b Y(t)) + \sum_{j=1}^Nd_{ij}(\b Y(t))\right), &i=j, \\

			0, & i=N+1\end{cases}\]
		for $j\leq N$.

		With that, we see from \eqref{eq:MPRK_PDRS_autonom} that the MPRK scheme takes the form of an NSARK method
		\begin{equation*}
			\begin{aligned}
				\b Y^{(k)} & = \b Y^n + \dt \sum_{\nu=1}^{k-1}  \sum_{j=1}^{N+1} a^{[j]}_{k\nu}(\b Y^n,\dt)  \Fj(	\b Y^{(\nu)}), \quad k=1,\dotsc,s,\\
				\b Y^{n+1} & = \b Y^n + \dt \sum_{k=1}^s \sum_{j=1}^{N+1} b^{[j]}_k(\b Y^n,\dt) \Fj(	\b Y^{(k)})
			\end{aligned}
		\end{equation*}
		with the solution-dependent RK coefficients
		\begin{equation*}
			a^{[j]}_{k\nu}(\b Y^n,\dt)=\begin{cases}
				a_{k\nu}\frac{Y_j^{(k)}}{\pi_j^{(k)}}, &j\leq N,\\
				a_{k\nu}, & j=N+1,
			\end{cases}\quad  b^{[j]}_k(\b Y^n,\dt)=\begin{cases}
				b_{k}\frac{Y_j^{n+1}}{\sigma_j}, &j\leq N,\\
				b_{k}, & j=N+1.
			\end{cases}
		\end{equation*}
		Altogether, it thus immediately follows from \cite{NSARK} that the order conditions for the generalized MPRK schemes in the case of non-autonomous PDRS coincide with those for autonomous PDS.
	\end{rem}
We want to point out that the above results are new and enable the application of MPRK schemes to general PDRS, while still maintaining their accuracy, positivity. Additionally, the conservative part of the PDRS will be also conservative on the fully discrete setting. In the following, we recall known MPRK methods that will be equipped with a customized time step controller.

	\subsubsection{Second-Order MPRK Schemes}
	The explicit two-stage RK method based on the Butcher array
	\begin{equation*}
		\begin{aligned}
			\def\arraystretch{1.2}
			\begin{array}{c|ccc}
				0 &  & \\
				\alpha & \alpha &  \\
				\hline
				& 1-\frac{1}{2\alpha} &\frac{1}{2\alpha}
			\end{array}
		\end{aligned}
	\end{equation*}
	is second-order accurate \cite[Section~320]{B16}. Moreover, the entries of the array are non-negative for $\alpha\geq~\!\! \frac12$. With that as a starting point, the authors from \cite{KM18} derived a one-parameter family of second-order accurate MPRK schemes, denoted by MPRK22($\alpha$), using the PWDs $\pi_i^{(2)}=y_i^n$ and $\sigma_i= (y_i^{(2)})^{\frac{1}{\alpha}}(y_i^n)^{1-\frac{1}{\alpha}}$ for $i=1,\dotsc,N$. For simplicity, we present the resulting MPRK22($\alpha$) scheme for the case of a conservative and autonomous PDS, that is
	\begin{equation}\label{eq:MPRK22b}
		\begin{aligned}
			y_i^{(1)} =&\, y_i^n,\\
			y_i^{(2)} =&\, y_i^n + \alpha\Delta t\sum_{j=1}^N\left(p_{ij}(\b y^{(1)})\frac{y_j^{(2)}}{y_j^n}-d_{ij}(\b y^{(1)})\frac{y_i^{(2)}}{y_i^n}\right),\\
			y_i^{n+1} =&\, y_i^n + \Delta t\sum_{j=1}^N\left( \Biggl(\biggl(1-\frac1{2\alpha}\biggr) p_{ij}(\b y^{(1)})+\frac1{2\alpha} p_{ij}(\b y^{(2)})\Biggr)\frac{y_j^{n+1}}{(y_j^{(2)})^{\frac{1}{\alpha}}(y_j^n)^{1-\frac{1}{\alpha}}}\right.\\
			& \left. - \Biggl(\biggl(1-\frac1{2\alpha}\biggr) d_{ij}(\b y^{(1)})+ \frac1{2\alpha} d_{ij}(\b y^{(2)})\Biggr)\frac{y_i^{n+1}}{(y_i^{(2)})^{\frac{1}{\alpha}}(y_i^n)^{1-\frac{1}{\alpha}}}\right)
		\end{aligned}
	\end{equation}
	for $i=1,\dots,N$ with $\alpha\geq\frac12$.

	\subsubsection{Third-Order MPRK Schemes}
\begin{figure}[htbp]
		\centering
		\includegraphics[width=.7\textwidth]{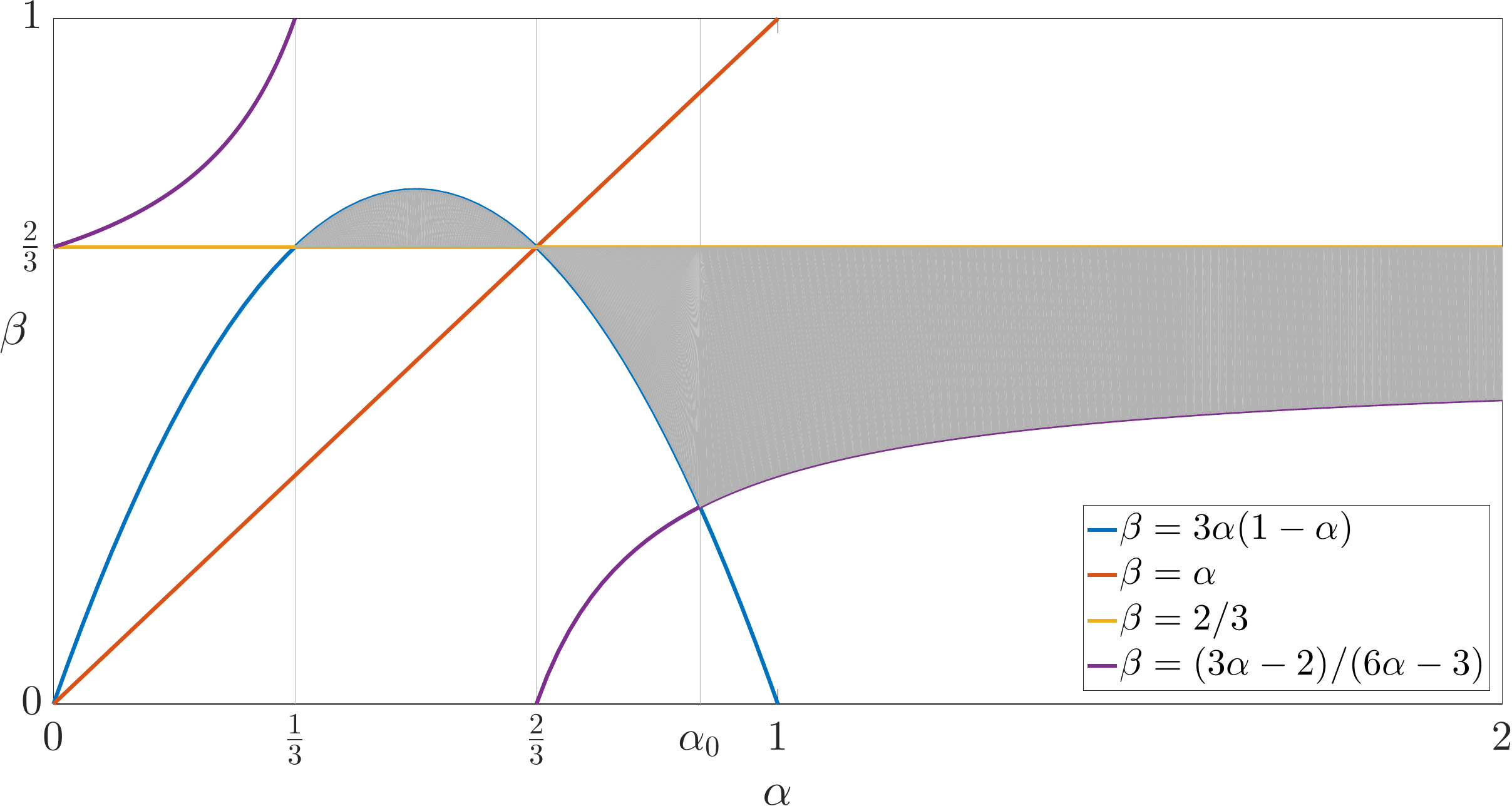}
		\caption{The gray area represents all $(\alpha,\beta)$ pairs which fulfill the conditions \eqref{eq:cond:alpha,beta}, \ie for which the Butcher tableau \eqref{array:MPRK43alphabeta} is non-negative \cite{KM18Order3}.}\label{fig:RK3_pos}
	\end{figure}

	Assuming a non-negative Butcher tableau from an explicit 3-stage RK method, third-order MPRK schemes have been constructed in \cite{KM18Order3} for conservative and autonomous PDS using the PWDs
	\begin{equation}
		\begin{aligned}
			\pi_i^{(2)} =\, &y_i^n,\\
			\pi_i^{(3)}  =\, &(y_i^{(2)})^{\frac{1}{p}}(y_i^n)^{1-\frac{1}{p}},\quad p=3a_{21}(a_{31}+a_{32})b_3,\\
			\sigma_i=\, &y_i^n + \dt \sum_{j=1}^N\left( \left( \beta_1p_{ij}(\b y^n)+\beta_2p_{ij}(\b y^{(2)})\right)\frac{\sigma_j}{(y_j^{(2)})^{\frac{1}{a_{21}}}(y_j^n)^{1-\frac{1}{a_{21}}}}\right. \\
			& - \left.\left( \beta_1d_{ij}(\b y^n)+\beta_2d_{ij}(\b y^{(2)})\right)\frac{\sigma_i}{(y_i^{(2)})^{\frac{1}{a_{21}}}(y_i^n)^{1-\frac{1}{a_{21}}}}\right)
		\end{aligned}
	\end{equation}
	for $i=1,\dotsc,N$, $\beta_1=1-\beta_2$ and $\beta_2=\frac{1}{2a_{21}}$. Note, that solving another system of linear equations is necessary to calculate $\bm \sigma=(\sigma_1,\dotsc,\sigma_N)$. Hence, the resulting MPRK scheme may be based on 3-stage RK methods but can be viewed as 4-stage schemes, whereby we note that $\bm \sigma$ can be computed simultaneously with $\b y^{(3)}$.  We also point out that there are no additional right-hand side evaluations required for computing $\bm\sigma$. The final scheme for conservative and autonomous PDS takes the form
	\begin{subequations}
		\label{eq:MPRK43-family}
		\begin{align}
			y^{(1)}_i &= y^n_i,\label{eq:MPRK43a}\\
			y^{(2)}_i &= y^n_i
			+  a_{21} \dt \sum_{j=1}^N \left(
			p_{ij}\bigl( \b y^n \bigr) \frac{y^{(2)}_j}{y^n_j}
			- d_{ij}\bigl( \b y^n \bigr) \frac{y^{(2)}_i}{y^n_i}
			\right),
			\\
			y^{(3)}_i &= y^n_i
			+\Delta t \sum_{j=1}^N
			\Biggl(\left(a_{31} p_{ij}\bigl(\b y^n\bigr)+ a_{32} p_{ij} \bigl(\b y^{(2)}\bigr) \right)  \frac{  y_j^{(3)}
			}{\bigl(y_j^{(2)}\bigr)^{\frac1p } \bigl(y_j^n\bigr)^{1-\frac1p} }
			\nonumber\\
			& \qquad\qquad\qquad
			-\left(a_{31} d_{ij}\bigl(\b y^n\bigr)+ a_{32} d_{ij} \bigl(\b y^{(2)}\bigr) \right)  \frac{  y_i^{(3)}
			}{\bigl(y_i^{(2)}\bigr)^{\frac1p } \bigl(y_i^n\bigr)^{1-\frac1p} }
			\Biggr),
			\\
			\sigma_i &= y_i^n + \Delta t \sum_{j=1}^N
			\Biggl(\left( \beta_1 p_{ij} \bigl( \b y^n\bigr) +\beta_2 p_{ij} \bigl(\b y^{(2)}\bigr)  \right) \frac{\sigma_j}{\bigl(y_j^{(2)} \bigr)^{\frac1q}
				\bigl(y_j^n\bigr)^{1-\frac1q}}
			\nonumber	\\ & \qquad-
			\left( \beta_1 d_{ij} \bigl( \b y^n\bigr) +\beta_2 d_{ij} \bigl(\b y^{(2)}\bigr)  \right) \frac{\sigma_i}{\bigl(y_i^{(2)} \bigr)^{\frac1q}
				\bigl(y_i^n\bigr)^{1-\frac1q}}\Biggr),\label{eq:MPRK43sigma}
			\\
			y^{n+1}_i &= y^n_i
			+ \dt \sum_{j=1}^N \Biggl(
			\left( b_1 p_{ij}\bigl( \b y^n \bigr) +b_2p_{ij}\bigl( \b y^{(2)} \bigr)
			+ b_3 p_{ij}\bigl( \b y^{(3)} \bigr)
			\right) \frac{y^{n+1}_j}{\sigma_j}
			\nonumber\\&\qquad\qquad\qquad
			- \left( b_1 d_{ij}\bigl( \b y^n \bigr) +b_2d_{ij}\bigl( \b y^{(2)} \bigr)
			+ b_3 d_{ij}\bigl( \b y^{(3)} \bigr)
			\right) \frac{y^{n+1}_i}{\sigma_i}
			\Biggr),
		\end{align}
	\end{subequations}
	where $p=3 a_{21}\left(a_{31}+a_{32} \right)b_3,\; q=a_{21},\;\beta_2=\frac{1}{2a_{21}}$ and $\beta_1= 1-\beta_2$.
	\paragraph{MPRK43($\alpha, \beta$)}\label{subsec:MPRK43alphabeta}
	All entries of the Butcher array
	\begin{equation}\label{array:MPRK43alphabeta}
		\begin{aligned}
			\def\arraystretch{1.2}
			\begin{array}{c|ccc}
				0 &  & & \\
				\alpha & \alpha & & \\
				\beta & \frac{3\alpha\beta (1-\alpha)-\beta^2}{\alpha(2-3\alpha)}& \frac{\beta (\beta-\alpha)}{\alpha(2-3\alpha)}& \\
				\hline
				& 1+\frac{2-3(\alpha+\beta)}{6 \alpha \beta } &\frac{3 \beta-2}{6\alpha (\beta-\alpha)} & \frac{2-3\alpha}{6\beta(\beta-\alpha)}
			\end{array}
		\end{aligned}
	\end{equation}
	with
	\begin{equation}\label{eq:cond:alpha,beta}
		\begin{cases}
			2/3 \leq \beta \leq 3\alpha(1-\alpha)\\
			3\alpha(1-\alpha)\leq\beta \leq 2/3 \\
			\tfrac{3\alpha-2}{6\alpha-3}\leq \beta \leq 2/3
		\end{cases}
		\text{ for }
		\begin{cases}
			1/3 \leq \alpha<\frac23,\\
			2/3 < \alpha<\alpha_0,\\
			\alpha>\alpha_0,
		\end{cases}
	\end{equation}
	and $\alpha_0\approx 0.89255$  are non-negative \cite[Lemma 6]{KM18Order3}. Figure~\ref{fig:RK3_pos} illustrates the feasible domain. Moreover, the corresponding RK method is third-order accurate \cite{RR01}.

	The resulting MPRK43$(\alpha,\beta)$ method is determined by \eqref{eq:MPRK43-family} using
	\eqref{array:MPRK43alphabeta} and
	\begin{equation*}
		\begin{aligned}
			p&=3a_{21}(a_{31}+a_{32})b_3=\alpha\frac{2-3\alpha}{2(\beta-\alpha)},&
			q&=a_{21}=\alpha,\\
			\beta_2&=\frac{1}{2a_{21}}= \frac{1}{2\alpha},& \beta_1&=1-\beta_2=1-\frac{1}{2\alpha}.
		\end{aligned}
	\end{equation*}

	\paragraph{MPRK43($\gamma$)}\label{subsec:MPRK43gamma}
	It was also proven in \cite[Lemma 6]{KM18Order3} that all entries of the tableau
	\begin{equation}\label{array:MPRK43gamma}
		\begin{aligned}
			\def\arraystretch{1.2}
			\begin{array}{c|ccc}
				0 &  & & \\
				\frac23 & \frac23 & & \\
				\frac23 & \frac23-\frac{1}{4\gamma}& \frac{1}{4\gamma}& \\
				\hline
				& \frac14 &\frac34-\gamma & \gamma
			\end{array}
		\end{aligned}
	\end{equation}
	are non-negative for $\frac38\leq \gamma \leq \frac34.$ Furthermore, in \cite{RR01} it was proven that the resulting RK method is third-order accurate. The corresponding third-order MPRK scheme is denoted by MPRK43$(\gamma)$ and can be obtained from \eqref{eq:MPRK43-family} by substituting \eqref{array:MPRK43gamma} and
	\begin{equation*}
		\begin{aligned}
			p&=3a_{21}(a_{31}+a_{32})b_3=\frac43\gamma,&
			q&=a_{21}=\frac23,\\
			\beta_2&=\frac{1}{2a_{21}}=\frac34, &\beta_1&=1-\beta_2=\frac14.
		\end{aligned}
	\end{equation*}
	\subsubsection{Embedded Methods of MPRK Schemes}\label{sec:embedded}
	It was proven in \cite{KM18,KM18Order3} and generalized in \cite[Lemma 5.4]{NSARK} that if the MPRK scheme is of order $k$, then the embedded method returning $\bm \sigma$ is of order $k-1$. In particular, in the case of MPRK22($\alpha$) the embedded method is first order and reads
	\begin{equation*}
		\begin{aligned}
			y_i^{(1)} = &\, y_i^n,\\
			y_i^{(2)} = &\, y_i^n + \alpha\Delta t\sum_{j=1}^N\left(p_{ij}(\b y^{(1)})\frac{y_j^{(2)}}{y_j^n}-d_{ij}(\b y^{(1)})\frac{y_i^{(2)}}{y_i^n}\right),\\
			\sigma_i = &\, (y_i^{(2)})^{\frac{1}{\alpha}}(y_i^n)^{1-\frac{1}{\alpha}}
		\end{aligned}
	\end{equation*}
	for $i=1,\dotsc,N$. Similarly, the embedded second-order scheme for the two MPRK43 families is given by \eqref{eq:MPRK43a}--\eqref{eq:MPRK43sigma} with the respective parameters specified in the Subsections~\ref{subsec:MPRK43alphabeta} and \ref{subsec:MPRK43gamma}.

	\subsubsection{Preferable Members of MPRK Families}\label{sec:ParasMPRK}
	Since several families of MPRK methods exist \cite{KM18,KM18Order3}, a first challenge is to determine which member to choose. An intuitive way to disqualify certain members of the family is based on a stability investigation. However, as MPRK schemes do not belong to the class of general linear methods, a new approach was used in \cite{IKM2122,izgin2022stability} to investigate their stability properties. The resulting theory is based on the center manifold theorem for maps \cite{carr1982,SH98} and was applied to second-order MPRK22($\alpha$) \cite{izgin2022stability} and to third-order MPRK methods MPRK43($\alpha,\beta$) and MPRK43($\gamma$) in \cite{IOE22StabMP}. To that end, the linear problem
	\begin{equation}\label{eq:Stabtestprob}
		\b y'(t)=\b \Lambda \b y(t),\quad \b y(0)=\b y^0\in \R^N_{>0}
	\end{equation}
	was considered, where $\b \Lambda=(\lambda_{ij})_{i,j=1,\dotsc,N} \in \R^{N\times N}$ is a Metzler matrix, \ie $\lambda_{ij}\geq 0$ for $i\neq j$, which additionally satisfies $\sum_{i=1}^N\lambda_{ij}=0$ for $j=1,\dotsc,N$. These two properties guarantee the positivity and conservativity of \eqref{eq:Stabtestprob}, see \cite{izgin2022stability}. In particular, the aim of the analysis was to derive restrictions for the time step size guaranteeing that Lyapunov stable steady states of \eqref{eq:Stabtestprob} are Lyapunov stable fixed points of the numerical scheme, see \cite{IKM2122} for the details. Moreover, in the case of stable fixed points, the iterates of the MPRK method are proved to locally converge towards the correct steady state solution. The stability properties as well as the local convergence could be observed in numerical experiments \cite{izgin2022stability,IOE22StabMP}. Indeed, it turned out that MPRK22$(\alpha$) and MPRK43$(\gamma$) are stable in this sense for all parameter choices \cite{izgin2022stability,IOE22StabMP}. Still, the work \cite{IssuesMPRK} favors $\alpha=1$ for the MPRK22$(\alpha$) scheme for reasons of oscillatory behavior and the existence of spurious steady steady states. In the case of MPRK43$(\gamma$) we restrict to $\gamma=0.563$, since the corresponding scheme has the largest $\dt$ bound while fulfilling the necessary condition for avoiding oscillatory behavior, see  \cite{ITOE22}.

Also, applying the stability analysis from \cite{IOE22StabMP} to the pairs $(\alpha,\beta)$ used so far \cite{KM18Order3}, the pair $(\alpha,\beta)=(0.5,0.75)$ is associated with the largest stability domain.

	As a result of that numerical analysis, we will consider MPRK43($\alpha,\beta$) with $\alpha=0.5$ and $\beta =0.75$ in the following.

	Altogether, we will consider in this work the schemes MPRK22$(1)$, MPRK43$(0.5,0.75)$, and MPRK43$(0.563)$.
	\section{The Controller}\label{sec:controller}
	Given an absolute tolerance \code{atol} and a relative tolerance \code{rtol}, we define the weighted error estimate
	\begin{equation}\label{eq:PIDepsomega}
		w_{n+1} = \left( \frac1N\sum_{i=1}^N\left(\frac{y_i^{n+1}-\sigma_i}{\code{atol} + \code{rtol}\max\{ \lvert y_i^{n+1}\rvert,\lvert\sigma_i\rvert\}}\right)^2\right)^{\frac{1}{2}},
	\end{equation}
	where $\b y^{n+1}=(y_1^{n+1},\dotsc,y_N^{n+1})^T$ is the numerical approximation at time $t_{n+1}$ and $\bm\sigma$ the corresponding output of the embedded $(k-1)$th order method.
	Next, we set
	\begin{equation}
		\epsilon_{n+1} = \frac{1}{\max\{\code{eps}, w_{n+1}\}},
	\end{equation}
	where \code{eps} denotes machine precision.
	The \emph{proportional-integral-derivative} (PID) controller with free parameters \[(\beta_1,\beta_2,\beta_3)\in [0.1,1]\times [-0.4, -0.05]\times [0, 0.1]\] is then given by
	\begin{equation*}
		\dt_{n+1} =\epsilon_{n+1}^\frac{\beta_1}{k}\epsilon_n^\frac{\beta_2}{k}\epsilon_{n-1}^\frac{\beta_3}{k} \dt_n,
	\end{equation*}
	with $\epsilon_0=\epsilon_{-1}=1$. Here, $k$ denotes the order of the numerical method.
	Furthermore, for solving stiff problems, we extend the controller using
	ideas from digital signal processing \cite{soderlind2003digital}. Thus, we use the additional factor $\left(\frac{\dt_n}{\dt_{n-1}}\right)^{-\alpha_2}$ with $\alpha_2\in[1/6,1/2]$, so that
	\begin{equation*}
		\dt_{n+1} =\epsilon_{n+1}^\frac{\beta_1}{k}\epsilon_n^\frac{\beta_2}{k}\epsilon_{n-1}^\frac{\beta_3}{k}\left(\frac{\dt_n}{\dt_{n-1}}\right)^{-\alpha_2} \dt_n.
	\end{equation*}  Also, it is common to limit the new time step by a bounded function. Altogether, we use the DSP controller
	\begin{equation}\label{eq:PIDlimit}
		\dt_{n+1} =\left(1 + \kappa_2\arctan\left(\frac{\epsilon_{n+1}^\frac{\beta_1}{k}\epsilon_n^\frac{\beta_2}{k}\epsilon_{n-1}^\frac{\beta_3}{k}\left(\frac{\dt_n}{\dt_{n-1}}\right)^{-\alpha_2}  - 1}{\kappa_2}\right)\right) \dt_n
	\end{equation}
	with $\kappa_2\in \{1,2\}$ as suggested in \cite{soderlind2006adaptive}, \ie
	\begin{equation*}
		(\beta_1,\beta_2,\beta_3,\alpha_2,\kappa_2)\in [0.1,1]\times [-0.4, -0.05]\times [0, 0.1]\times [1/6,1/2]\times\{1,2\}.
	\end{equation*}
	We also note that if the coefficient of $\dt_n$ is smaller than a safety value of $s_f\coloneqq0.81$, we reject the current step since it led to a small coefficient of $\dt_n$ and use $\dt_{n+1}$ to recalculate it.

\begin{figure}[!htbp]
		\centering
		\begin{scaletikzpicturetowidth}{0.5\textwidth}
			\begin{tikzpicture}[scale=\tikzscale]
				\begin{axis}[xmode=log,
					ytick distance=1,
					enlarge x limits=false,ymin=0,ymax=8.4,xmax=1500,xmin=0.001,
					minor tick num=1, axis lines = left,xlabel=$x$,ylabel=$y$, ylabel near ticks,
					ylabel style={rotate=-90},
					every axis plot/.append style={thick}]
					\addplot [scale=1, domain=0.001:1000, smooth, variable=\x, blue, solid, mylabel=at 0.9 below right with {$L_1$}] plot ({\x}, {1+rad(atan(\x-1))});
					\addplot [scale=1, domain=0.001:1000, smooth, variable=\x, purple, dashed, mylabel=at 0.95 below with {$L_2$}]  plot ({\x}, {1+2*rad(atan((\x-1)/2))});
					\addplot [scale=1, domain=0.001:1000, smooth, variable=\x, red, dashed, mylabel=at 0.95 below with {$L_3$}]  plot ({\x}, {1+3*rad(atan((\x-1)/3))});
					\addplot [scale=1, domain=0.001:1000, smooth, variable=\x, orange, dashed, mylabel=at 0.95 below with {$L_4$}]  plot ({\x}, {1+4*rad(atan((\x-1)/4))});
					\addplot [scale=1, domain=0.001:1000, variable=\x, black, dotted, mylabel=at 0.95 above with {$s_f$}]  plot ({\x}, {0.81});
				\end{axis}
			\end{tikzpicture}
		\end{scaletikzpicturetowidth}
		\caption{Plots of $s_f=0.81$ and $L_{\kappa_2}(x)=1+\kappa_2\arctan\left(\frac{x - 1}{\kappa_2}\right)$ for $\kappa_2\in \{1,2,3,4\}$ with logarithmic $x$-axis. }\label{fig:PIDlim}
	\end{figure}
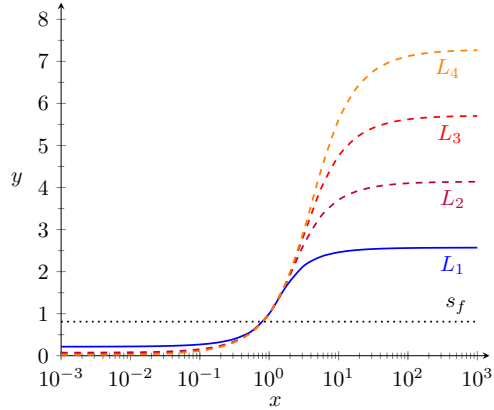

	As we want to test the reliability of the cost function, we rather consider the larger domain
	\begin{equation}\label{eq:PIDparas2}
	 D\coloneqq[-5,5]\times [-3, 3]\times [-2, 2]\times [-3,3]\times\{1,2,3,4\}
\end{equation}
including further controllers from \cite{soderlind2003digital,Z1964, Watts1984,G1994,AS2017}.
	To give an insight on the effect of the limiter, we plot the function \[L_{\kappa_2}(x)=1+\kappa_2\arctan\left(\frac{x - 1}{\kappa_2}\right)\] for the mentioned values of $\kappa_2$ in Figure~\ref{fig:PIDlim}.

	As it can be observed from Figure~\ref{fig:PIDlim}, $L_{\kappa_2}$ is steeper and allows larger changes in the step size for increasingly values of $\kappa_2\in\{1,2,3,4\}$ while maintaining approximately the same threshold for rejecting a step, meaning the intersection points of $L_{\kappa_2}$ with $s_f=0.81$ all lie in the interval $[0.8076,0.81]$.

	\subsection{Criteria for Good Controllers}

	Designing good step size controllers is difficult in general.
	If we consider a fixed combination of main method and embedded method,
	we could sample the space of all possible controller parameters. Then,
	for each problem and each tolerance, we could find an ``optimal''
	controller resulting in the least error with the lowest amount of
	computational work. However, we do not want to make the controller
	dependent on the ODE to be integrated or the tolerance chosen by the
	user. Thus, we need to find controllers that perform well for a range
	of problems and tolerances.

	The classical (deadbeat) I controller \cite[Section~II.4]{HNW1993}
	is derived for the asymptotic regime of small time step sizes and
	can be argued to be optimal there. However, practical applications
	require controllers that also work well for bigger tolerances.
	For explicit Runge--Kutta methods and mildly stiff problems, step size
	control stability is important \cite{hall1988analysis,higham1990embedded}
	and led to the development of more advanced controllers such as
	PI controllers \cite{gustafsson1988pi}. However, this theory of
	step size control stability is based on a linear stability analysis
	not applicable to MPRK methods. Moreover, it only gives necessary
	bounds on the controller parameters and still requires further tools
	to derive efficient controllers. See also
	\cite{optimizedRK,ranocha2023stability}
	for some recent studies in the context of computational fluid dynamics.

	Further criteria for good controllers are discussed in
	\cite{soderlind2006adaptive}. In particular, we would like to achieve
	\emph{computational stability}, i.e., a continuous dependence of
	the computed results on the given parameters. An important aspect
	for computational stability is to avoid discontinuous effects in the
	controller and already imposed by the form of the DSP controller
	\eqref{eq:PIDlimit}.
	Furthermore, we would like to achieve \emph{tolerance convergence},
	i.e., to get better results when decreasing the tolerance
	(up to restrictions imposed by machine accuracy).
	Finally, a good controller should be able to achieve an error
	imposed by a given tolerance with as little work as possible. Here,
	we mainly measure the work by the number of accepted and rejected steps
	(proportional to the number of function evaluations and linear system
	solves for MPRK methods).

	\subsection{Methodology}\label{sec:methodology}
	We apply the second-order Rosenbrock method as well as the MPRK schemes introduced in Section~\ref{sec:schemes}, and further specified in Section~\ref{sec:ParasMPRK}, to several training problems which we will describe in Appendix~\ref{sec:training}. We choose $\code{atol}=\code{rtol}\eqqcolon \code{tol}$ in \eqref{eq:PIDepsomega} and consider values
	\begin{equation}\label{eq:TOL}
		\code{tol}\in\code{TOL}\coloneqq\{10^{-j}\mid j=1,\dotsc, 8\}.
	\end{equation}
	We formally evaluate the numerical solution $\b y_{\mathrm{num}}(t)$ for the times $t_1,t_2,\dotsc,t_K$.
	We measure the number of successful and rejected steps at time $t_n$ denoted by $S_n$ and $R_n$, respectively. Thus, $S_K$ and $R_K$ represent the number of successful and rejected steps at the end of the calculations. Moreover, we use the trapezoidal rule to approximate the relative $L_2$ error by introducing
	\begin{equation}\label{eq:L2errsol}
		\code{L2err\_rel}(\b t, \b Y,\byref)\coloneqq \left(\frac{\sum_{k=0}^{K-1} \frac{\dt_{k+1}}{2} \left(\norm{\byref(t_k)-\b y^k}_2^2 +\norm{\byref(t_{k+1})-\b y^{k+1}}_2^2 \right)}{\sum_{k=0}^{K-1} \frac{\dt_{k+1}}{2} \left(\norm{\byref(t_k)}_2^2 +\norm{\byref(t_{k+1})}_2^2 \right)}\right)^{\frac{1}{2}},
	\end{equation}
	where $\b Y=\left(\b y^0\dashline\dotsc\dashline \b y^K\right)$, $\b t=(t_0,\dotsc, t_K)^T$, and $\dt_{k+1}=t_{k+1}-t_k$. Then the trapezoidal rule yields
	\[ \left(\frac{\int_{t_0}^{t_K} \norm{\byref(t)-\b y_{\mathrm{num}}(t)}_2^2\dd t}{\int_{t_0}^{t_K} \norm{\byref(t)}_2^2\dd t}\right)^{\frac{1}{2}}= \code{L2err\_rel}(\b t, \b Y,\byref) + \O\left(\left(\max_{k=0,\dotsc,K-1}\dt_{k+1}\right)^3\right)\]
	for smooth enough integral kernels.
	If no analytical solution is available, a reference solution will be computed using the built-in function \code{ode15s} in \MATLAB R2023a \cite{matlabR2023b,shampine1997matlab} together with the inputs \code{AbsTol = 1e-13} and \code{RelTol = 1e-13} for the absolute and relative tolerances.

	A first task is to find abortion criteria to reduce the computational cost. To that end, we ran preliminary experiments, which suggest that any of the numerical methods equipped with standard controllers do not require more than $10^6$ total steps to solve any of the test problems (reviewed in Appendix~\ref{sec:training}). Thus, the calculations are aborted if $\Smax\coloneqq10^6$ steps are successful or if  $\Rmax\coloneqq10^4$ steps were rejected or $R_n\geq \frac{\Smax}{\Rmax}(S_n+1)=10^2 (S_n+1)$ holds at some point during the calculations. With that in mind, it might be possible that $t_K\neq \tend$ or $K=0$. In the latter case, we divide by $0$ in \eqref{eq:L2errsol}, which is why we introduce
	\begin{equation}\label{eq:err}
		\code{err}(\b t, \b Y,\byref)\coloneqq\begin{cases}
\code{L2err\_rel}(\b t, \b Y,\byref),& K>0,\\
\code{NaN}, & K=0.
		\end{cases}
	\end{equation}
	We then search for the optimal parameters of the time step controller introduced in Section~\ref{sec:controller} by means of minimizing the following cost function.

	Given a parameter $s>0$ (see Remark~\ref{rem:costfun} below for its influence), we introduce a cost function \begin{equation*}
		\begin{aligned}
			C_s\colon D&\to \R^+, \\
			\b x\coloneqq (\beta_1,\beta_2,\beta_3,\alpha_2,\kappa_2) &\mapsto C_s(\b x),
		\end{aligned}
	\end{equation*} where the domain $D$ is specified in \eqref{eq:PIDparas2}.
				For a  given set $\Tset$ of test problems and a $k$th order method we propose the cost function $C_s$ given by
				\begin{equation}\label{eq:cost_fun}
						C^\ast_s(\b x)=\sum_{\test\in\Tset}\psi\left(\sum_{\code{tol}\in \code{TOL}}(C_\step(\test,\code{tol},\b x)+	C_{\tol,s}(\test,\code{tol},\b x))\right),
				\end{equation}
			where $\psi(x)=\left(\arctan\bigl(\frac{x}{100}\bigr)\right)^2$ and
				\begin{equation}\label{eq:cost_fun_parts}
					\begin{aligned}
						C_\step(\test,\code{tol},\b x)=& k \cdot \ln\Bigl(S_K^\ast(\test,\code{tol},\b x)+ R_K^\ast(\test,\code{tol},\b x)\Bigr)\\&+\ln\left(\frac{ \code{err}(\code{tol},\test,\b x)}{\code{tol}}\right),\\
						C_{\tol,s}(\test,\code{tol},\b x)=& \max\left(0,\ln\left(\frac{ \code{err}(\code{tol},\test,\b x)}{s\cdot\code{tol}}\right)\right).\\
					\end{aligned}
				\end{equation}
				Here, the superscript asterisk indicates that the value of the related function is replaced by a penalty value, if the corresponding calculation was aborted. The particular penalty values and some properties of the cost function are summarized in the following remark.
				\begin{rem}\label{rem:costfun}
					Let us start discussing $C_\step$ by noting that in logarithmic scale, the slope between two consecutive points in the work-precision (WP) diagram equals $-k$ as $\dt\to 0$. Hence, $C_\step$ is such that any point on that straight line corresponds to the same cost, and everything below is cheaper.

					 However, this way it might happen that all points for the different tolerances are clustered in the upper left corner of the WP diagram. To punish this behavior, we also add $C_{\tol,s}$, which returns $0$ if $\emph{\code{err}}(\emph{\code{tol}},\test,\b x)\leq s\cdot\emph{\code{tol}}$. A sketch is given in Figure~\ref{fig:cost}, where all points on the black line will be associated with the same cost. However, if the error exceeds $s\cdot\code{\tol}$, the corresponding point lies on the red line and will yield higher costs. If we now add further lines for different tolerances into Figure~\ref{fig:cost}, we see that the presence of the red segments ensures that clustered points are expensive.
					 		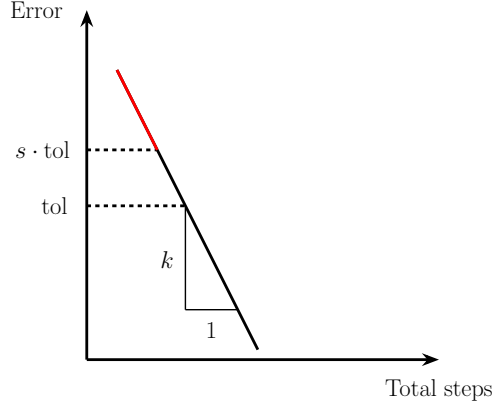
\begin{figure}[!ht]
					 			\centering
					 			\resizebox{0.5\textwidth}{!}{%
					 				\begin{circuitikz}
					 					\tikzstyle{every node}=[font=\LARGE]
					 					\draw [line width=2pt, dashed] (5,12) .. controls (6,12) and (6,12) .. (6.75,12);
					 					\draw [line width=2pt, dashed] (5,10.6) .. controls (6.25,10.6) and (6.25,10.6) .. (7.5,10.6);
					 					\draw [line width=2pt, ->, >=Stealth] (5,6.75) .. controls (9.5,6.75) and (9.5,6.75) .. (13.75,6.75);
					 					\draw [line width=2pt, ->, >=Stealth] (5,6.75) .. controls (5,11) and (5,11) .. (5,15.5);

					 					\draw [line width=1pt, short] (7.45,10.55) .. controls (7.45,9.25) and (7.45,9.25) .. (7.45,8);
					 					\draw [line width=1pt, short] (7.45,8) .. controls (8.25,8) and (8.25,8) .. (8.75,8);
					 					\node [font=\LARGE] at (7,9.25) {$k$};
					 					\node [font=\LARGE] at (8.1,7.5) {$1$};
										\node [font=\LARGE] at (4.2,10.6) {$\tol$};
										\node [font=\LARGE] at (3.9,12) {$s\cdot \tol$};
										\node [font=\LARGE] at (13.75,6) {Total steps};
										\node [font=\LARGE] at (3.75,15.5) {Error};

					 						\draw [ color={rgb,255:red,0; green,0; blue,0}, line width=2pt, short] (5.75,14) .. controls (6,13.5) and (6,13.5) .. (9.25,7);
					 					\draw [ color={rgb,255:red,255; green,0; blue,0}, line width=2pt, short] (6.75,12) .. controls (6,13.5) and (6,13.5) .. (5.75,14);
					 				\end{circuitikz}
					 			}%
					 			\caption{Sketch of a WP diagram of a $k$th order method (in logarithmic scale).}
					 			\label{fig:cost}
					 		\end{figure}
				 		 In particular, assume that $\emph{\code{err}}(\emph{\code{tol}},\test,\b x)= 10^m s\cdot\emph{\code{tol}}$ for some $m\in \R^+$. Then we find $C_{\tol,s}=m\ln(10)$, which means that exceeding a given tolerance by a factor of $10^ms$ will be punished proportionally to the exponent $m$.

					Additionally, we have the following penalty values and conditions.
					\begin{itemize}
						\item  If the calculation is aborted because of $S_K=\Smax$ or $R_K=\Rmax$ we punish this by setting \[S^\ast=10\max\{\Smax,\Rmax\}\quad \text{ or } \quad R^\ast=10\max\{\Smax,\Rmax\},\] respectively.
						\item If $\dt_n<10^{-100}$ at some point, the calculation is aborted and we set \[S^\ast=10\max\{\Smax,\Rmax\}.\]
						\item If the slope of the straight line through two consecutive points in the WP diagram does not lie in the interval $(-\infty,-0.7)$, the DSP parameter combination is disqualified. An exception is the first slope in the WP diagram, which we only require to lie in the interval $(-\infty,-0.35)$. This approach is based on the fact that the slope should tend to $-k$ as $\dt\to 0$. The disqualification is done in our case by canceling the calculations and adding $M=10$ to the current value of the cost function $C_s$, which represents the value $C_s^\ast$. The penalty $M=10$ is a reasonable choice since we will consider four test problems, eight different tolerances, and methods of order at most three, so that
						\begin{equation*}
							\begin{aligned}
								\sum_{\code{\tol}\in \code{TOL}}(C_\step(\test,\code{\tol},\b x)+	C_{\tol,s}(\test,\code{\tol},\b x))&\leq 8 \cdot(3\ln(20\Smax)+20\ln(10))\\&\leq 8\cdot(3\cdot 7\ln(20)+60)<10^3
							\end{aligned}
						\end{equation*}
					for any controller satisfying $\emph{\code{err}}(\emph{\code{tol}},\test,\b x)\leq 10^{10} s\cdot\emph{\code{tol}}$.
						In view of the transformation $\psi$, see Figure~\ref{fig:atan} for a plot of the graph, we thus observe that for each test, we obtain a value between $0$ and $2.5$ justifying the penalty addend $M=4\cdot 2.5$. With this transformation we make sure that a single test problem is not dominating the remaining, for instance due to a large value of $C_{\tol,s}$.
					\end{itemize}
	\begin{figure}[!h]
	\centering
	\begin{scaletikzpicturetowidth}{0.5\textwidth}
		\begin{tikzpicture}[scale=\tikzscale]

			\begin{axis}[				ytick distance=1,
				xtick={0,100,300,500,1000},
				xticklabels = {0,100,300,500,1000},
				enlarge x limits=false,ymin=0,ymax=2.4,xmax=1050,xmin=0.001,
				minor tick num=1, axis lines = left,xlabel=$x$,ylabel=$y$, ylabel near ticks,
				ylabel style={rotate=-90},
				every axis plot/.append style={thick}]
				\addplot [scale=1, domain=0.001:1000, smooth, variable=\x, black, solid, mylabel=at 0.9 below right with {$\psi$}] plot ({\x}, {rad(atan(\x/100))*rad(atan(\x/100))});
			\end{axis}
		\end{tikzpicture}
	\end{scaletikzpicturetowidth}
	\caption{Plot of $\psi(x)=\left(\arctan\left(\frac{x}{100}\right)\right)^2$. }\label{fig:atan}
\end{figure}
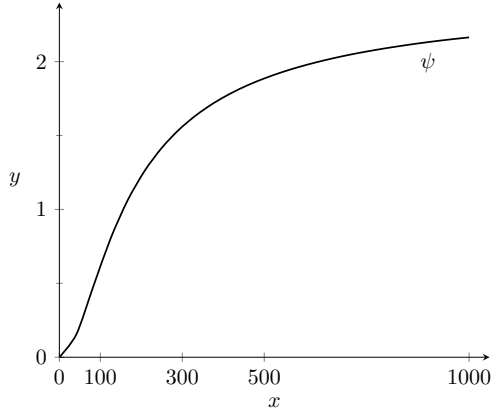
\end{rem}

\begin{rem}
	The last ingredient disqualifying points with a wrong slope has to
	be adapted for methods with reduced stability properties such as
	explicit RK methods. Indeed, good controllers typically lead to a
	clustering of points in a work-precision diagram as observed for example
	in \cite{bleecke2023step}.
\end{rem}

				The minimum of the cost function $C_{s}$ is searched using the Bayesian optimization function \code{bayesopt} in \MATLAB \cite{matlabR2023b} with $500$ iterations and the input arguments \code{'IsObjectiveDeterministic'} set to \code{true} as well as \code{'AcquisitionFunctionName'} set to \code{'expected-improvement-per-second-plus'} for a somewhat balanced trade off between exploration and exploitation, followed by a run with another $500$ iterations using the \code{'AcquisitionFunctionName'} set to \code{'probability-of-improvement'} enhancing the local search. Since we validate the results by means of plotting the WP diagram, the value $s=1$ is a natural choice here, meaning that there is no safety gap between $\code{\tol}$ and $s\cdot \code{\tol}$ in Figure~\ref{fig:cost}. The search of the best set of hyperparameters continues as long as the cost function has improved by more than $10^{-3}$ within the cycle of 1000 iterations.

\section{Controller Parameters}\label{sec:PIDpara}
As described in Section~\ref{sec:methodology}, we use Bayesian optimization to determine customized parameters for the time step controller reviewed in Section~\ref{sec:controller}.

\begin{figure}[!hbtp]
	\centering
	\includegraphics[width=0.5\textwidth]{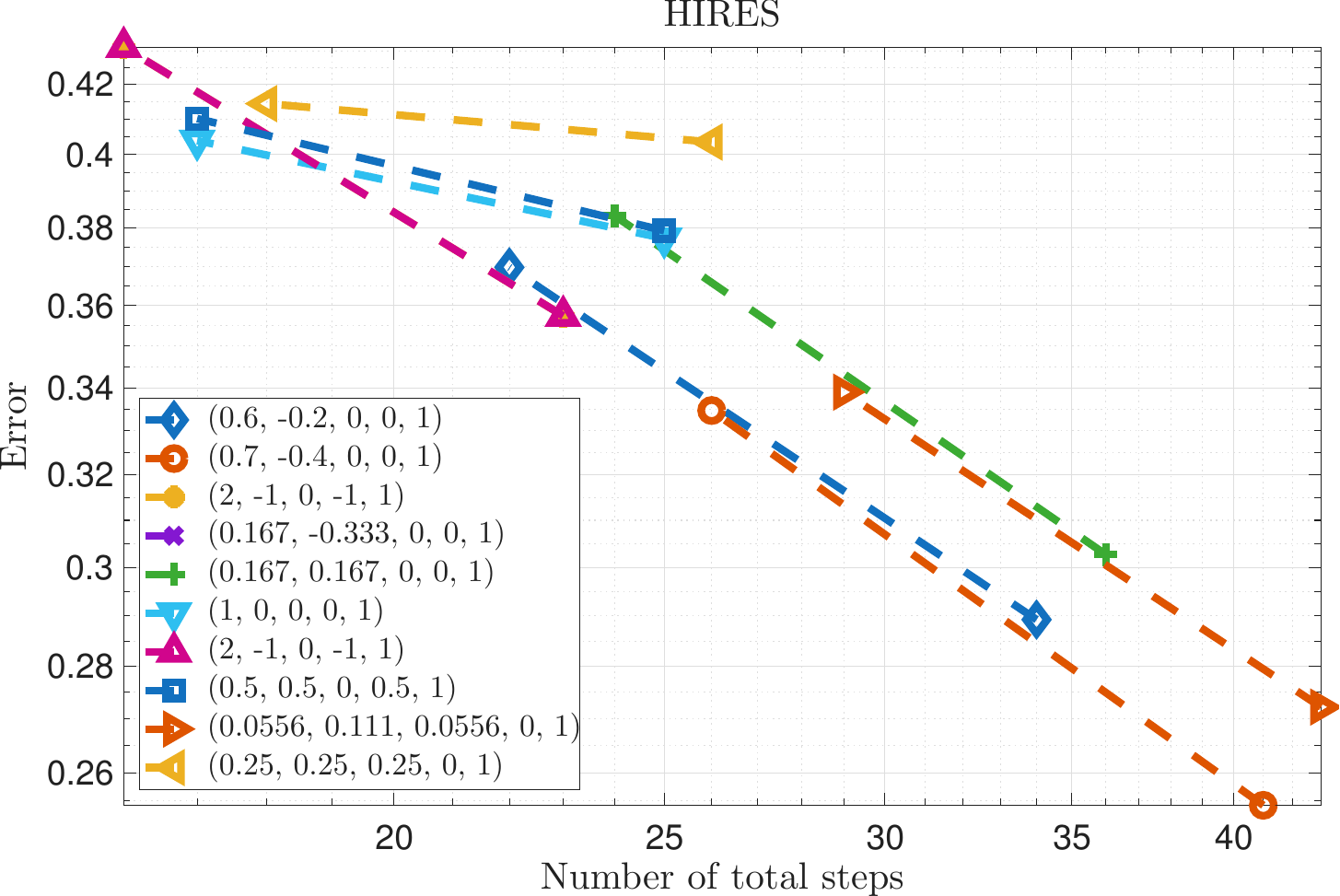}
	\caption{WP diagram for the HIRES problem using MPRK22(1) with several standard controllers and tolerance set $\{10^{-1},10^{-2}\}$.}\label{fig:zoomhires}
\end{figure}

\begin{figure}[!htbp]
	\begin{subfigure}[t]{0.505\textwidth}
		\includegraphics[width=\textwidth]{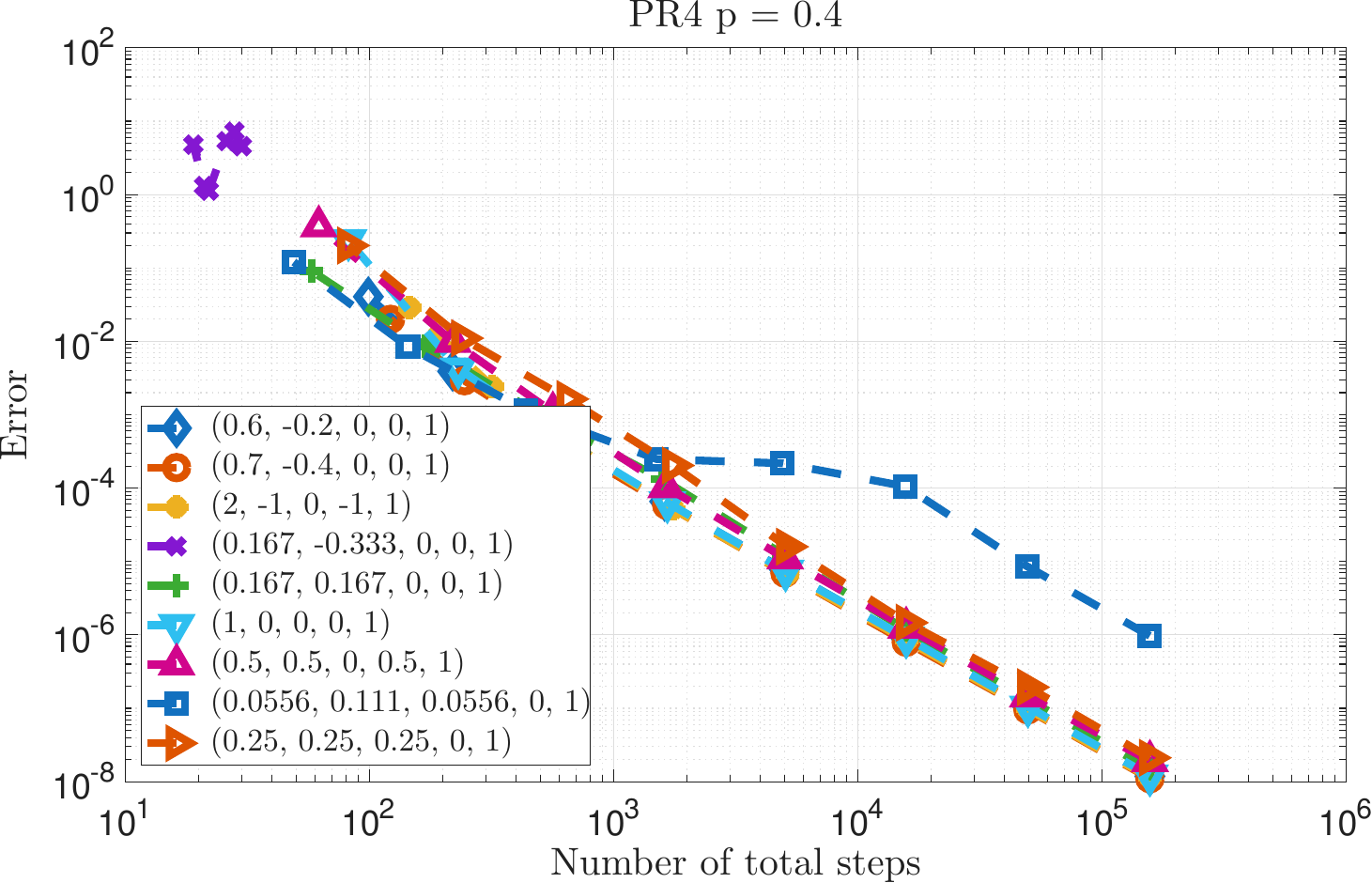}
	\end{subfigure}
	\begin{subfigure}[t]{0.485\textwidth}
		\includegraphics[width=\textwidth]{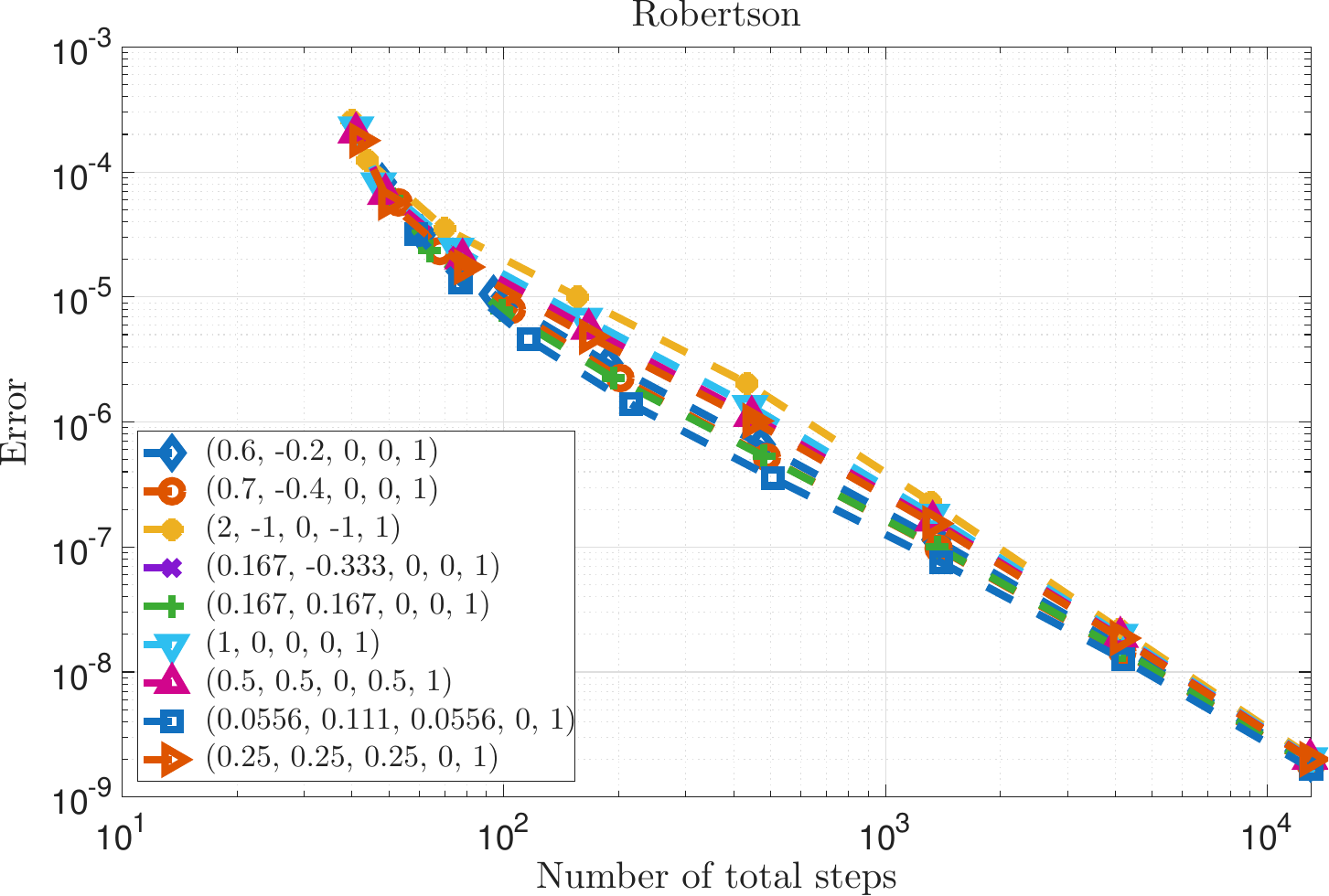}
	\end{subfigure}\\
	\begin{subfigure}[t]{0.495\textwidth}
		\includegraphics[width=\textwidth]{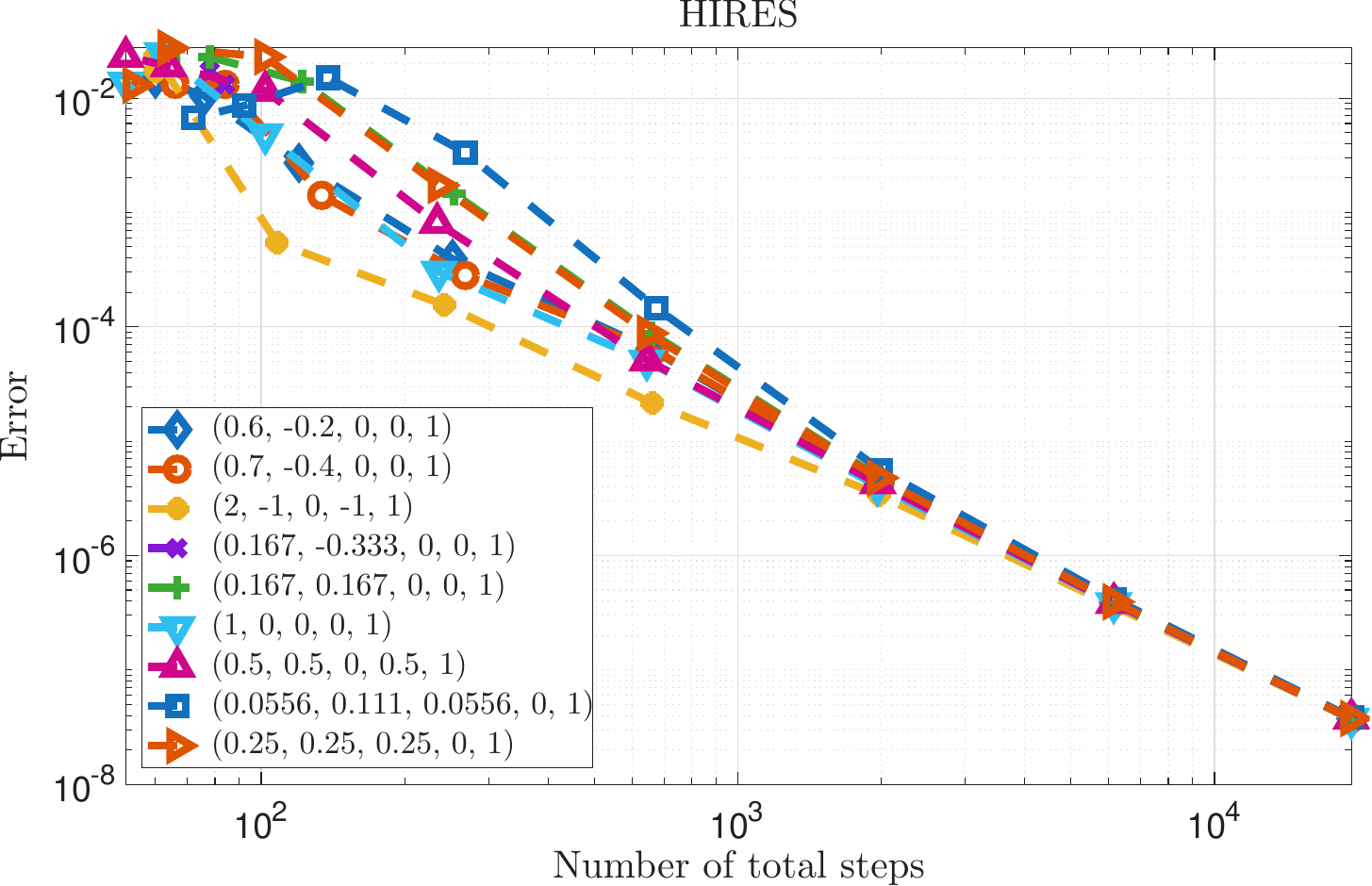}
	\end{subfigure}
	\begin{subfigure}[t]{0.495\textwidth}
		\includegraphics[width=\textwidth]{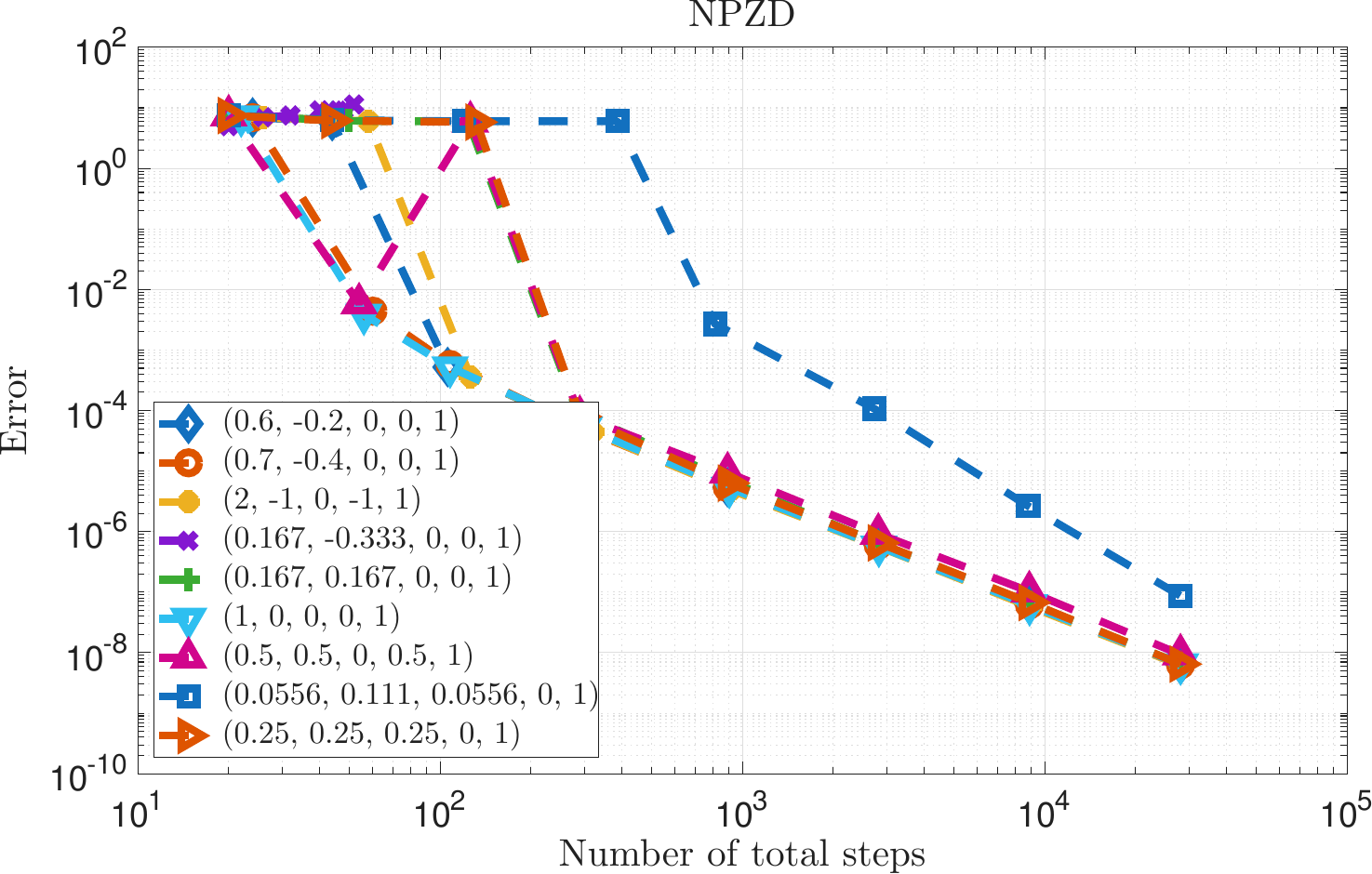}
	\end{subfigure}
	\caption{WP diagram for the test problems from Section~\ref{sec:training} using ROS2(1/$(2+\sqrt{2})$) with several standard controllers and tolerance set \code{TOL} from \eqref{eq:TOL}.}\label{Fig:ROS2_std}
\end{figure}

\begin{figure}[!htbp]
	\begin{subfigure}[t]{0.495\textwidth}
		\includegraphics[width=\textwidth]{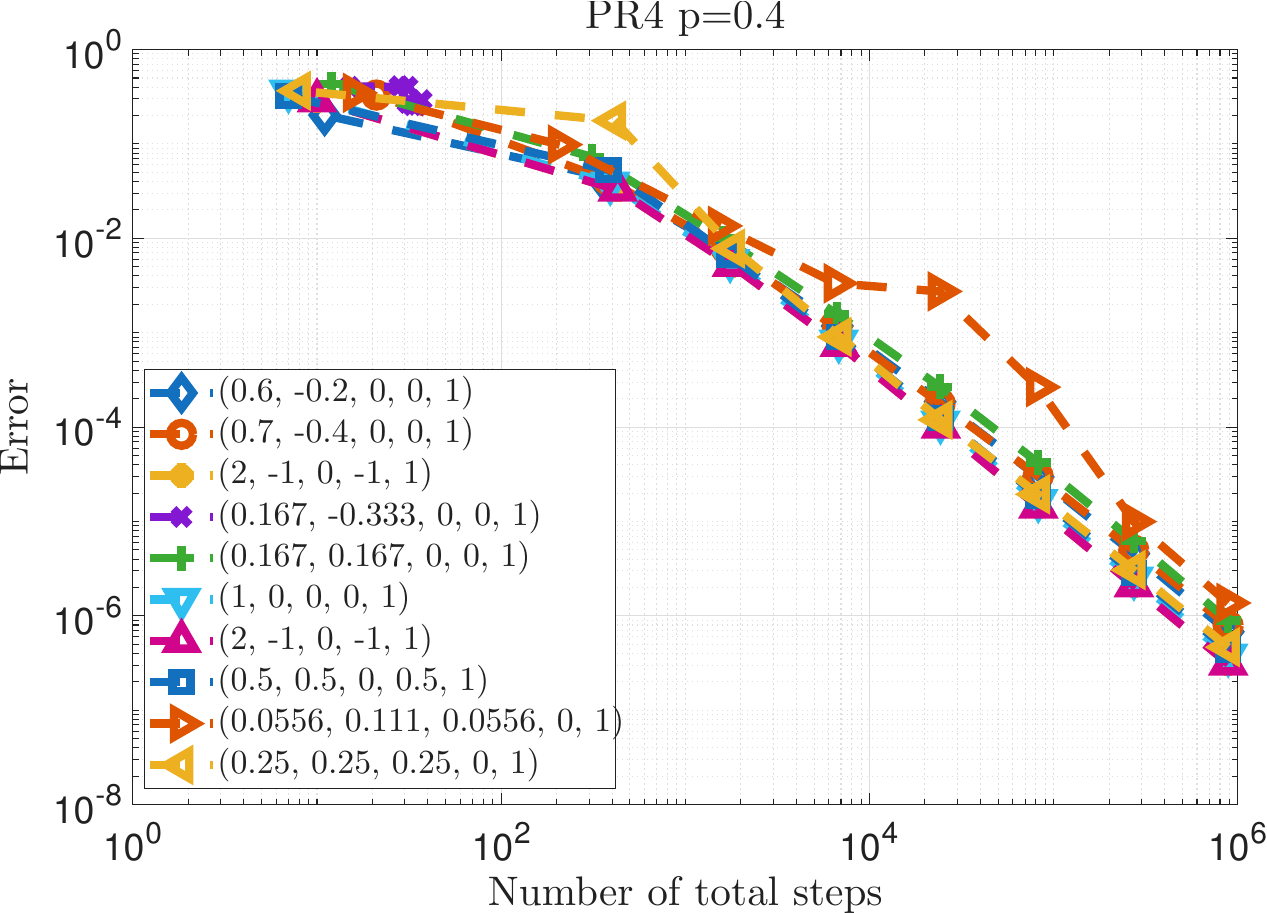}
	\end{subfigure}
	\begin{subfigure}[t]{0.5\textwidth}
		\includegraphics[width=\textwidth]{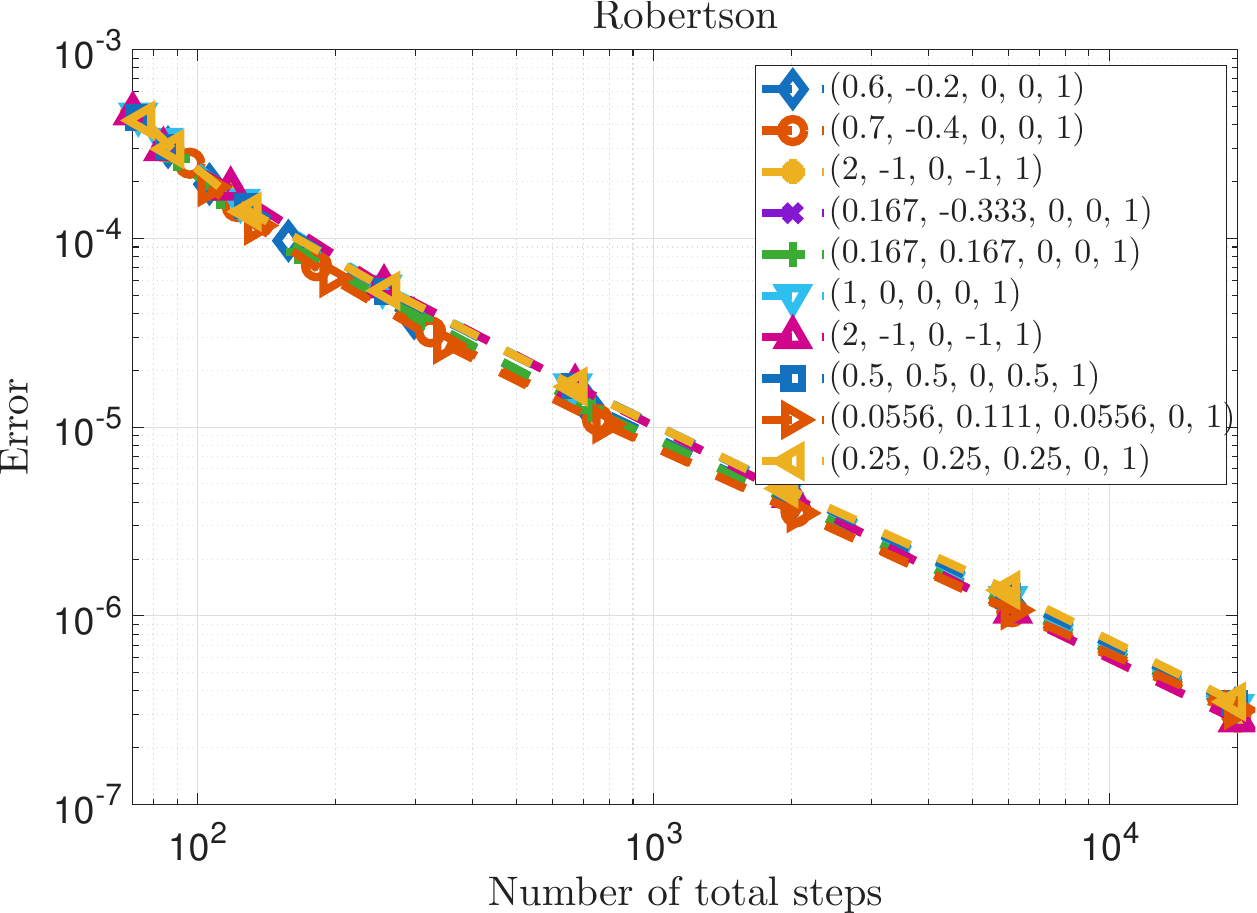}
	\end{subfigure}\\
	\begin{subfigure}[t]{0.5\textwidth}
		\includegraphics[width=\textwidth]{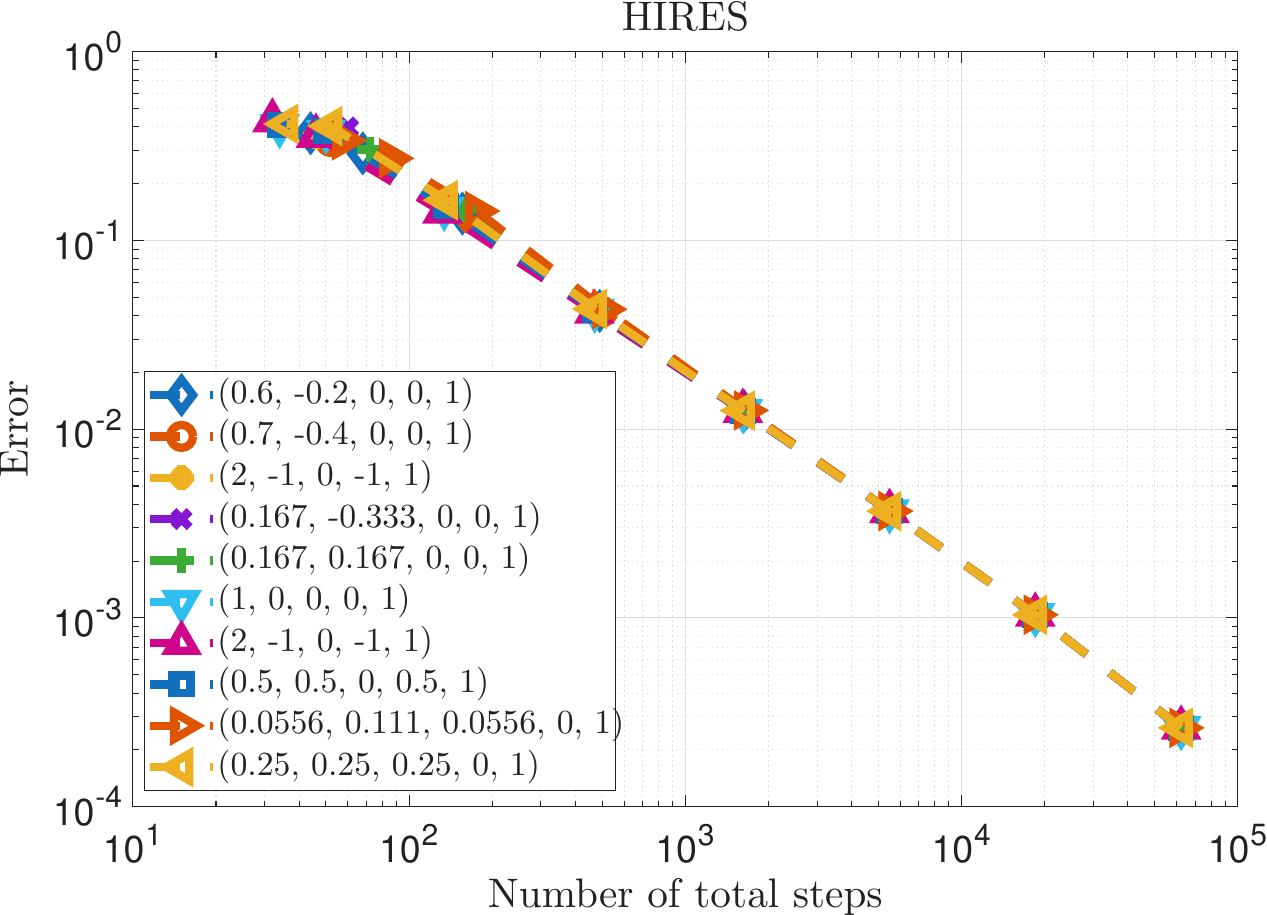}
	\end{subfigure}
	\begin{subfigure}[t]{0.5\textwidth}
		\includegraphics[width=\textwidth]{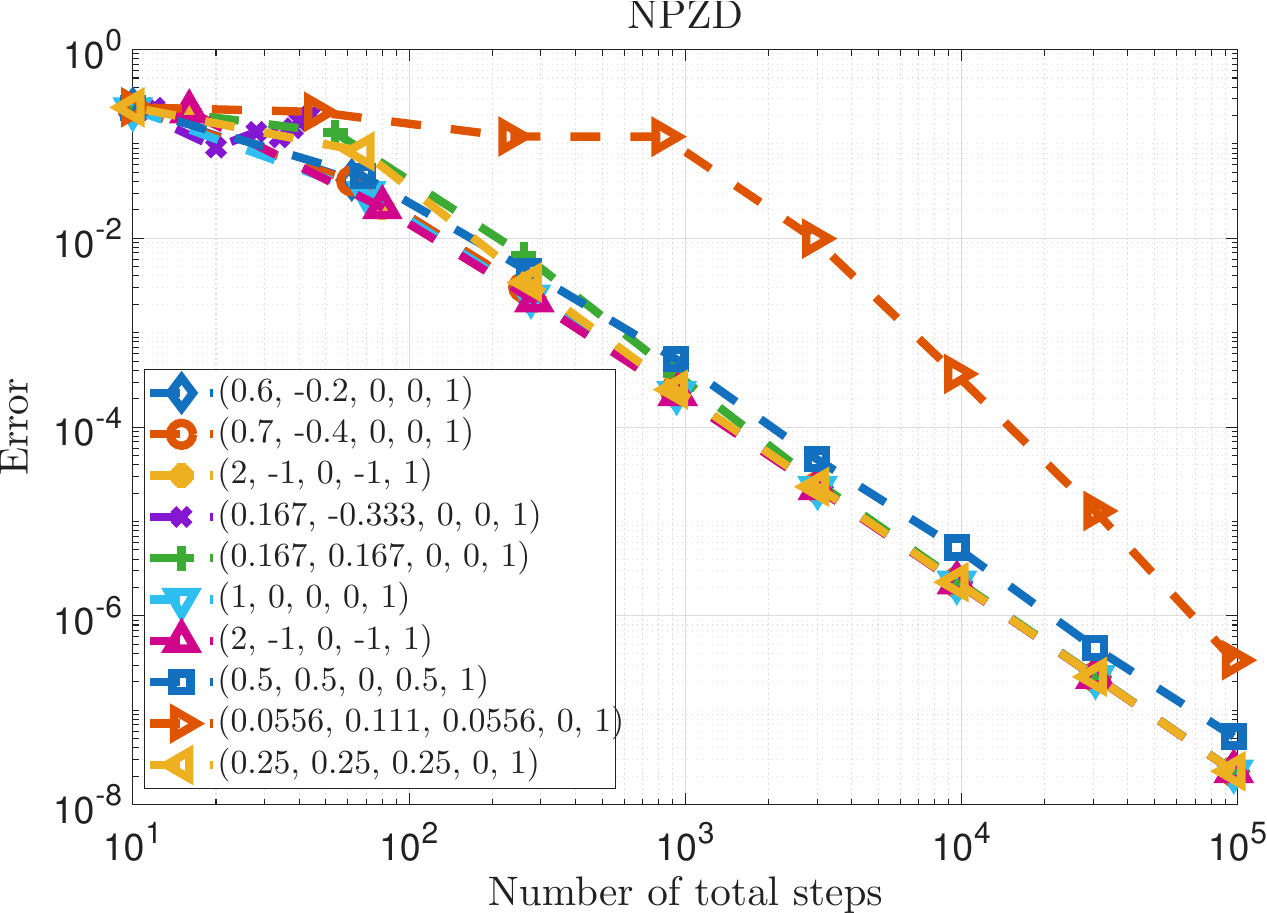}
	\end{subfigure}
	\caption{WP diagram for the test problems from Section~\ref{sec:training} using MPRK22(1) with several standard controllers and tolerance set \code{TOL} from \eqref{eq:TOL}.}\label{Fig:MPRK22_std}
\end{figure}

\begin{figure}[!htbp]
	\begin{subfigure}[t]{0.495\textwidth}
		\includegraphics[width=\textwidth]{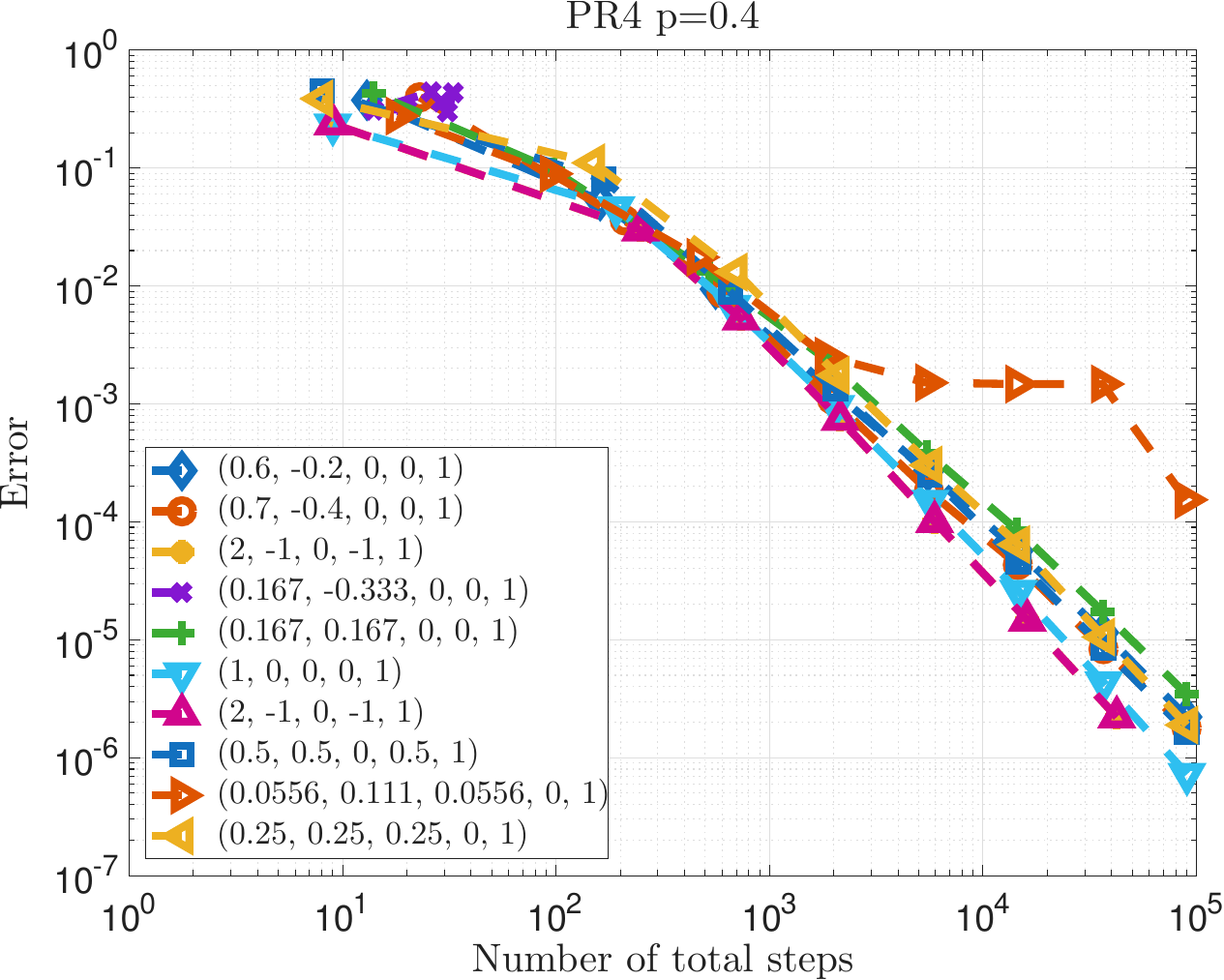}
	\end{subfigure}
	\begin{subfigure}[t]{0.495\textwidth}
		\includegraphics[width=\textwidth]{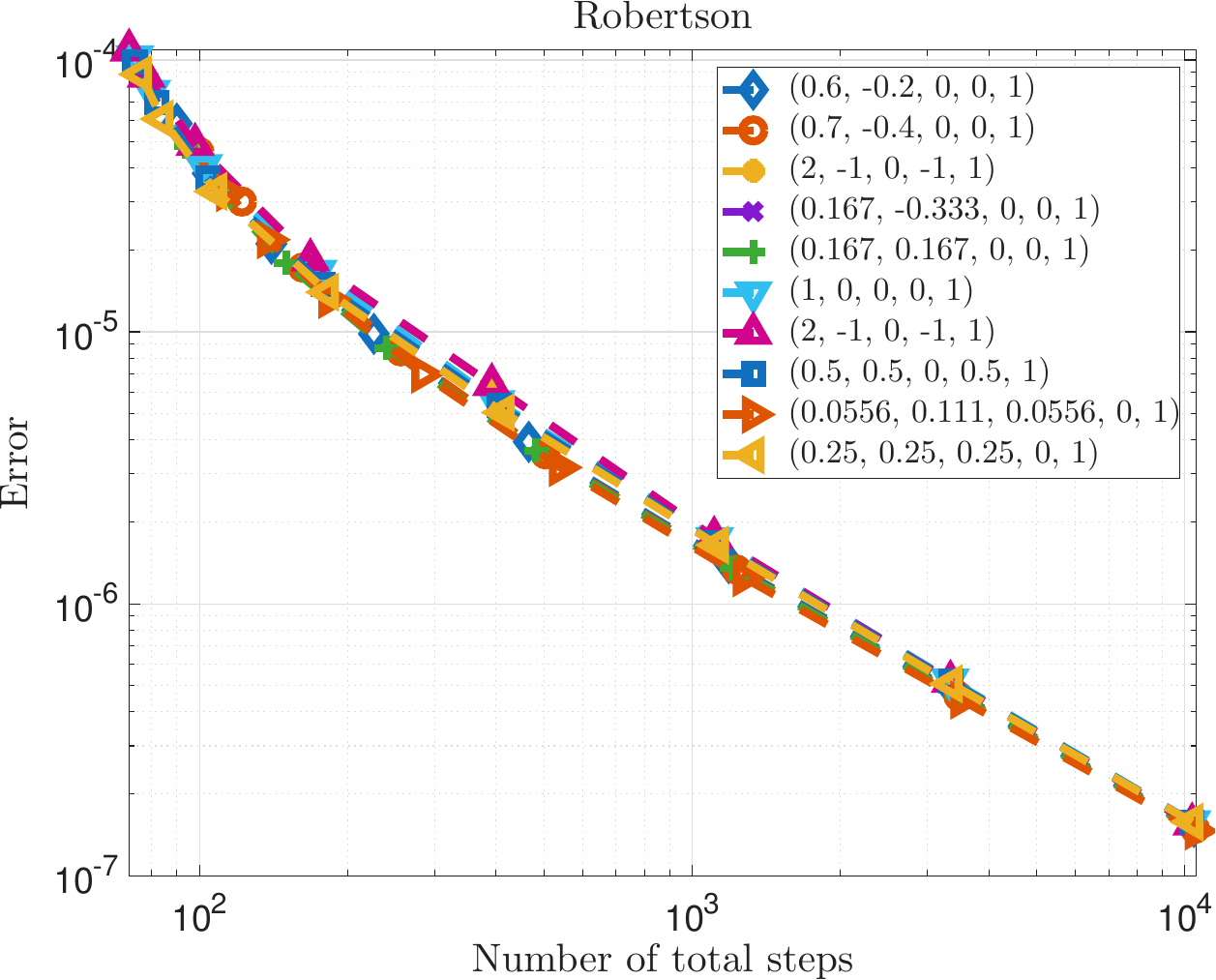}
	\end{subfigure}\\
	\begin{subfigure}[t]{0.495\textwidth}
		\includegraphics[width=\textwidth]{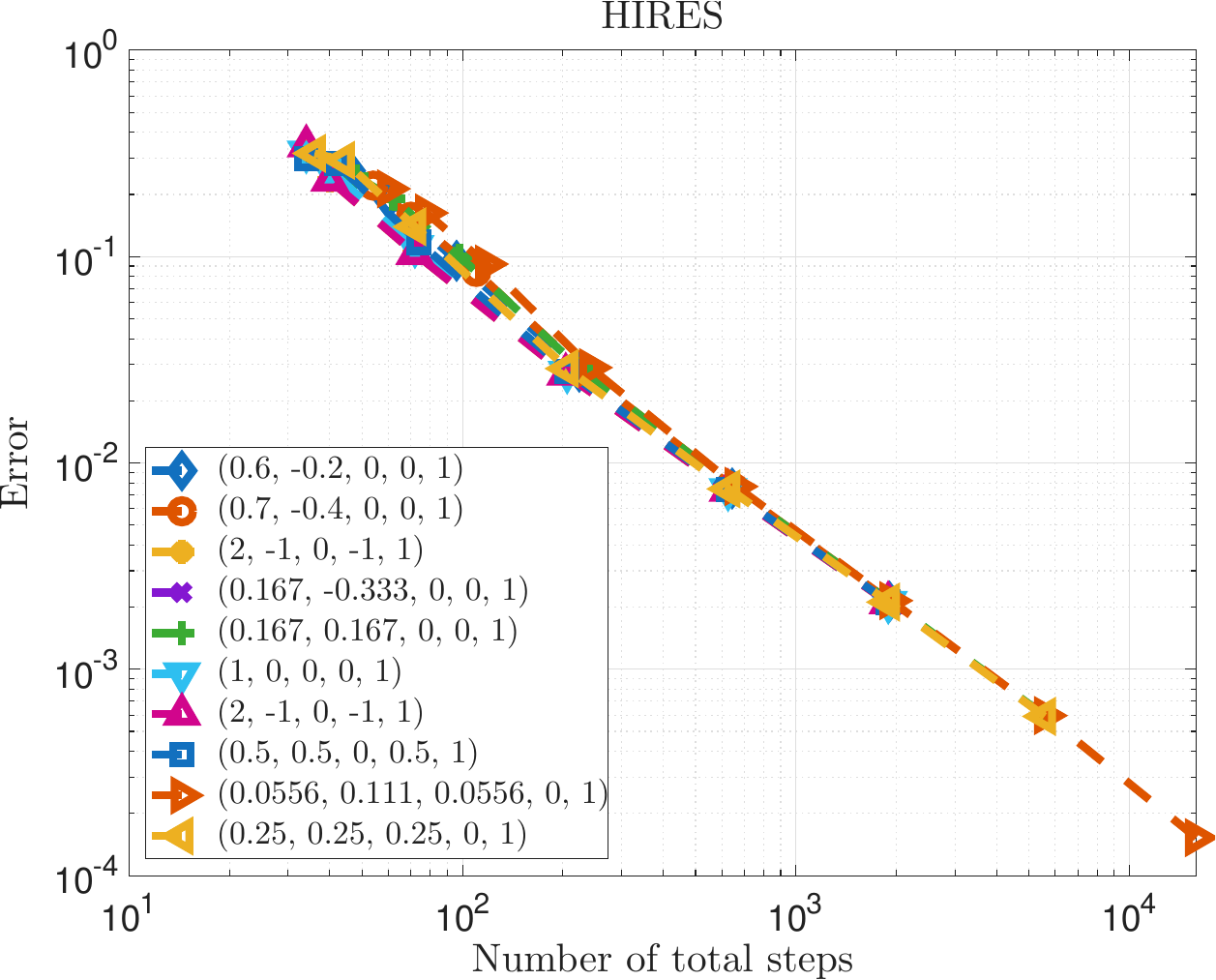}
	\end{subfigure}
	\begin{subfigure}[t]{0.495\textwidth}
		\includegraphics[width=\textwidth]{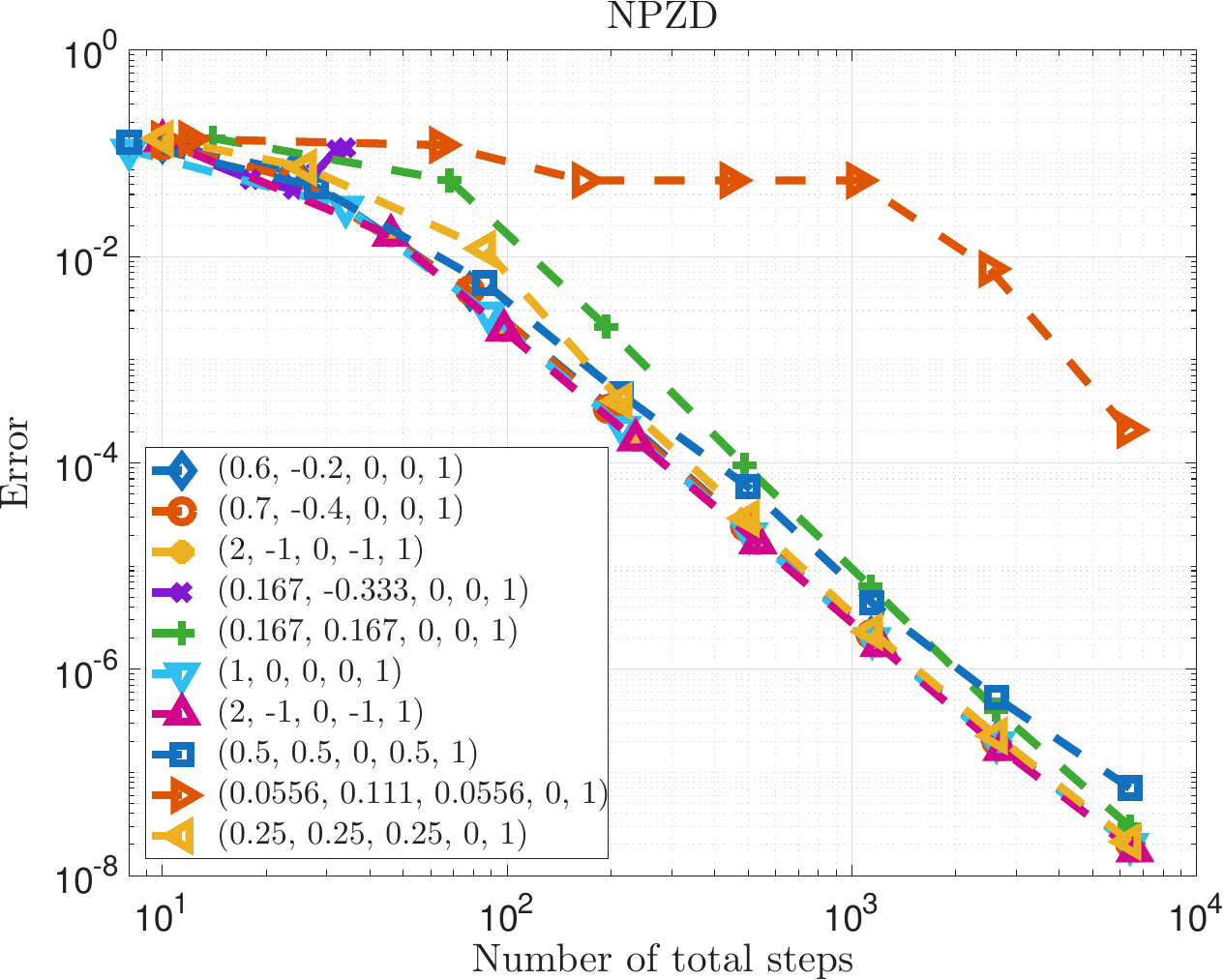}
	\end{subfigure}
	\caption{WP diagram for the test problems from Section~\ref{sec:training} using MPRK43(0.5,0.75) with several standard controllers and tolerance set \code{TOL} from \eqref{eq:TOL}.}\label{Fig:MPRK43I_std}
\end{figure}

\begin{figure}[!htbp]
	\begin{subfigure}[t]{0.495\textwidth}
		\includegraphics[width=\textwidth]{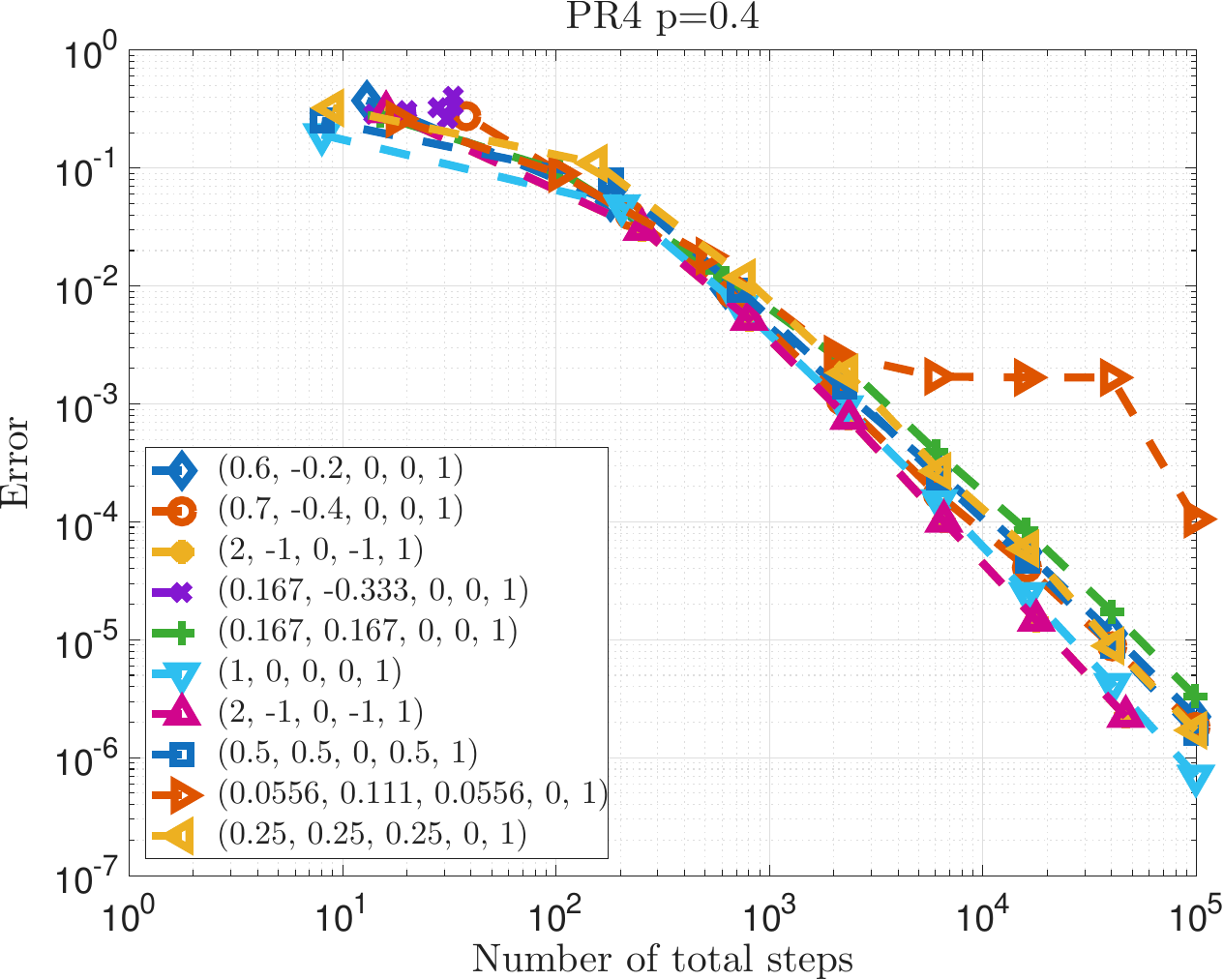}
	\end{subfigure}
	\begin{subfigure}[t]{0.495\textwidth}
		\includegraphics[width=\textwidth]{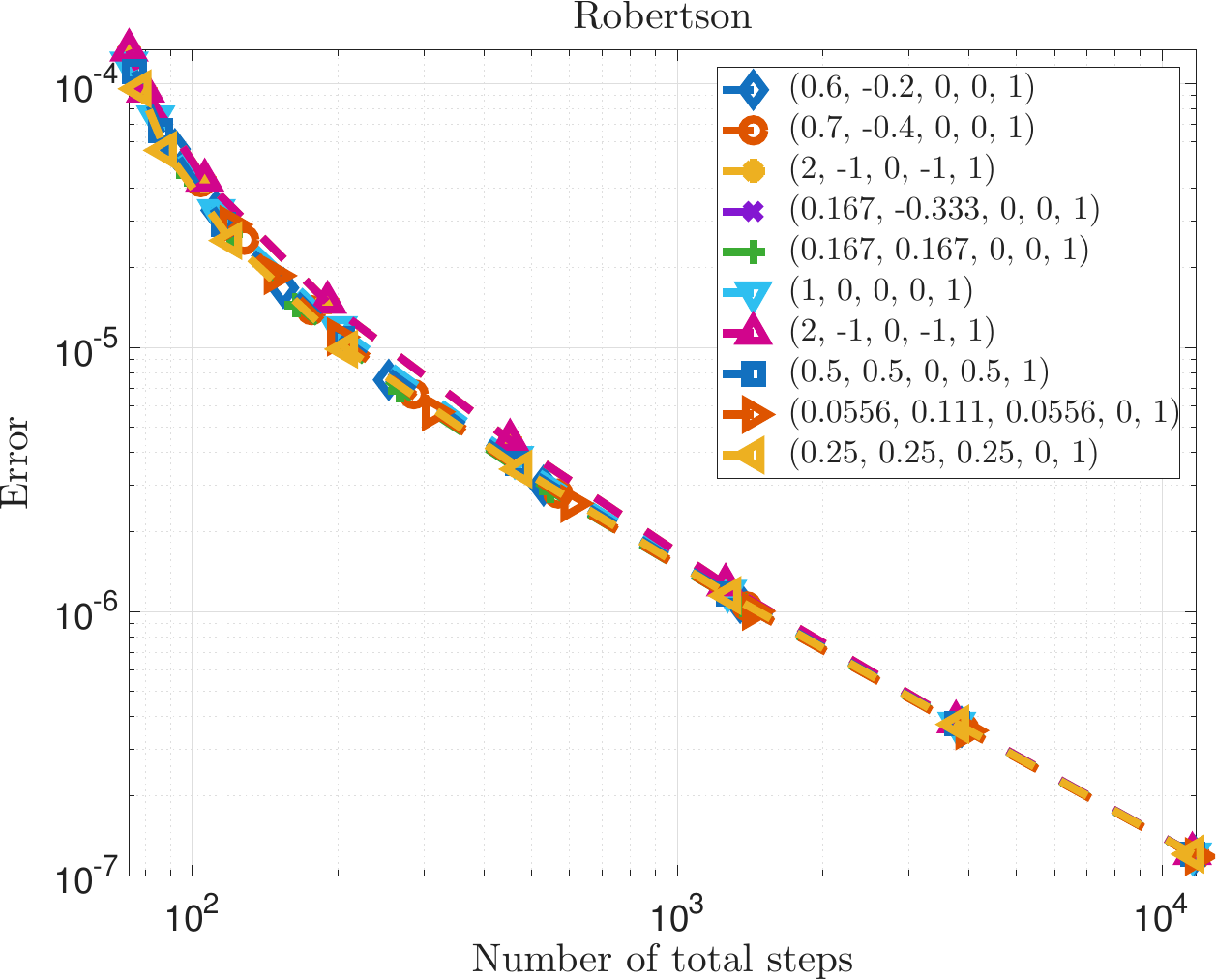}
	\end{subfigure}\\
	\begin{subfigure}[t]{0.495\textwidth}
		\includegraphics[width=\textwidth]{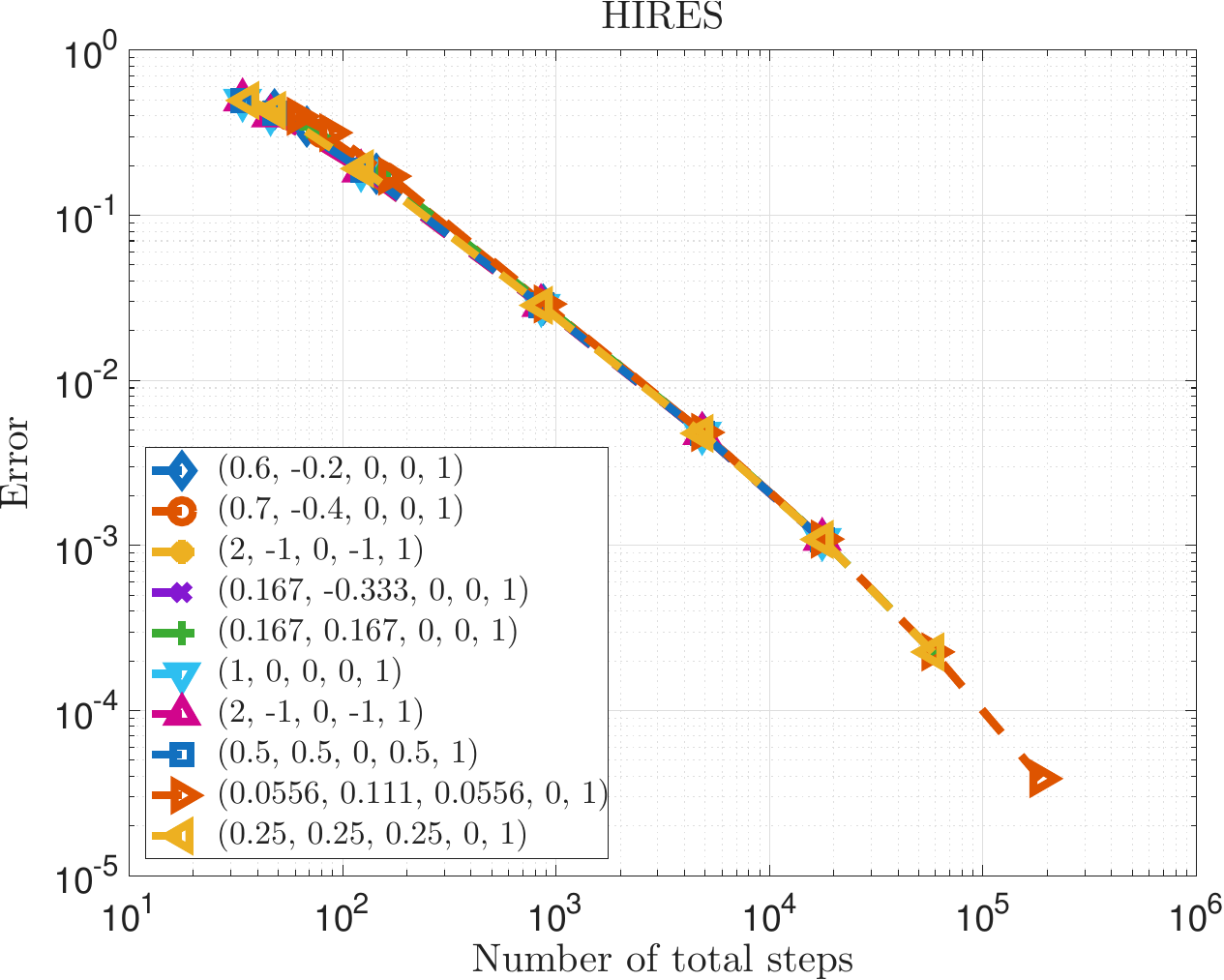}
	\end{subfigure}
	\begin{subfigure}[t]{0.495\textwidth}
		\includegraphics[width=\textwidth]{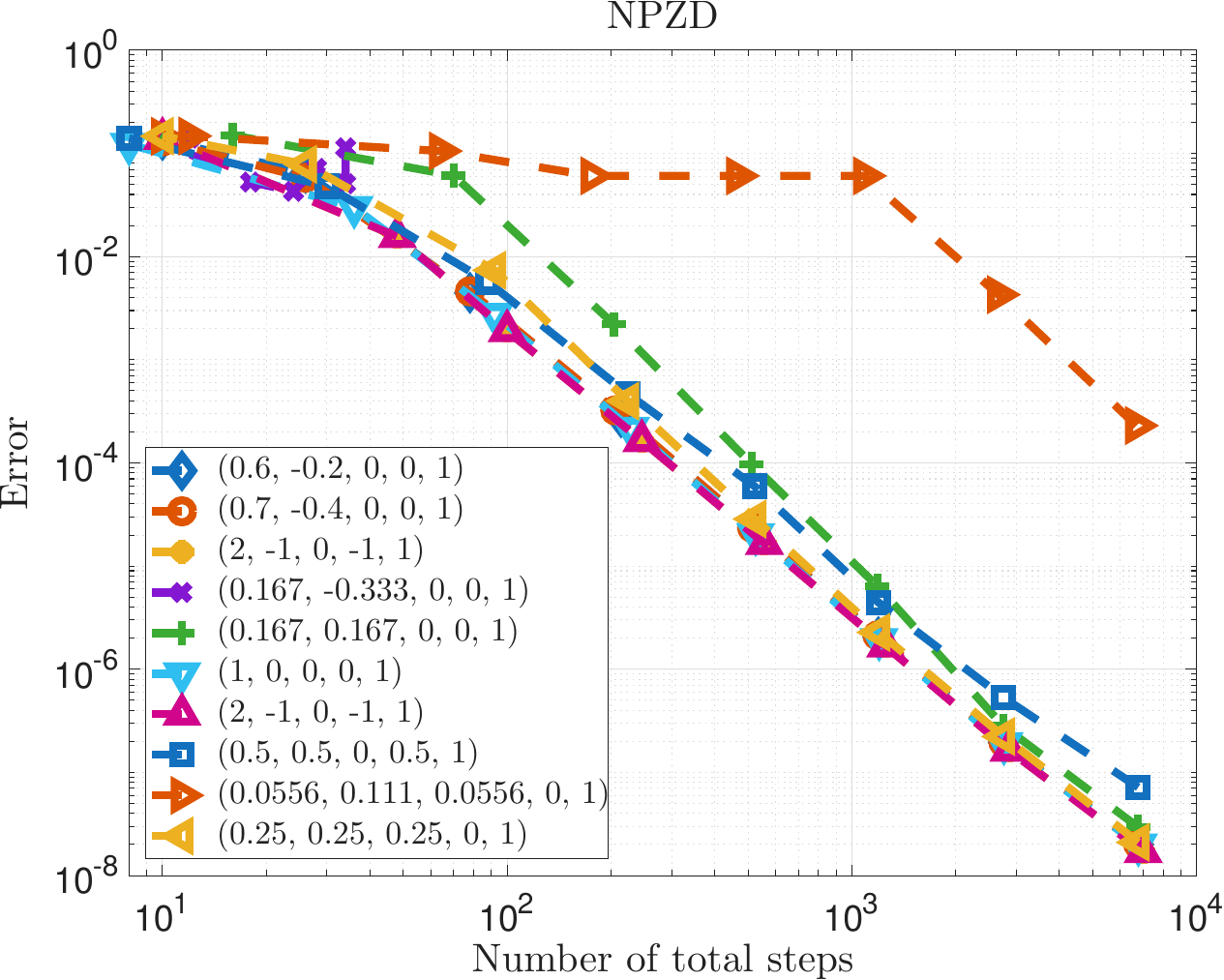}
	\end{subfigure}
	\caption{WP diagram for the test problems from Section~\ref{sec:training} using MPRK43(0.563) with several standard controllers and tolerance set \code{TOL} from \eqref{eq:TOL}.}\label{Fig:MPRK43II_std}
\end{figure}

We first focus on the performance of the standard parameters from \cite{gustafsson1988pi,gustafsson1991control,soderlind2006time,soderlind2006adaptive,S2002,soderlind2003digital,Z1964, G1994, AS2017} of the form \eqref{eq:PIDlimit}, see Figures~\ref{Fig:ROS2_std}, \ref{Fig:MPRK22_std}, \ref{Fig:MPRK43I_std}, and \ref{Fig:MPRK43II_std}.
The parameters
 $(\tfrac{1}{18},\tfrac19,\tfrac{1}{18},0,1)\approx(0.056,0.111,0.056,0,1)$ \cite{soderlind2003digital} and $(\tfrac16, -\tfrac13,0,0,1)$  are by far the worst. Indeed, the fact that the top three standard parameters in Figure~\ref{fig:zoomhires} are also the best in terms of the costs from Table~\ref{tab:sd_paras} supports our approach.
While all standard controllers were disqualified for ROS2, see Table~\ref{tab:sd_paras} and Figure~\ref{Fig:ROS2_std}, we observe that, besides the parameters with lowest costs, the parameters $p_1=(0.7, -0.4,0,0,1)$ and $p_2=(0.6,-0.2,0,0,1)$ seem to be promising candidates.
\begin{table}[!htbp]
	\centering
	\caption{Costs of standard parameters using the cost function $C_{1}$.}\label{tab:sd_paras}
	\begin{tabular}{lllll}
		Controller & MPRK22(1) & MPRK43Iadap(0.5,075) & MPRK43IIadap(0.563) & ROS2 \\
		\hline
		(0.6,-0.2,0,0,1) & 3.7463 & 4.3631 & 4.3659 & 13.6323 \\
		(0.7,-0.4,0,0,1) & 3.7557 & 4.36 & 4.3712 & 12.5167 \\
		($\tfrac16$,-$\tfrac13$,0,0,1) & 14.3934 & 14.5452 & 14.5454 & 14.387 \\
		($\tfrac16$,$\tfrac16$,0,0,1) & 14.3523 & 14.6537 & 14.654 & 12.5202 \\
		(1,0,0,0,1) & 12.5475 & \textbf{4.3014} & \textbf{4.3124} & 12.5172 \\
		(2,-1,0,-1,1) & \textbf{3.7188} & 4.4895 & 4.4916 & 12.5188 \\
		(0.5,0.5,0,0.5,1) & 12.5477 & 12.5671 & 4.4146 & 13.7713 \\
		(0.056,0.111,0.056,0,1) & 14.352 & 14.6664 & 14.6669 & 12.5234 \\
		(0.25,0.25,0.25,0,1) & 12.5491 & 4.394 & 16.7426 & 12.518 \\
		\hline
	\end{tabular}
\end{table}

As a next step, we run our experiments using the cost function $C_1$ introduced in Section~\ref{sec:methodology}.
The resulting customized parameters are summarized for the four schemes under consideration in Table~\ref{tab:optiparas2}. For more detailed information, we refer to Tables~\ref{tab:optiparasROS2},~\ref{tab:optiparasMPRK22},~\ref{tab:optiparasMPRK43I}, and \ref{tab:optiparasMPRK43II} in Appendix~\ref{sec:data}. Therein, the column "Bayes"  has entries of the form \enquote{scheme\_number} and refers to the scheme for which the optimization was run and how many evaluations of the cost functions were needed until there was no significant improvement after 1000 iterations.



		\begin{table}[!hbtp]
			\centering
			\caption{Overall cheapest parameters and cost $C_1$}\label{tab:optiparas2}
		\bgroup
		\def\arraystretch{1.5}
			\begin{tabular*}{\linewidth}{@{\extracolsep{\fill}}ccccccc}
				\toprule
				\multicolumn{1}{c}{Method} &
				\multicolumn{1}{c}{$\beta_1$} &
				\multicolumn{1}{c}{$\beta_2$} &
				\multicolumn{1}{c}{$\beta_3$} &
				\multicolumn{1}{c}{$\alpha_2$} &
				\multicolumn{1}{c}{$\kappa_2$} &
				$C_1$
				\\
				\midrule
					ROS2(1/(2+$\sqrt2$)) &
			1.4783 & -1.027 & -0.27743 & -0.23112 & 3 & 2.279
				\\
				MPRK22(1) &
				1.5193 & -0.42576 & -0.078535 & -0.29465 & 2 & 3.6998
				\\
				MPRK43(0.5,0.75) &
				1.8476 & -0.11068 & -0.2863 & -0.24606 & 2		& 4.2631
				\\
				MPRK43(0.563) &
								1.5193 & -0.42576 & -0.078535 & -0.29465 & 2 & 4.2745
				\\
				\bottomrule
			\end{tabular*}
		\egroup
		\end{table}

\begin{figure}[!htbp]
		 	\begin{subfigure}[t]{0.5\textwidth}
		 		\includegraphics[width=\textwidth]{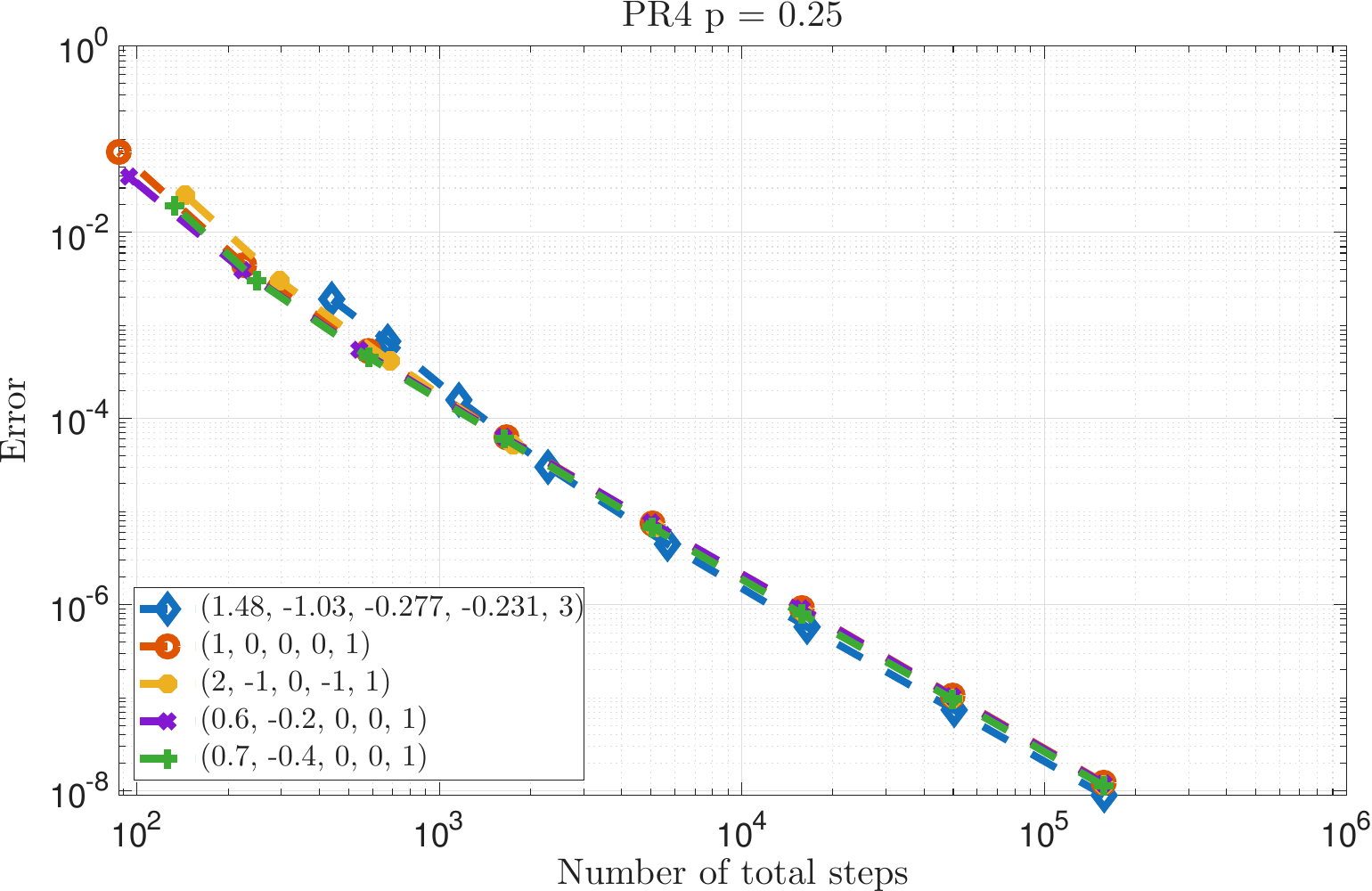}
		 	\end{subfigure}
		 	\begin{subfigure}[t]{0.5\textwidth}
		 		\includegraphics[width=\textwidth]{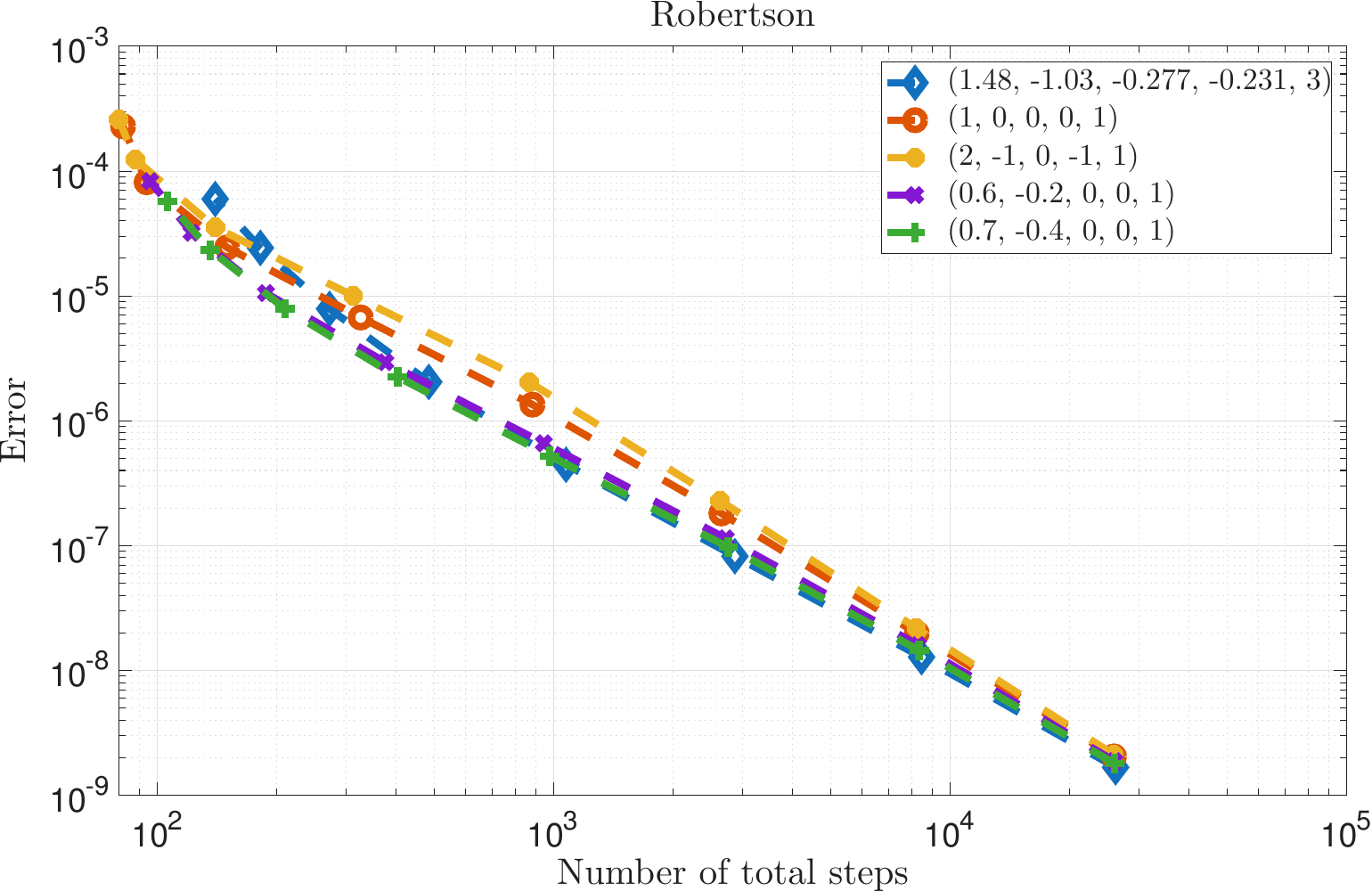}
		 	\end{subfigure}\\
		 	\begin{subfigure}[t]{0.5\textwidth}
		 		\includegraphics[width=\textwidth]{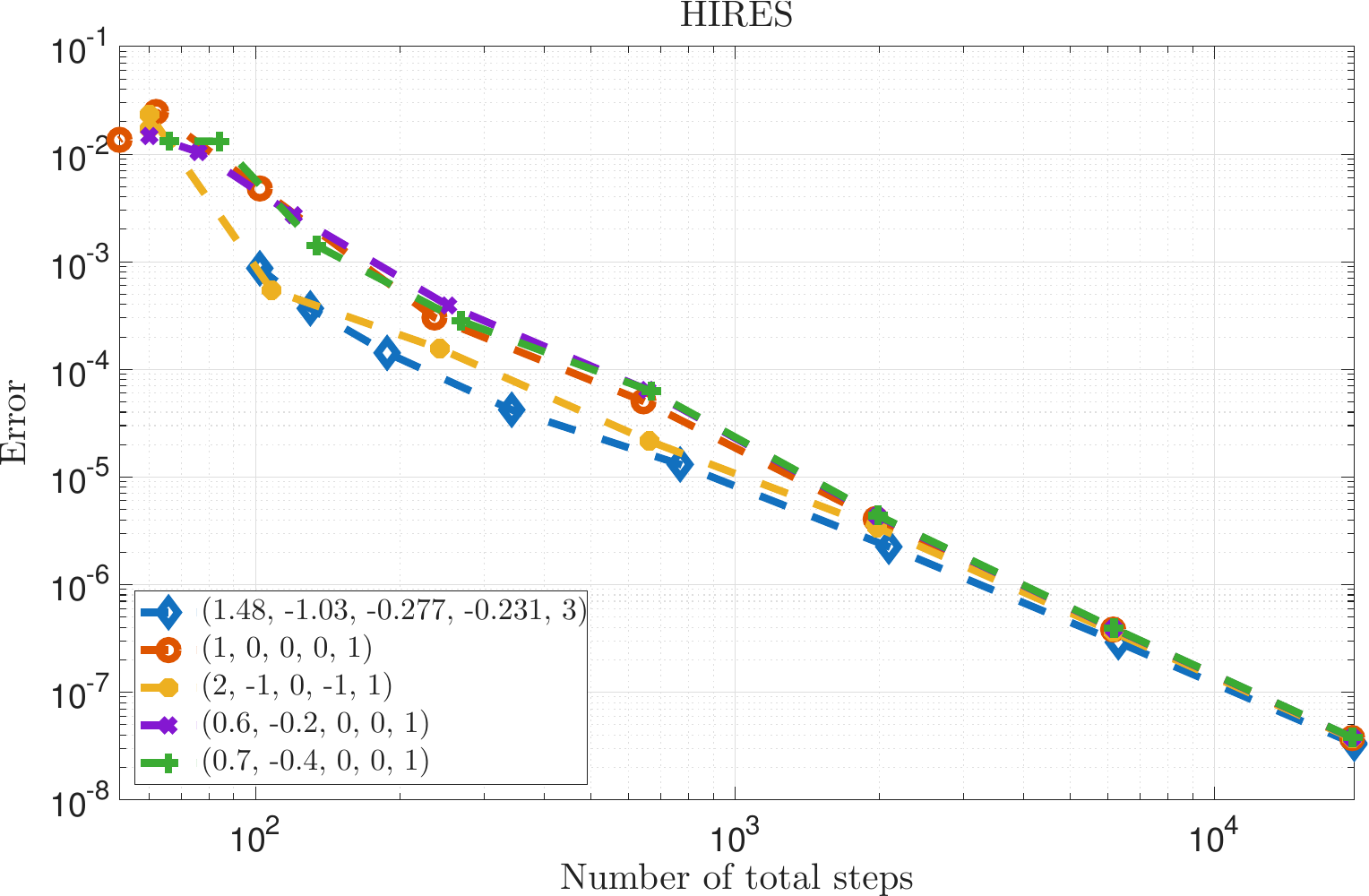}
		 	\end{subfigure}
		 	\begin{subfigure}[t]{0.5\textwidth}
		 		\includegraphics[width=\textwidth]{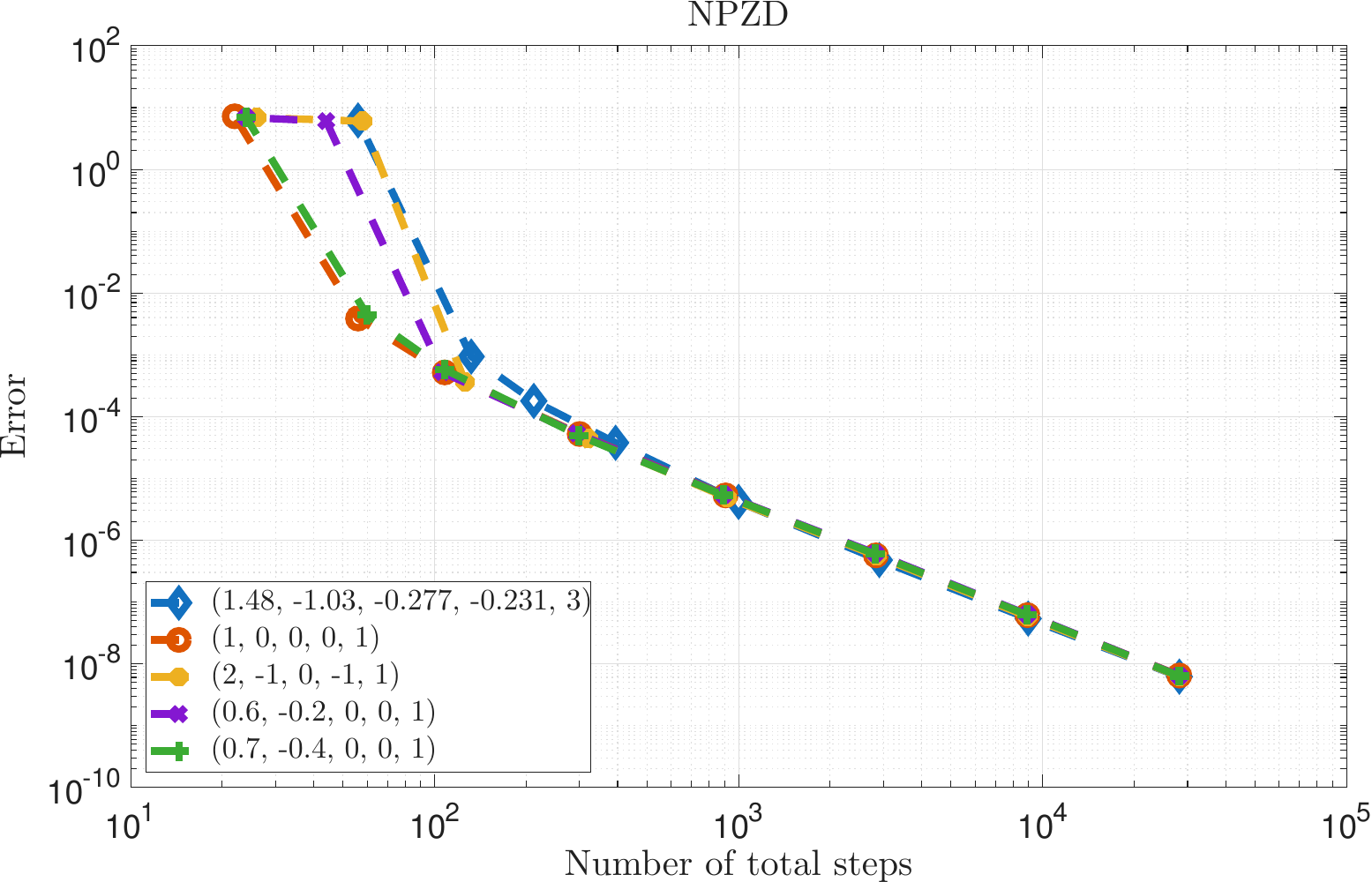}
		 	\end{subfigure}\\
		 	\begin{subfigure}[t]{0.5\textwidth}
		 		\includegraphics[width=\textwidth]{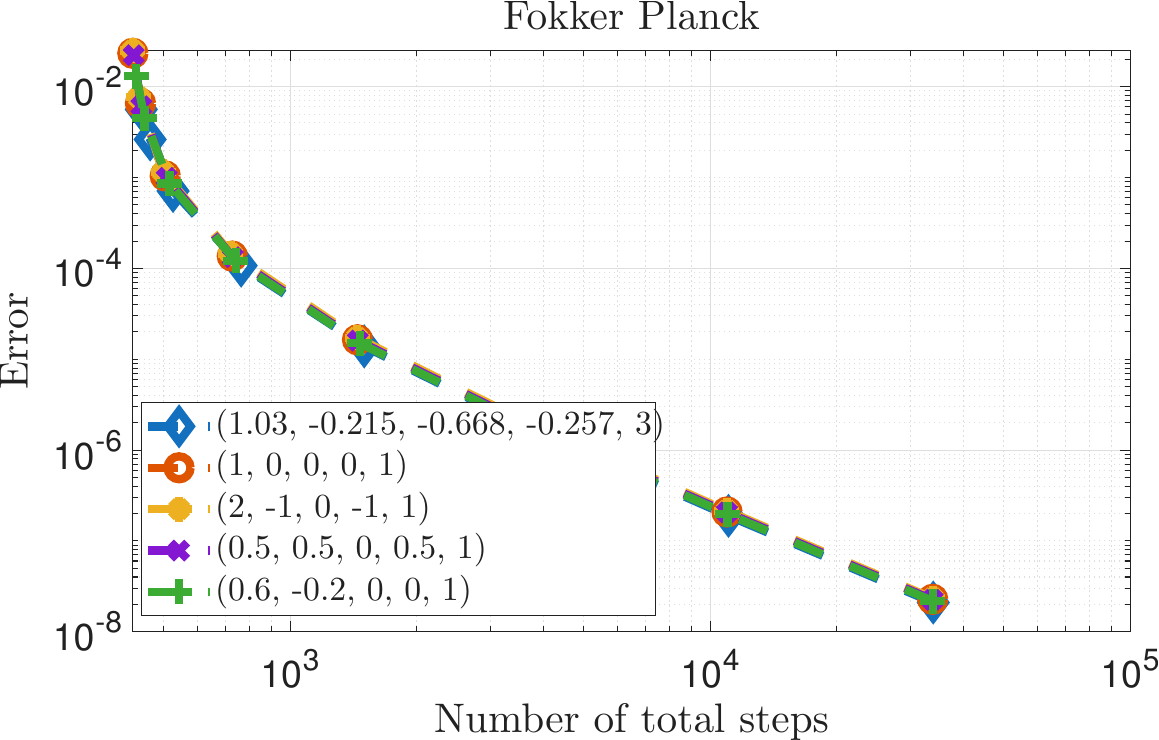}
		 	\end{subfigure}
		 	\begin{subfigure}[t]{0.5\textwidth}
		 		\includegraphics[width=\textwidth]{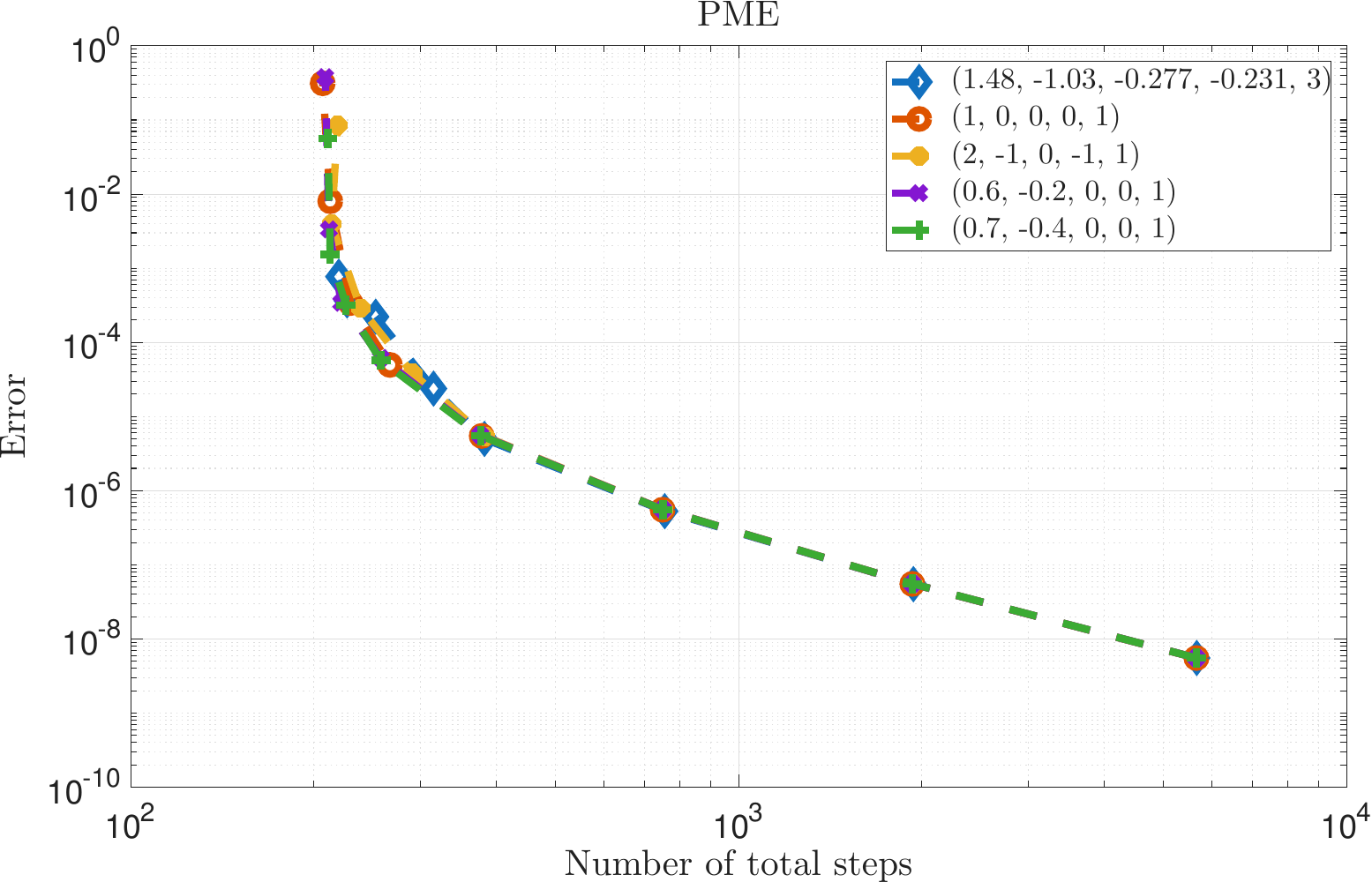}
		 	\end{subfigure}\\
		 	\begin{subfigure}[t]{0.5\textwidth}
		 		\includegraphics[width=\textwidth]{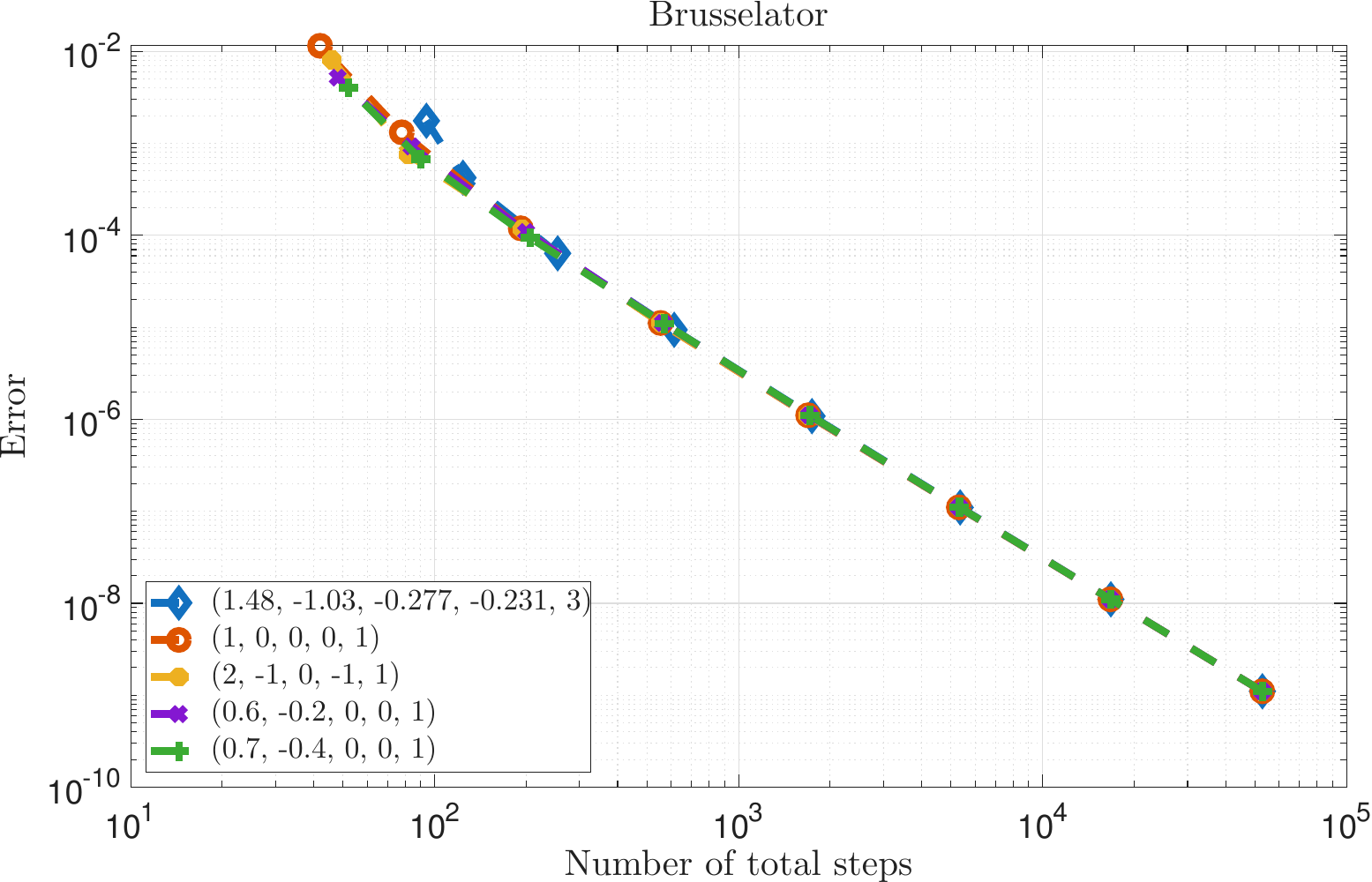}
		 	\end{subfigure}
		 	\caption{WP diagram of the problems from Section~\ref{sec:training} and Section~\ref{sec:validation} using ROS2($1/(2+\sqrt2)$) with the overall cheapest customized parameter and the best standard controllers.}\label{Fig:ROS2_val}
		 \end{figure}

		 Now let us compare the performances of the  parameters found with those of the best fitting standard controllers for the respective method. We start with the ROS2$(\gamma)$ with $\gamma=(2+\sqrt2)^{-1}$, see Figure~\ref{Fig:ROS2_val}, where we see that the standard parameters are disqualified in our procedure as they do not obey the slop criterion for the NPZD problem~\eqref{eq:NPZD} and HIRES problem \eqref{eq:hires}. In contrast, although all parameters are comparable for most of the problems, the customized parameter is consistently good and sometimes is the best in achieving the required level of accuracy, as in HIRES and PR4(0.25)~\eqref{eq:PR}.

\begin{figure}[!htbp]
	\begin{subfigure}[t]{0.5\textwidth}
		\includegraphics[width=\textwidth]{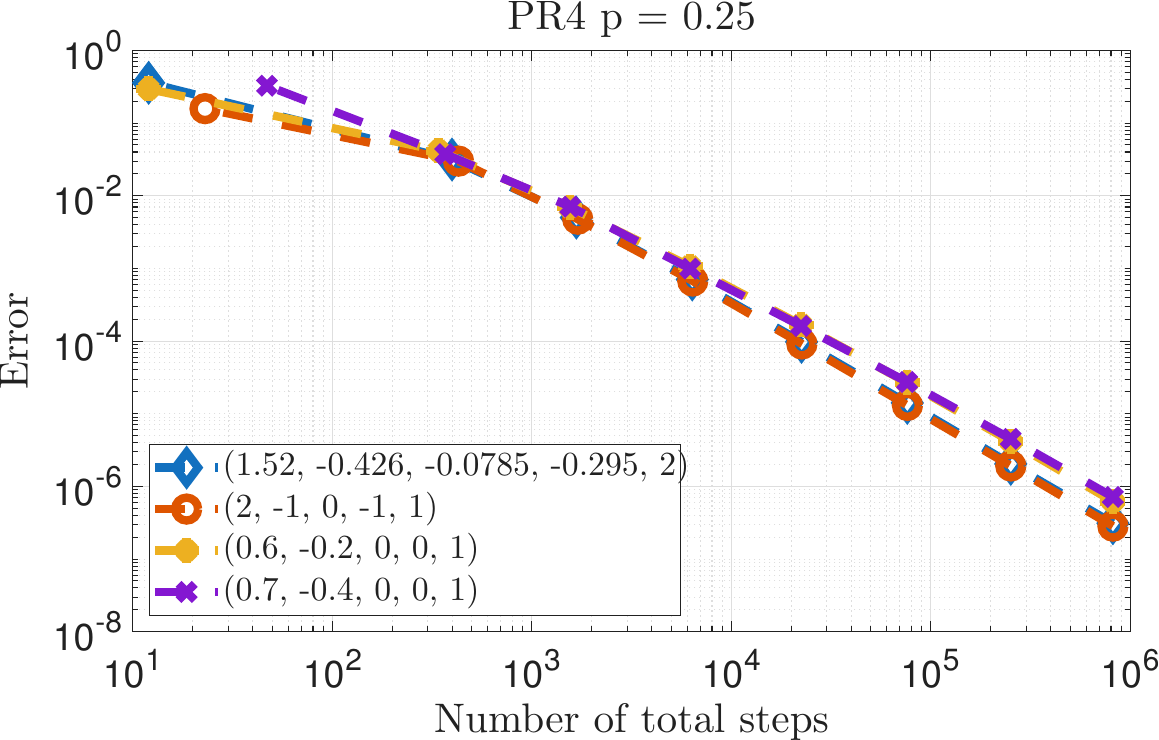}
	\end{subfigure}
	\begin{subfigure}[t]{0.5\textwidth}
		\includegraphics[width=\textwidth]{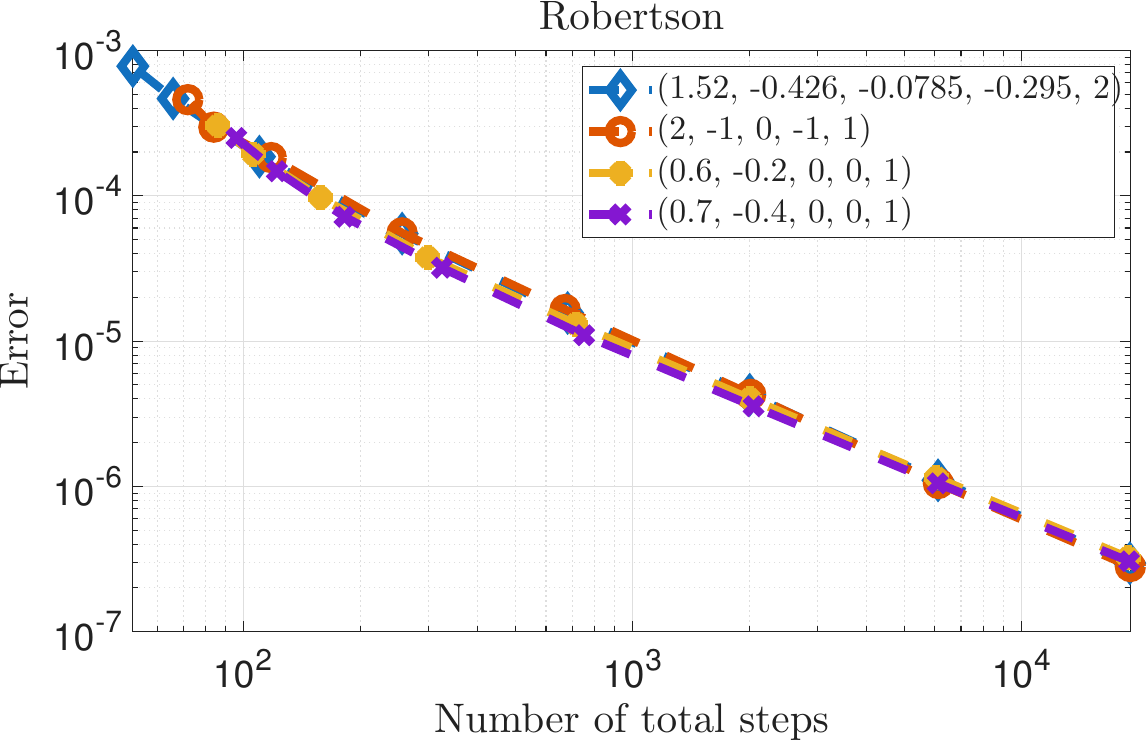}
	\end{subfigure}\\
	\begin{subfigure}[t]{0.5\textwidth}
		\includegraphics[width=\textwidth]{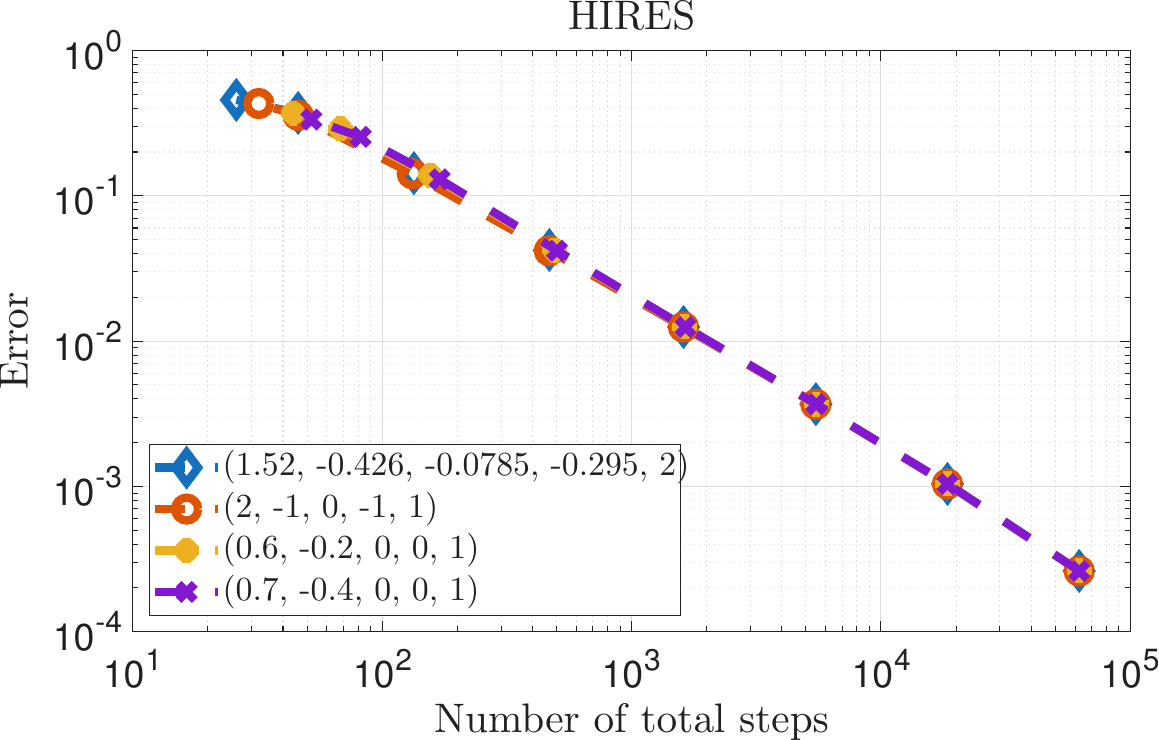}
	\end{subfigure}
	\begin{subfigure}[t]{0.5\textwidth}
		\includegraphics[width=\textwidth]{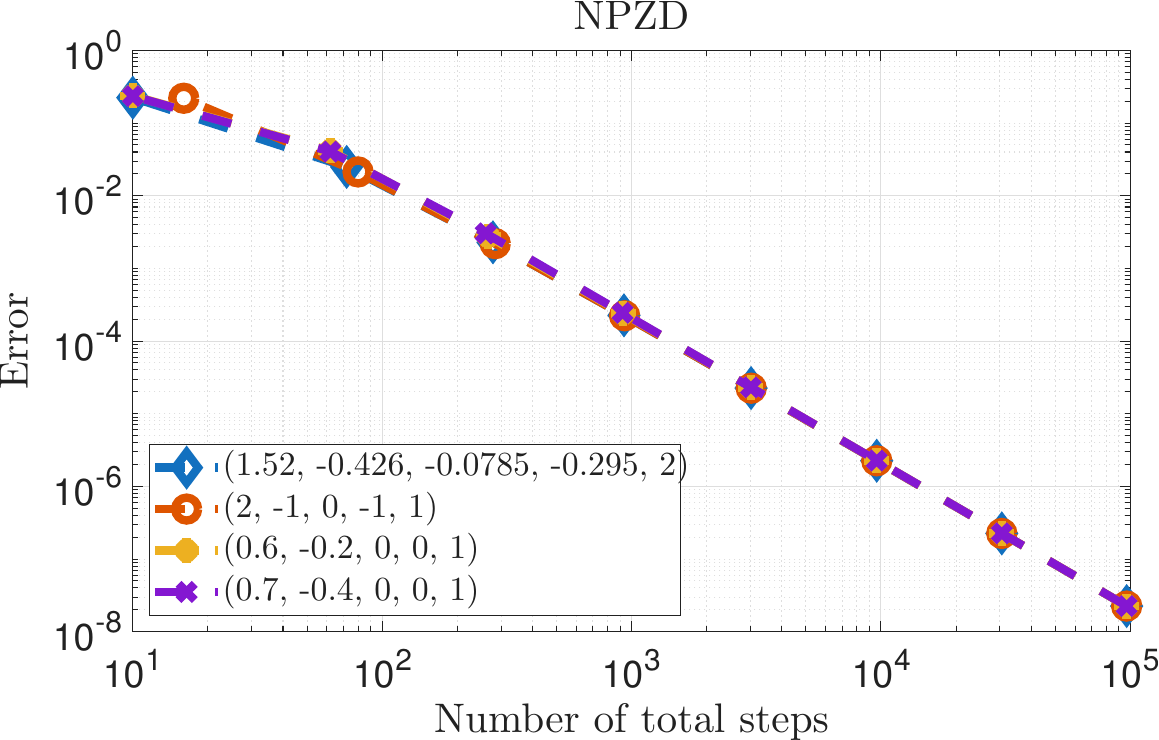}
	\end{subfigure}\\
	\begin{subfigure}[t]{0.5\textwidth}
		\includegraphics[width=\textwidth]{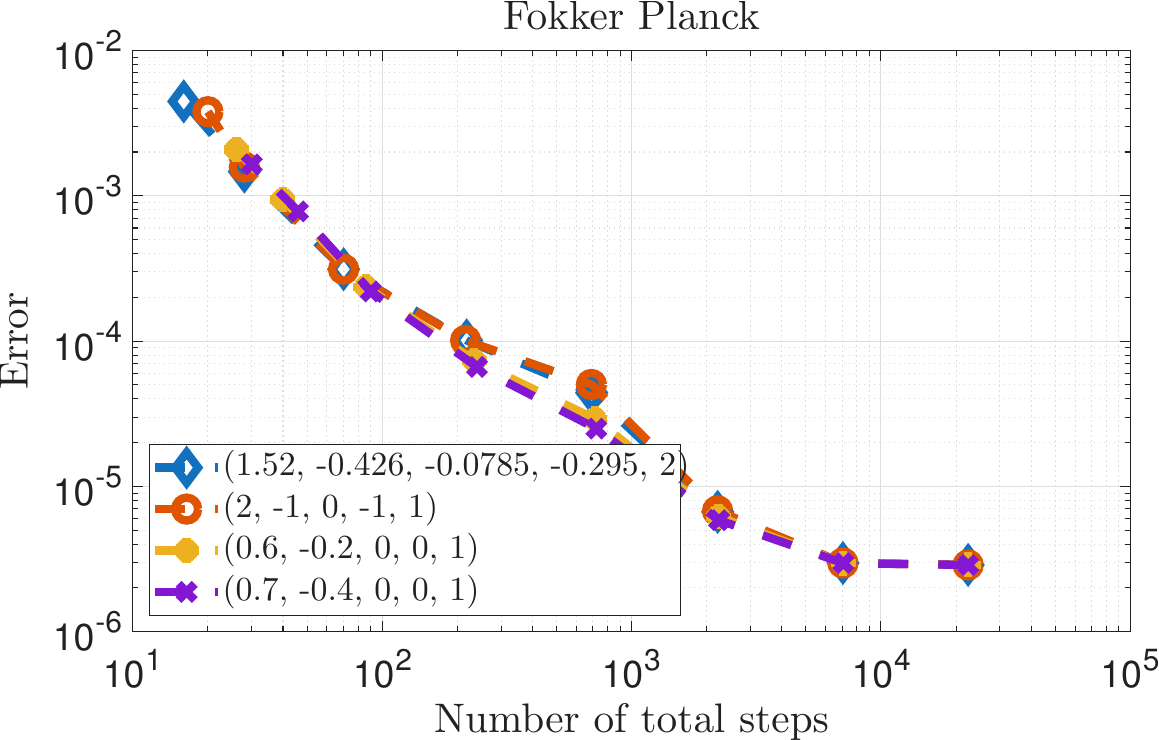}
	\end{subfigure}
	\begin{subfigure}[t]{0.5\textwidth}
		\includegraphics[width=\textwidth]{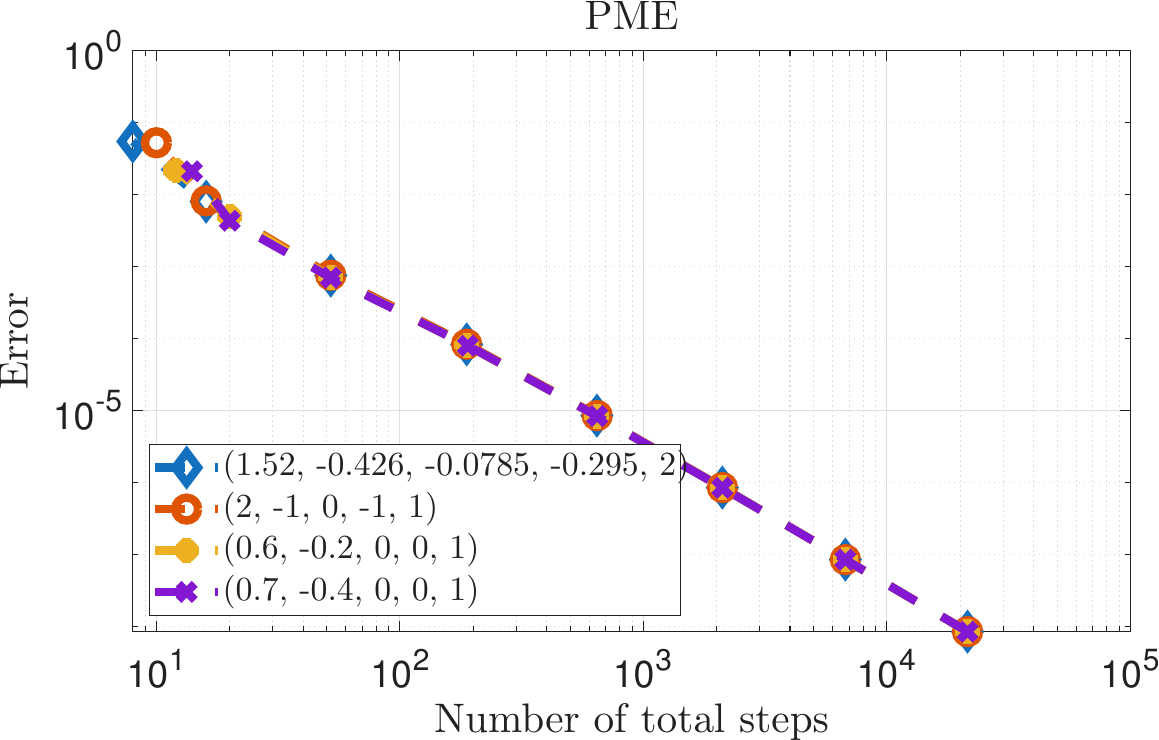}
	\end{subfigure}\\
	\begin{subfigure}[t]{0.5\textwidth}
		\includegraphics[width=\textwidth]{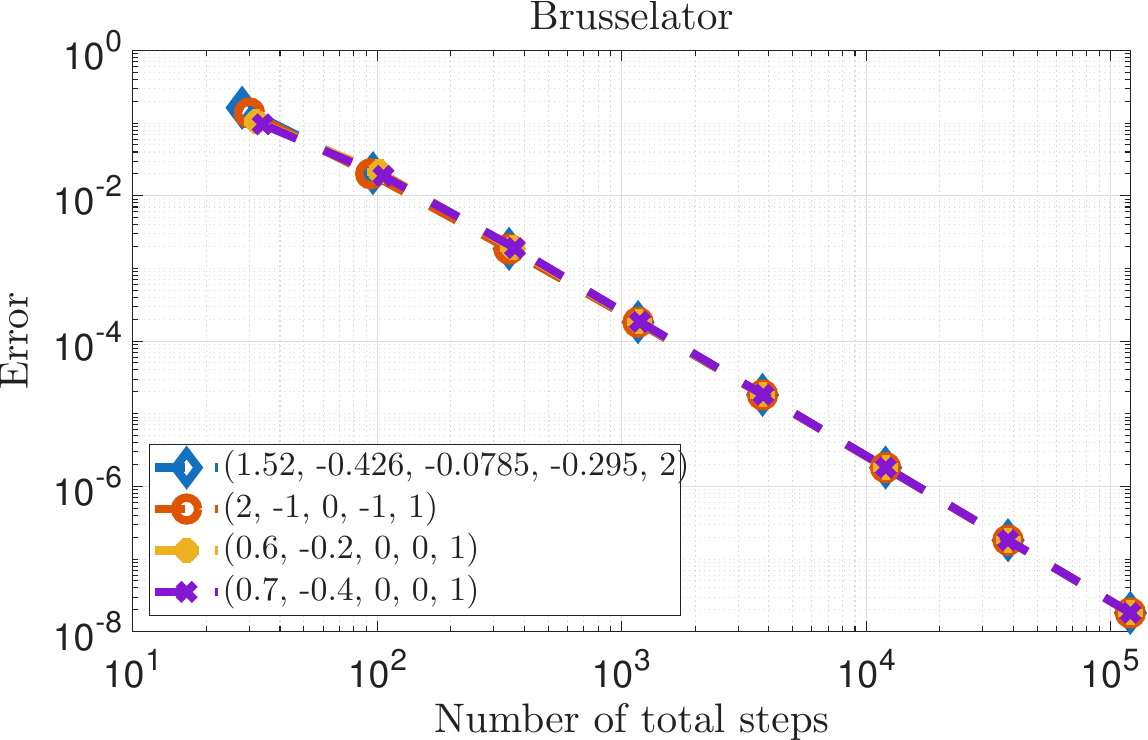}
	\end{subfigure}
	\caption{WP diagram of the problems from Section~\ref{sec:training} and Section~\ref{sec:validation} using MPRK22(1) with the overall cheapest customized parameter and the best standard controllers.}\label{Fig:MPRK22_val}
\end{figure}

		 Next, looking at the second-order MPRK scheme, see Figure~\ref{Fig:MPRK22_val}. Here, it becomes evident that the new controller performs slightly better than the best standard controllers as, besides being highly accurate for small tolerances, it is also cheaper for the larger ones. The only exception is found in the WP diagram for the Fokker--Planck equation \eqref{eq:FP}, where $p_1$ and $p_2$ work better for middle-ranged tolerances.

\begin{figure}[!htbp]
	\begin{subfigure}[t]{0.5\textwidth}
	\includegraphics[width=\textwidth]{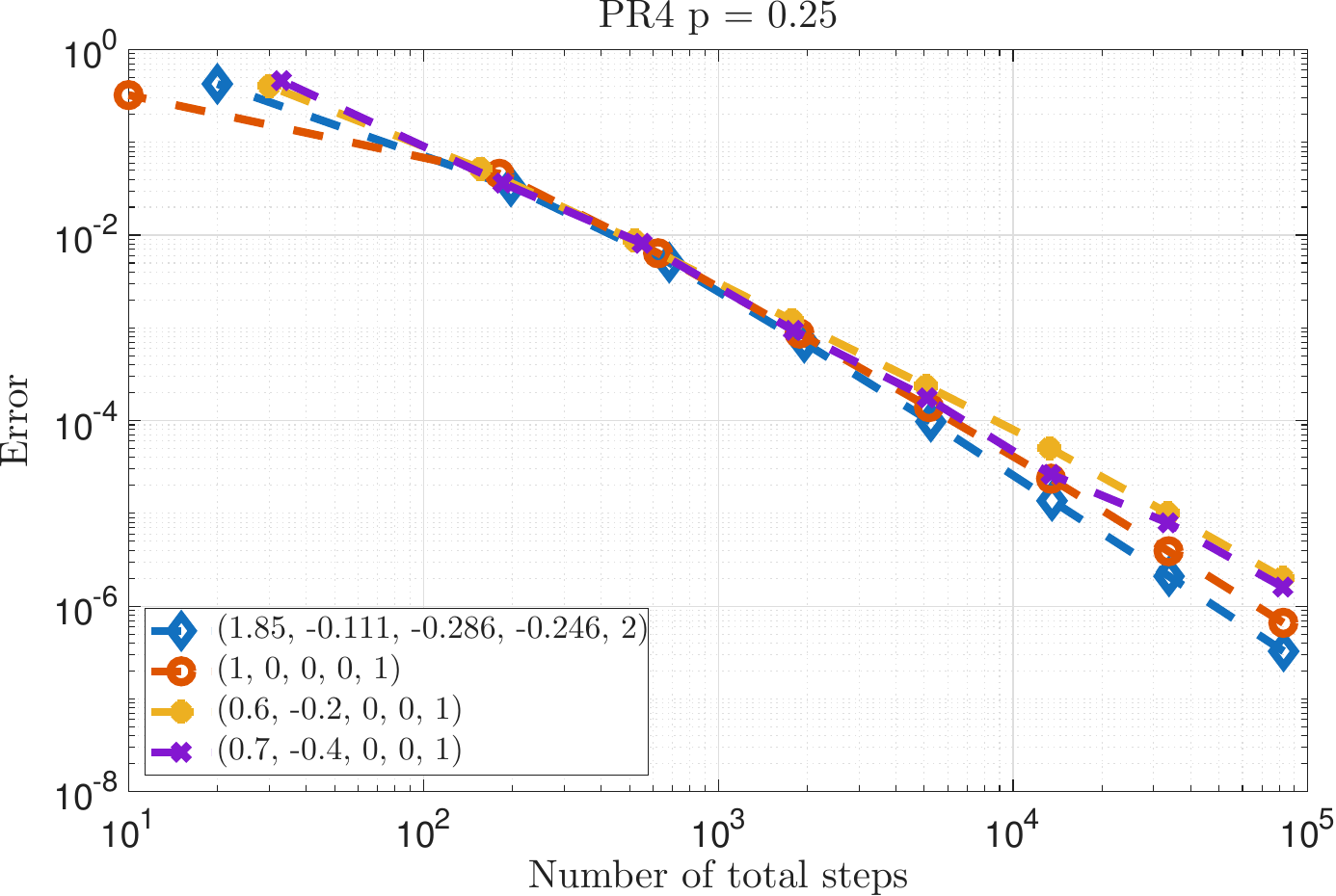}
\end{subfigure}
\begin{subfigure}[t]{0.5\textwidth}
	\includegraphics[width=\textwidth]{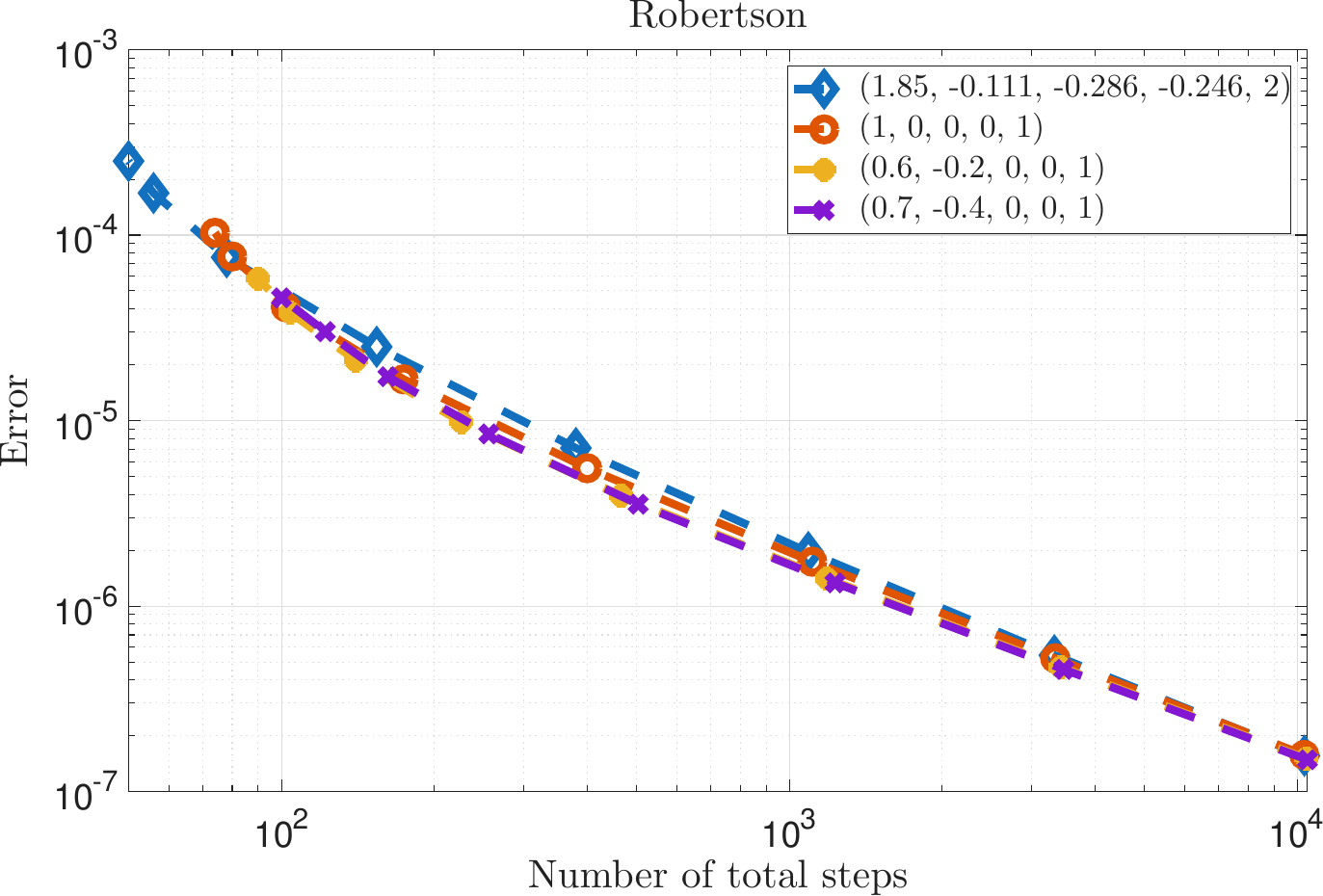}
\end{subfigure}\\
\begin{subfigure}[t]{0.5\textwidth}
	\includegraphics[width=\textwidth]{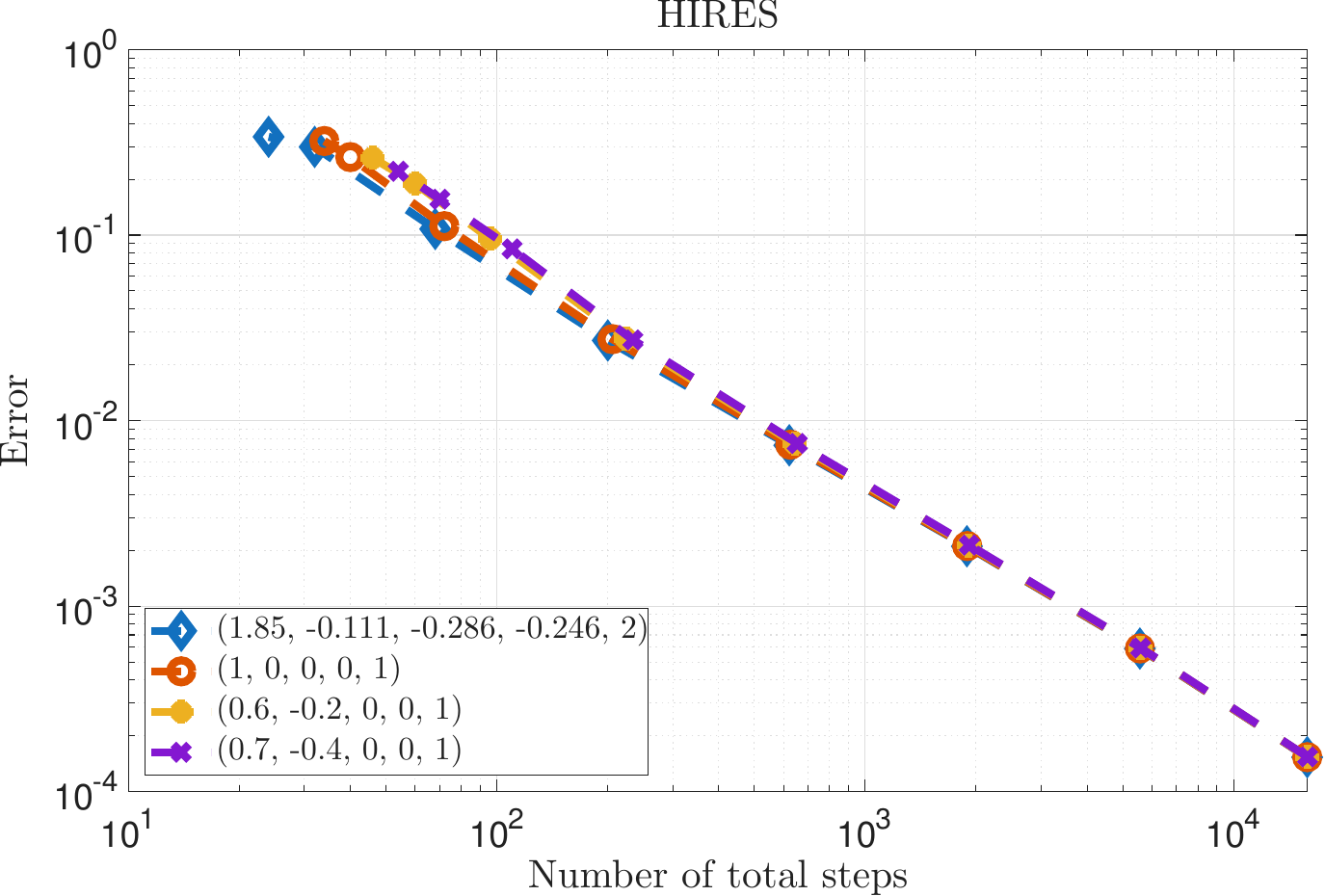}
\end{subfigure}
\begin{subfigure}[t]{0.5\textwidth}
	\includegraphics[width=\textwidth]{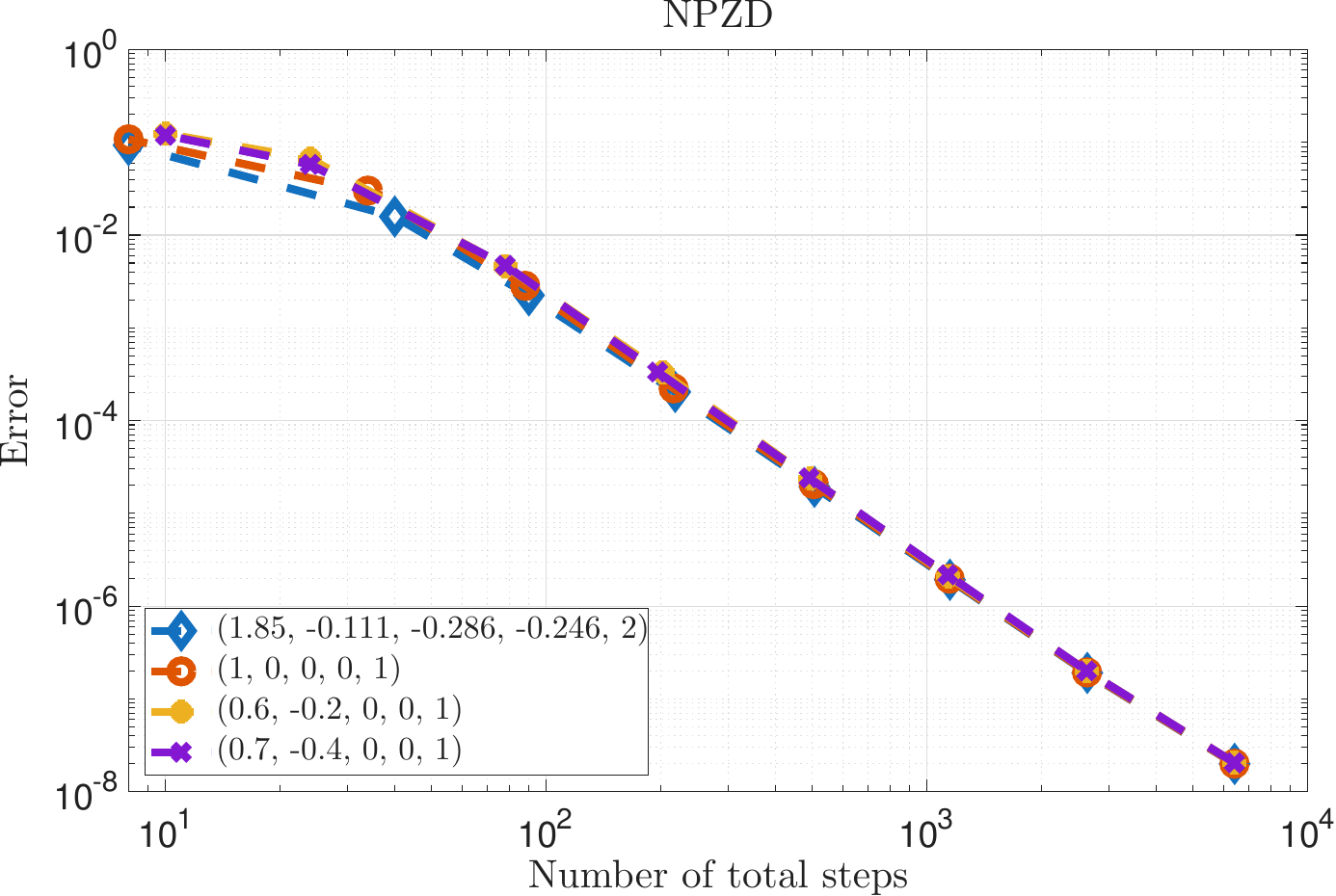}
\end{subfigure}\\
\begin{subfigure}[t]{0.5\textwidth}
	\includegraphics[width=\textwidth]{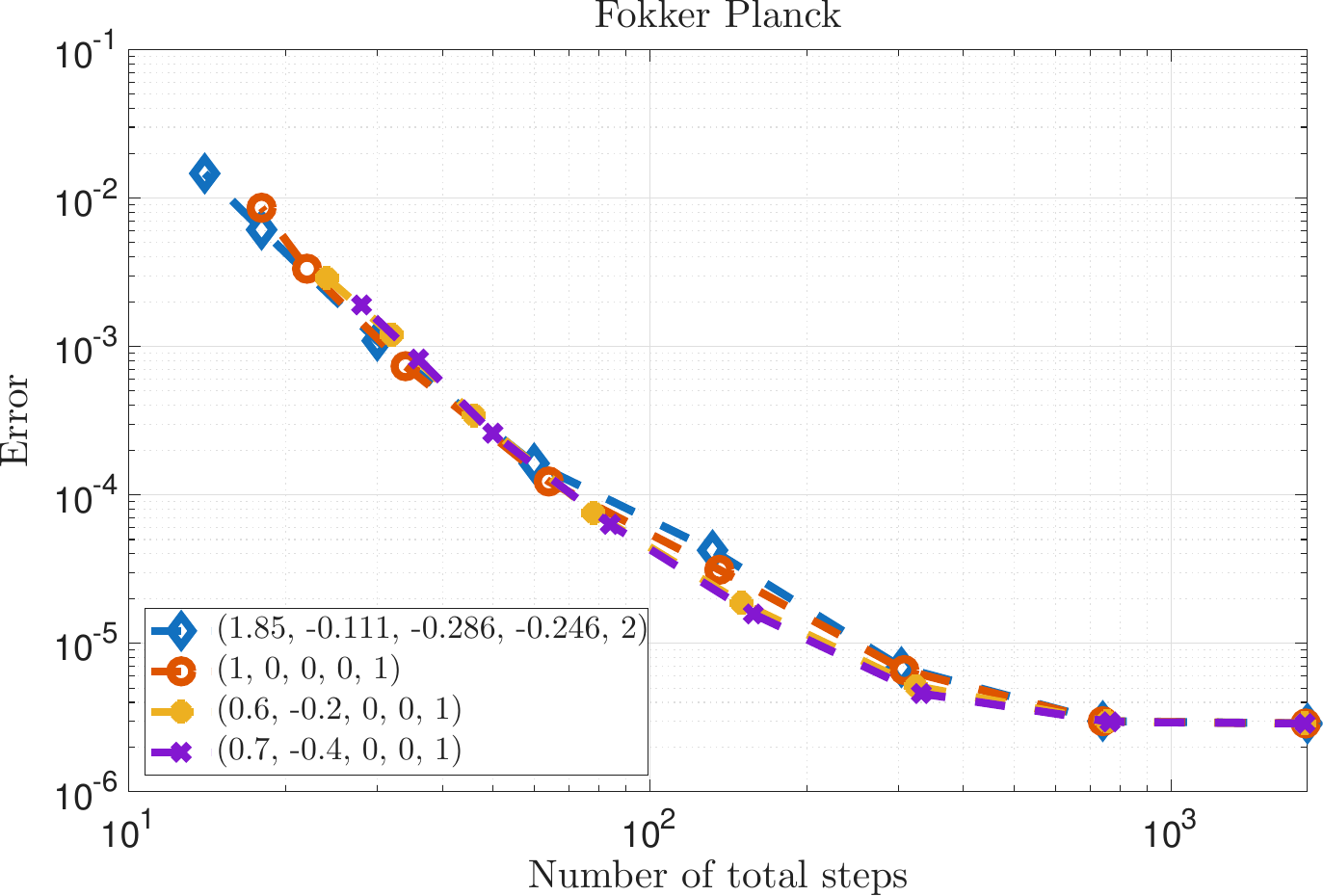}
\end{subfigure}
\begin{subfigure}[t]{0.5\textwidth}
	\includegraphics[width=\textwidth]{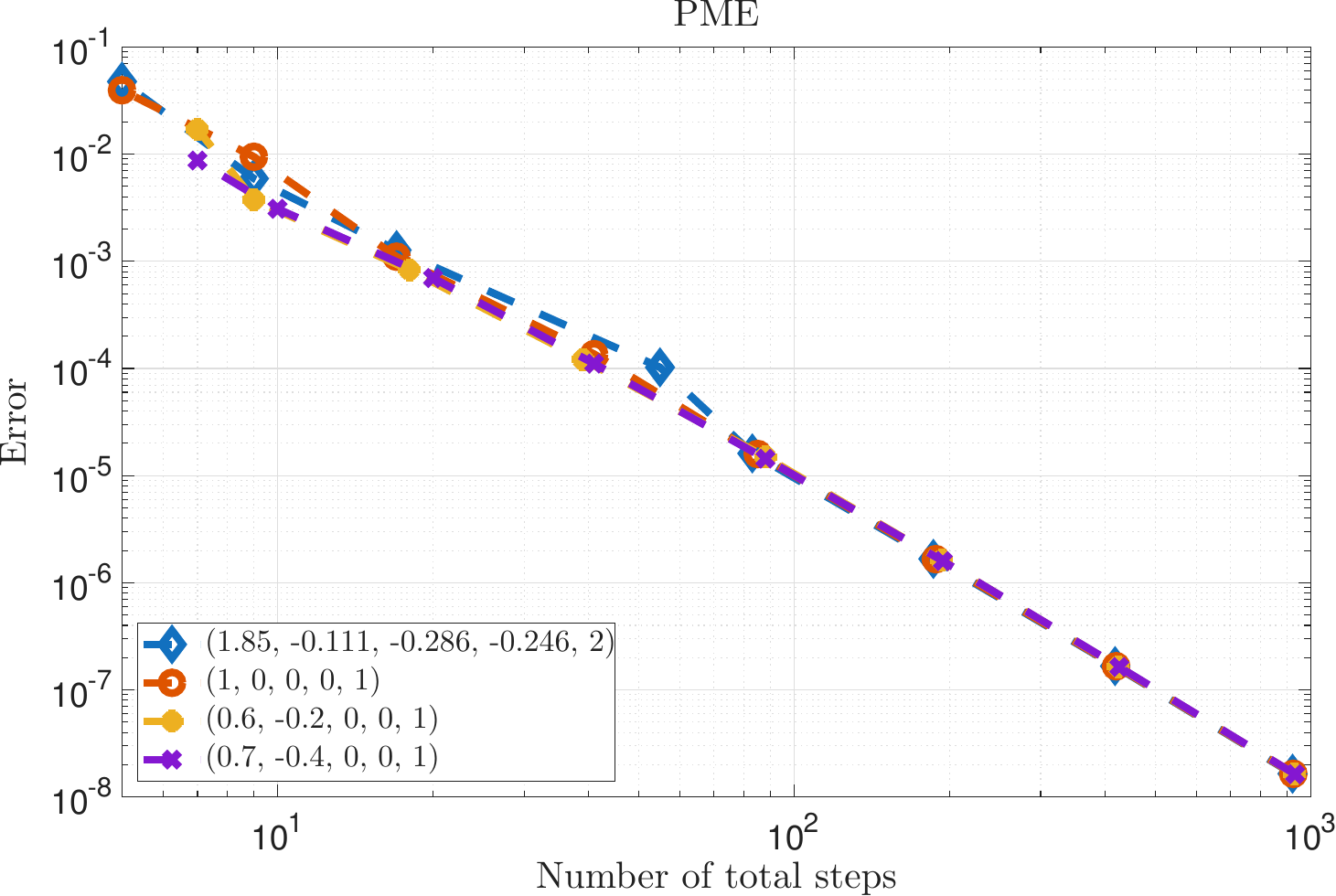}
\end{subfigure}\\
\begin{subfigure}[t]{0.5\textwidth}
	\includegraphics[width=\textwidth]{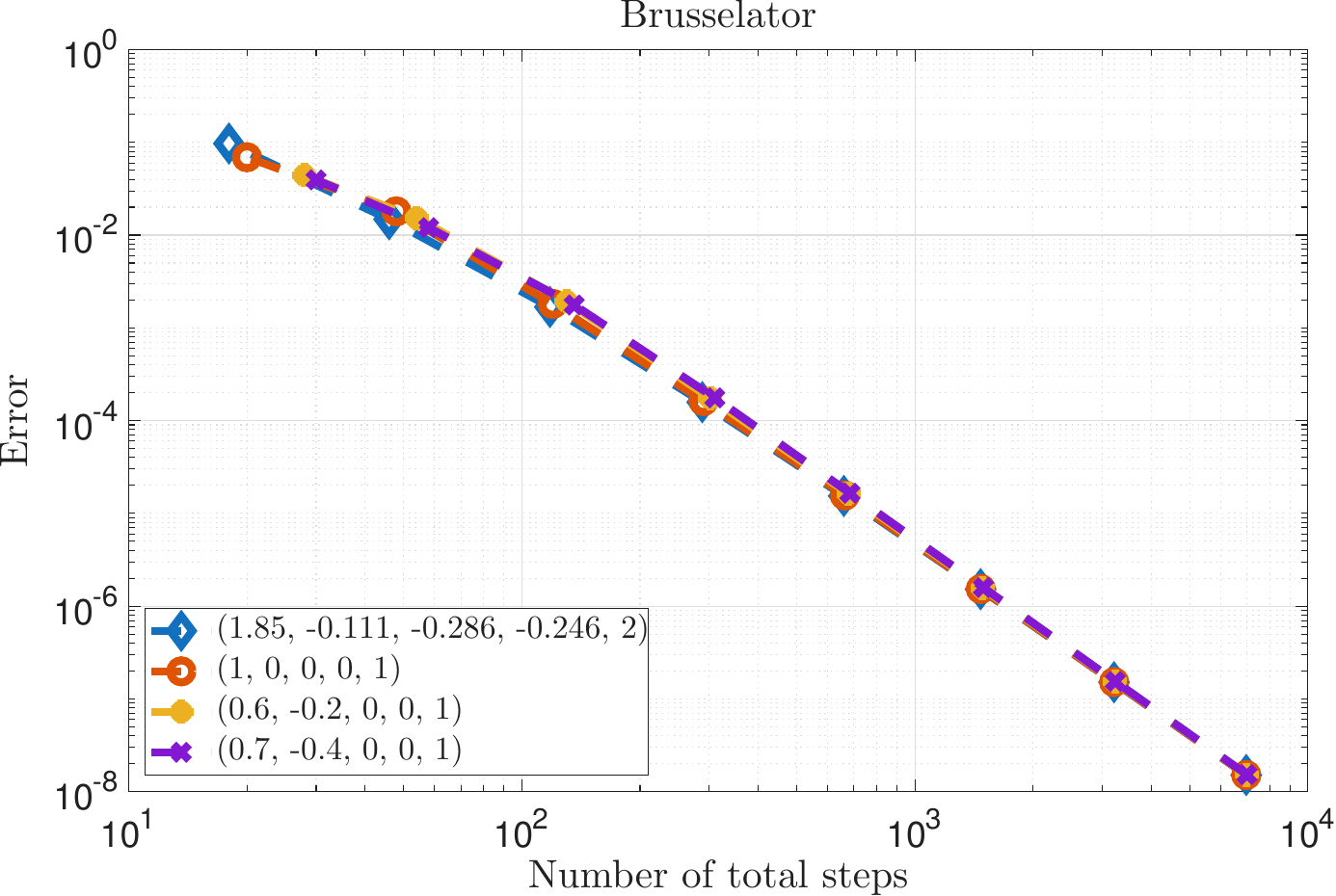}
\end{subfigure}
	\caption{WP diagram of the problems from Section~\ref{sec:training} and Section~\ref{sec:validation} using MPRK43(0.5,0.75) with the overall cheapest customized parameter and the best standard controllers $p_1$, $p_2$ and $p_3$.}\label{Fig:MPRK43I_val}
\end{figure}

\begin{figure}[!htbp]
		\begin{subfigure}[t]{0.5\textwidth}
		\includegraphics[width=\textwidth]{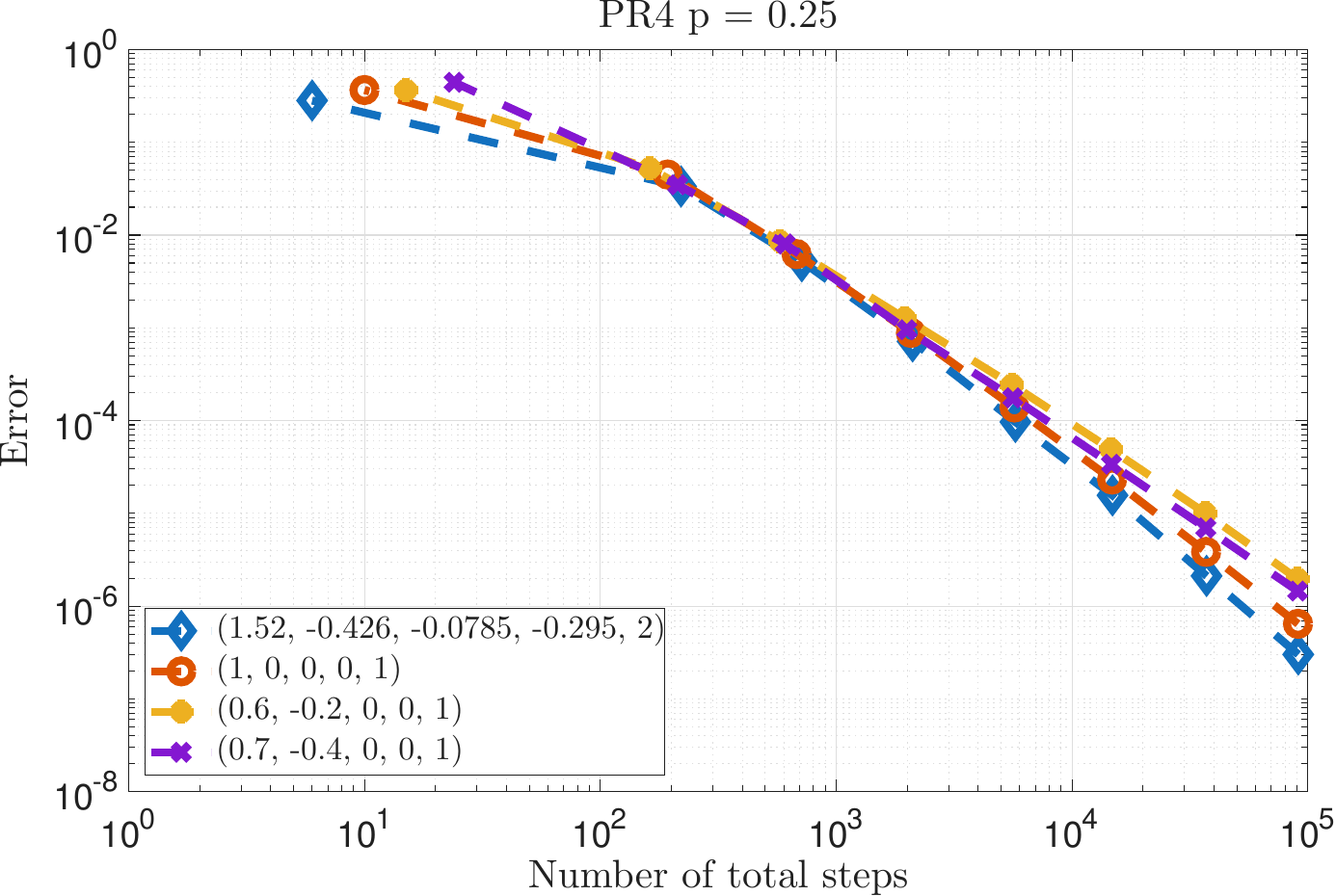}
	\end{subfigure}
	\begin{subfigure}[t]{0.5\textwidth}
		\includegraphics[width=\textwidth]{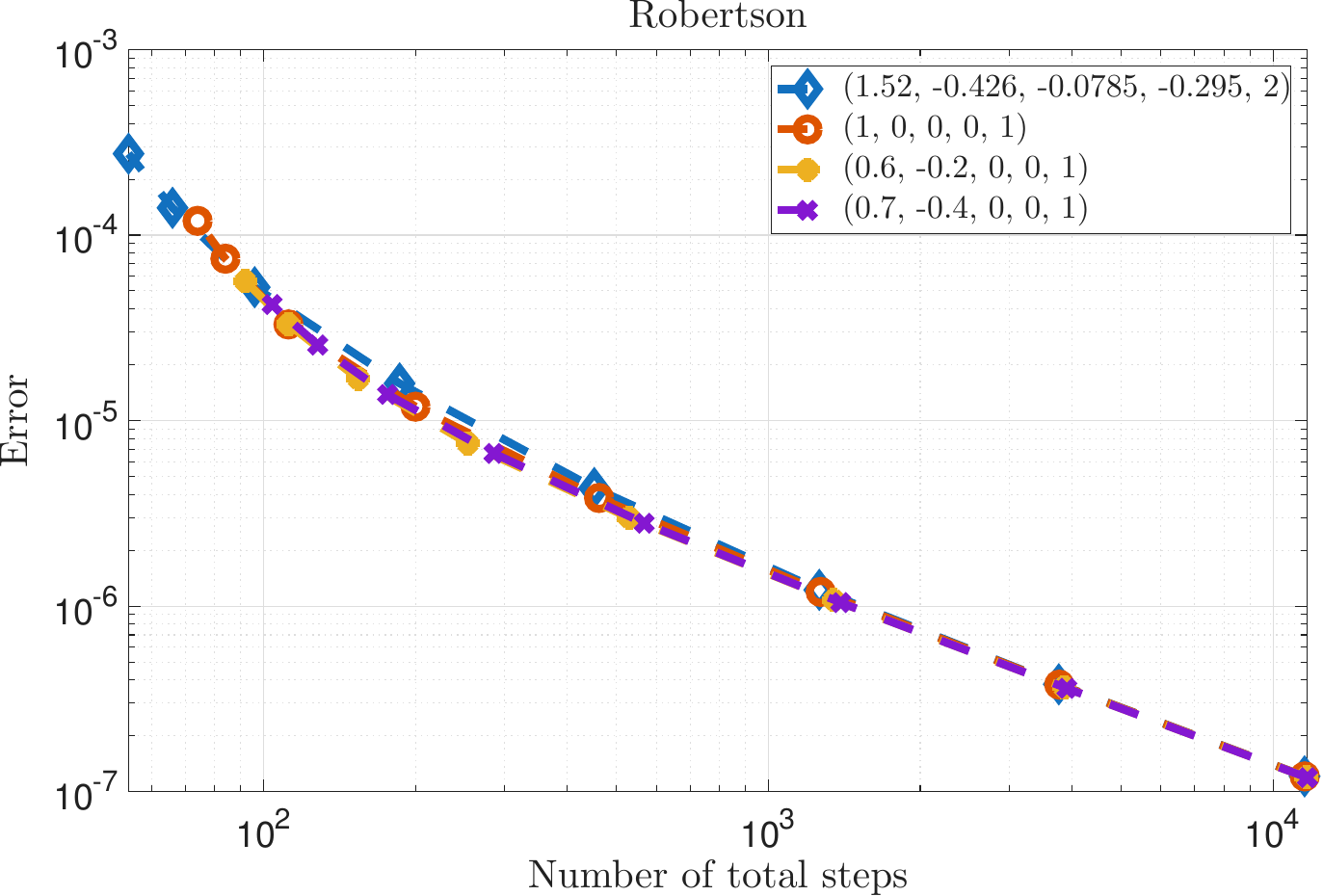}
	\end{subfigure}\\
	\begin{subfigure}[t]{0.5\textwidth}
		\includegraphics[width=\textwidth]{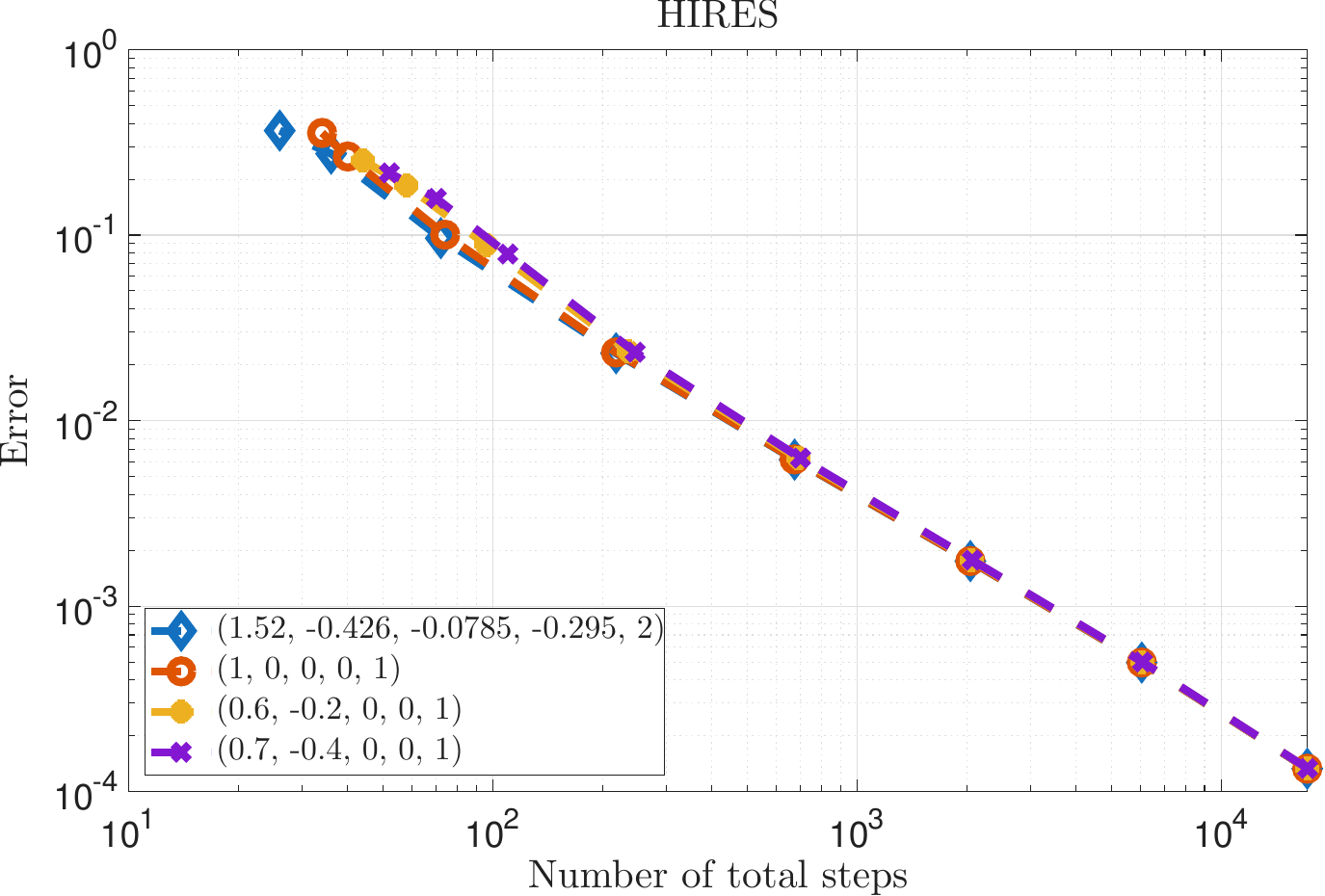}
	\end{subfigure}
	\begin{subfigure}[t]{0.5\textwidth}
		\includegraphics[width=\textwidth]{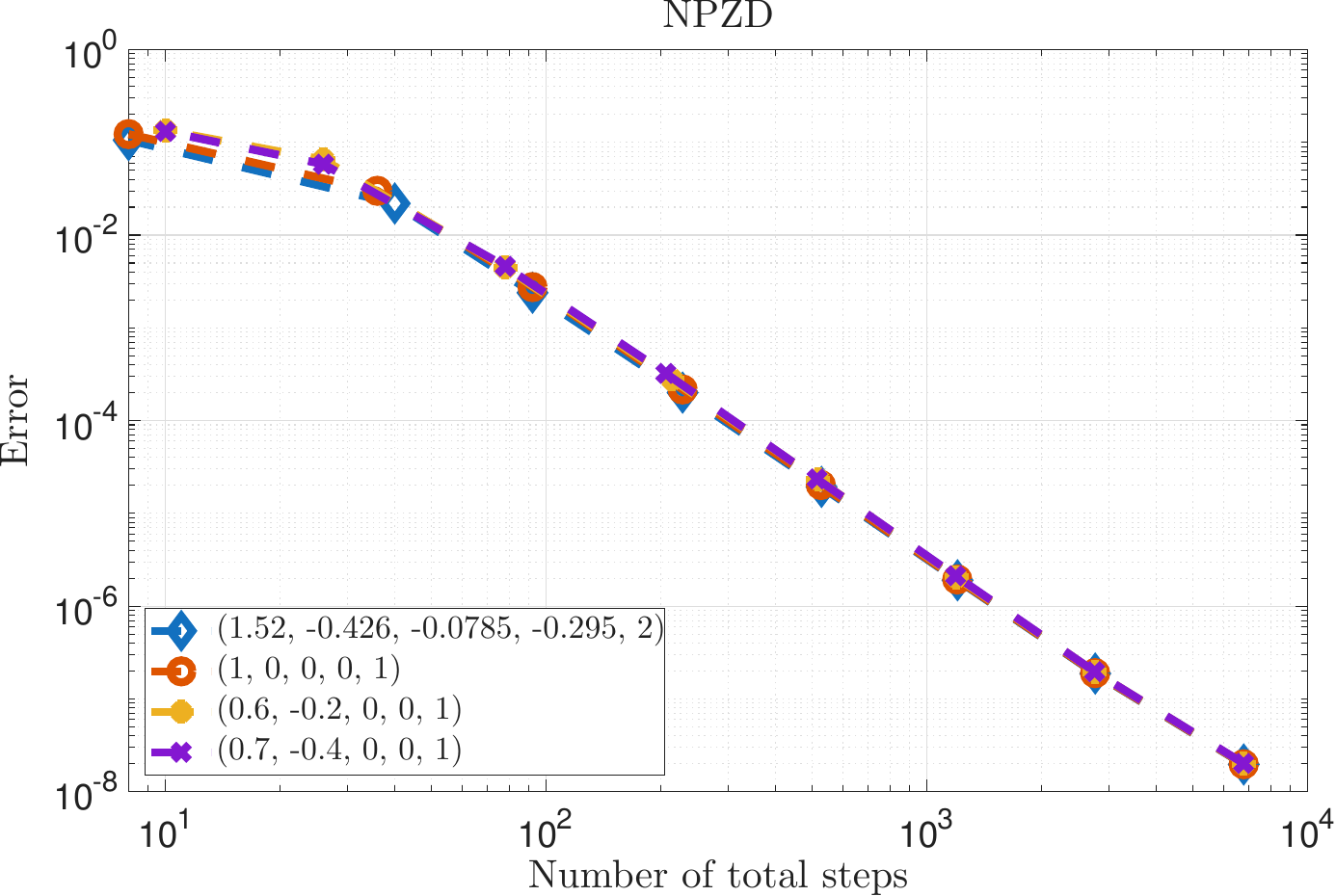}
	\end{subfigure}\\
	\begin{subfigure}[t]{0.5\textwidth}
		\includegraphics[width=\textwidth]{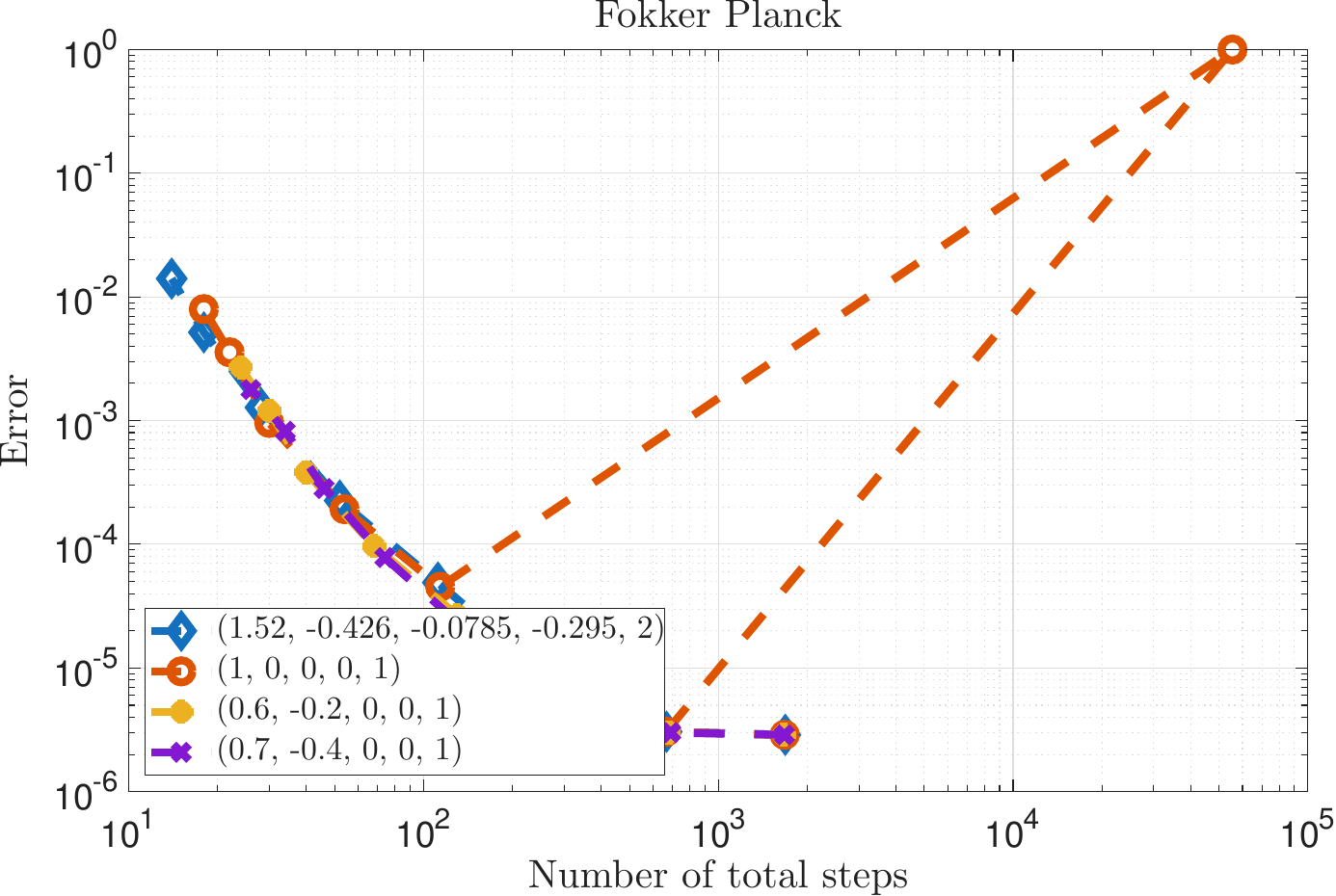}
	\end{subfigure}
	\begin{subfigure}[t]{0.5\textwidth}
		\includegraphics[width=\textwidth]{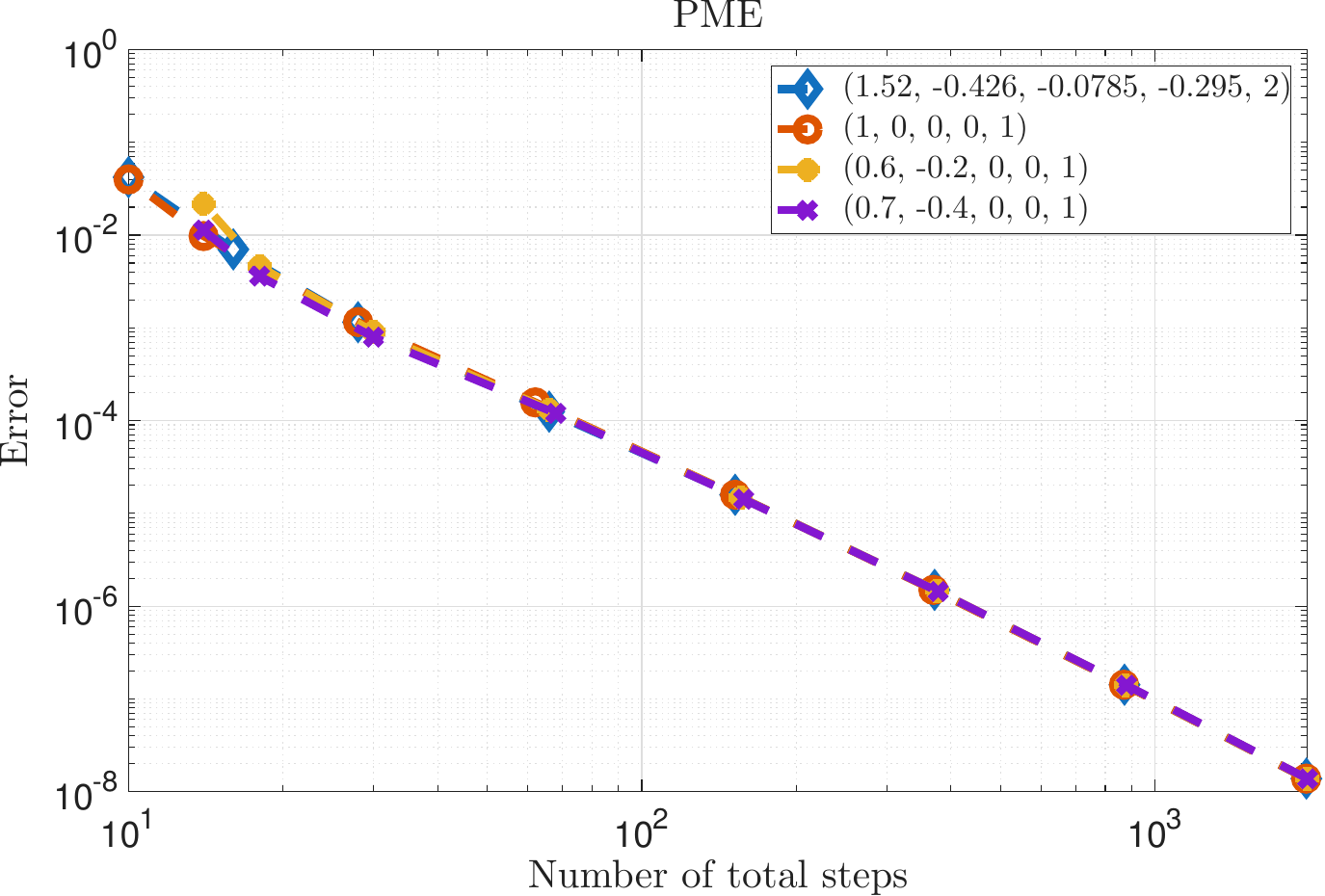}
	\end{subfigure}\\
	\begin{subfigure}[t]{0.5\textwidth}
		\includegraphics[width=\textwidth]{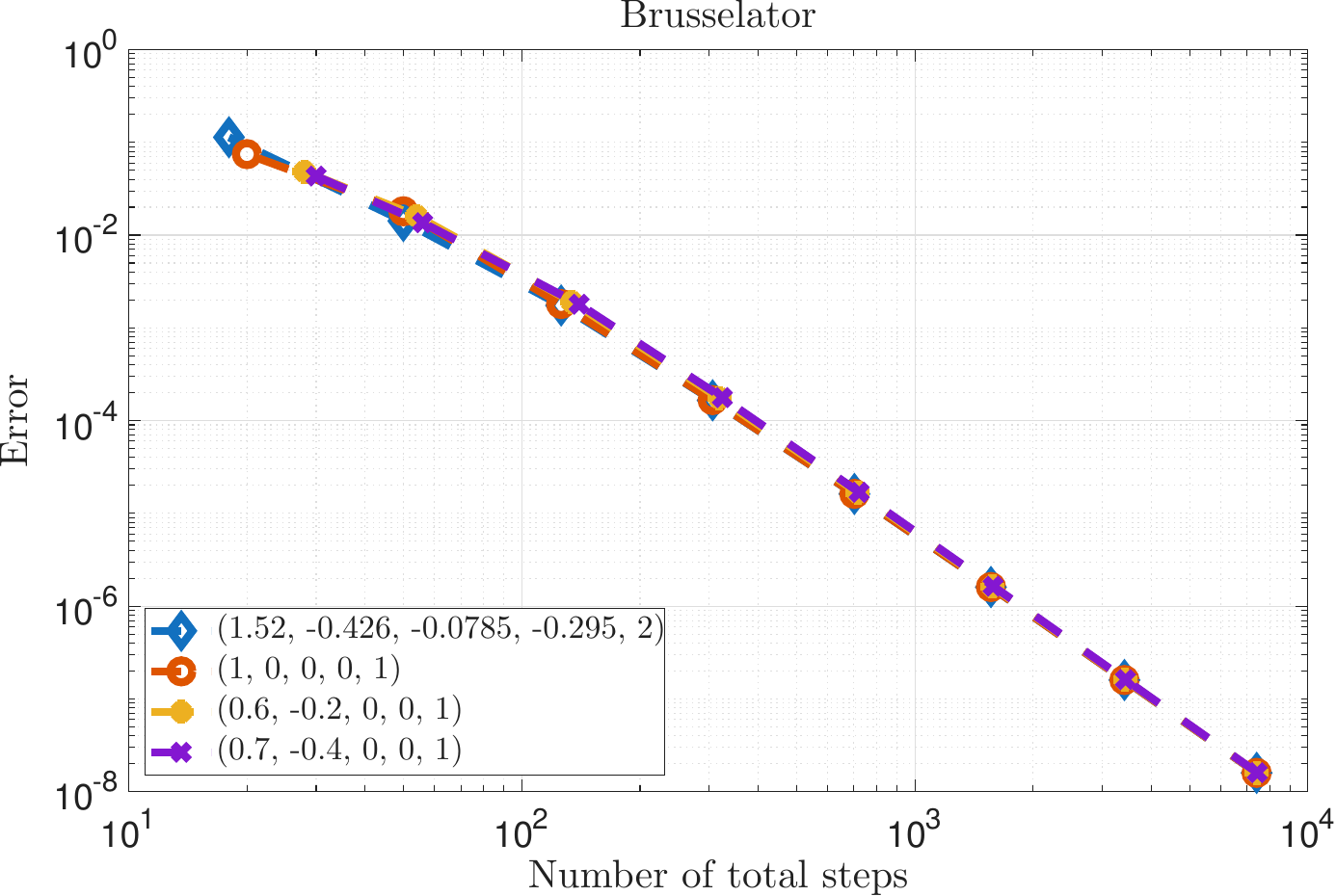}
	\end{subfigure}
	\caption{WP diagram of the problems from Section~\ref{sec:training} and Section~\ref{sec:validation} using MPRK43(0.563) with the overall cheapest customized parameter and the best standard controllers $p_1$, $p_2$ and $p_3$.  We note here, that the error was set to 1 whenever our abortion criteria from Section~\ref{sec:methodology} were met -- except for the slope criterion.}\label{Fig:MPRK43II_val}
\end{figure}

This becomes even clearer when looking at the third-order schemes, see Figure~\ref{Fig:MPRK43I_val} and Figure~\ref{Fig:MPRK43II_val}, respectively. Even though the found parameter is slightly less accurate for mid-ranged tolerances, the found controller performs clearly better at the largest and smallest tolerances. Also note that in the case of the Robertson problem the calculated error is below the given tolerance for \mbox{$\code{\tol}\in \{ 10^{-j}\mid j=1,\dotsc,5\}$}, \ie is meeting our requirements of a good controller. We note here, that the error was set to 1 whenever our abortion criteria from Section~\ref{sec:methodology} were met -- except for the slope criterion.

				\subsection{Comparison with Adaptive Runge--Kutta Methods}
				We also compare the performance of the methods with parameters presented in Table~\ref{tab:optiparas2} with some standard time integration methods.
				A comparison of MPRK22(1) with constant step size to built-in solvers in \MATLAB for the NPZD problem \eqref{eq:NPZD} was already done in \cite{KMpos}, where MPRK schemes were more efficient for coarser tolerances.
				In this work, we compare the MPRK schemes with adaptive time step control to the built-in solvers. In particular we compare the methods with the second-order, linearly implicit Rosenbrock scheme \code{ode23s} and the variable-order BDF method \code{ode15s} \cite{shampine1997matlab}. In both cases, we use \code{options.NonNegative = 1} in order to compare the positivity-preserving variants of the built-in solvers with our schemes.
				We compare the number of right-hand side (RHS) evaluations and the total runtime. For the latter we note that we use our proof of concept implementation which means that the comparison may be unfair given that these built-in functions are optimized in \MATLAB.
				The elapsed time is measured using the built-in functions \code{tic} and \code{toc} from \MATLAB. Furthermore, we averaged the elapsed time using ten runs on a Dell Precision 7680 notebook with an 13th Gen Intel(R) Core(TM) i7-13850HX processor and \MATLAB R2025b.

\begin{figure}[!htbp]
					\begin{subfigure}[t]{0.495\textwidth}
						\includegraphics[width=\textwidth]{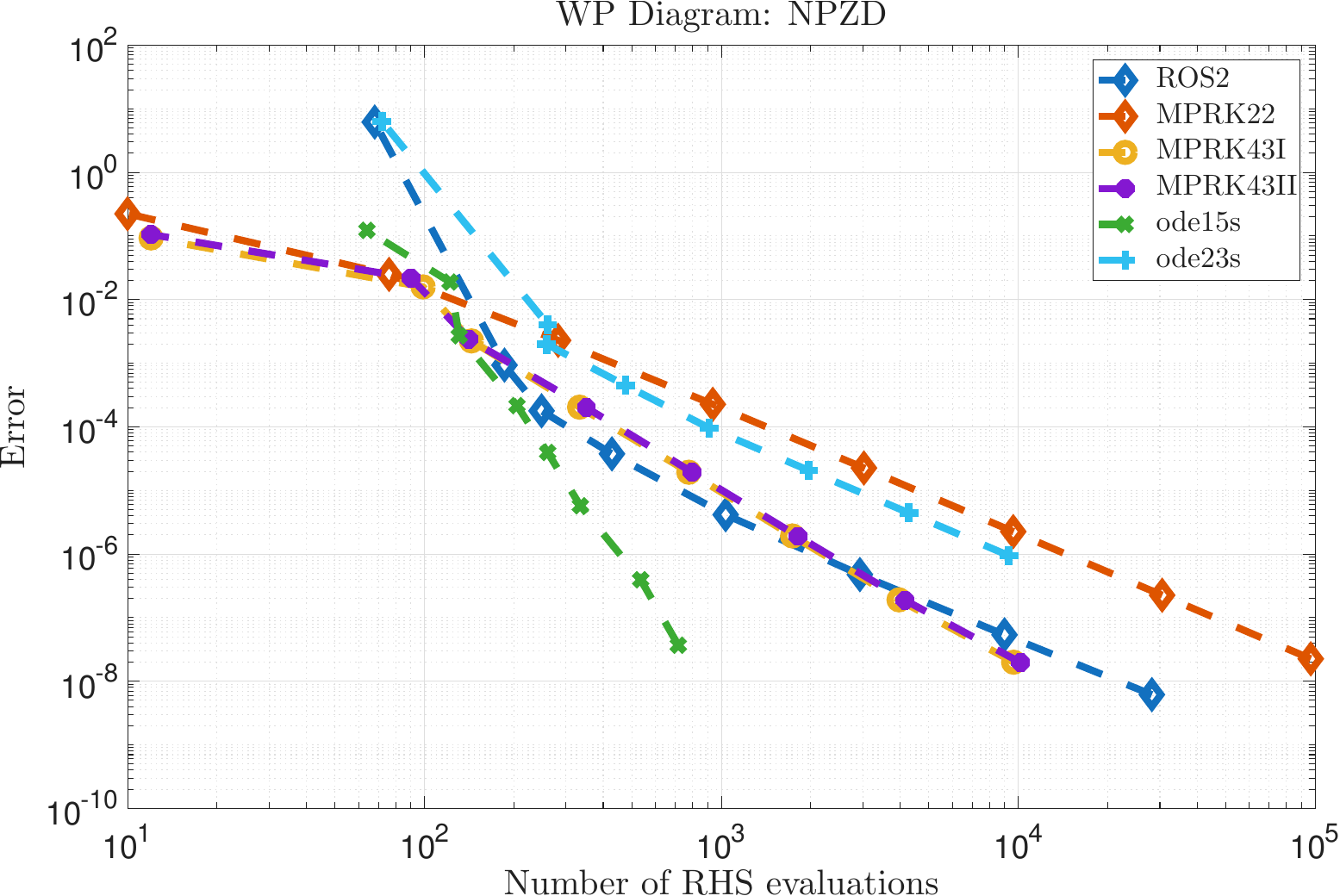}
					\end{subfigure}
					\begin{subfigure}[t]{0.495\textwidth}
						\includegraphics[width=\textwidth]{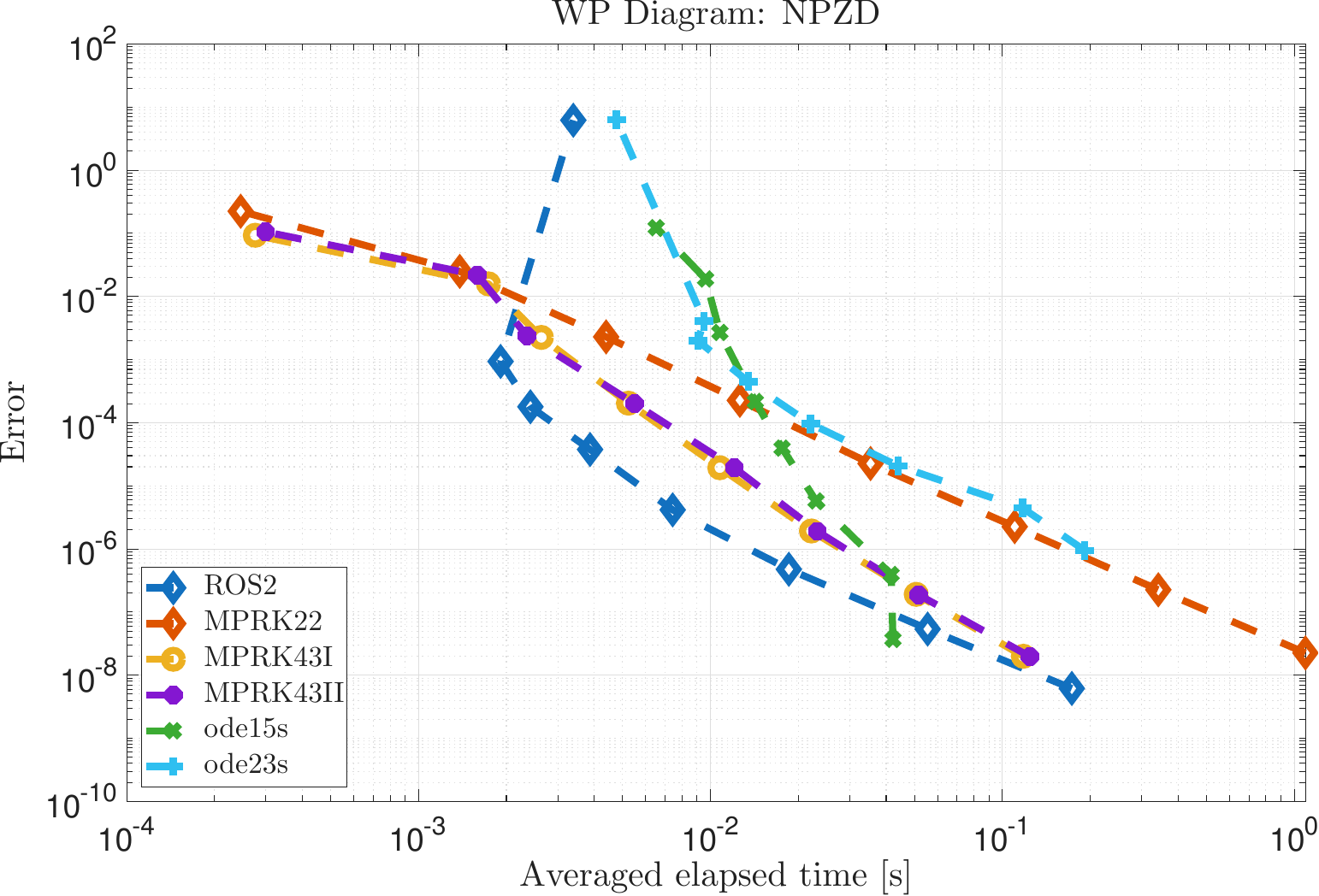}
					\end{subfigure}\\
									\begin{subfigure}[t]{0.495\textwidth}
					\includegraphics[width=\textwidth]{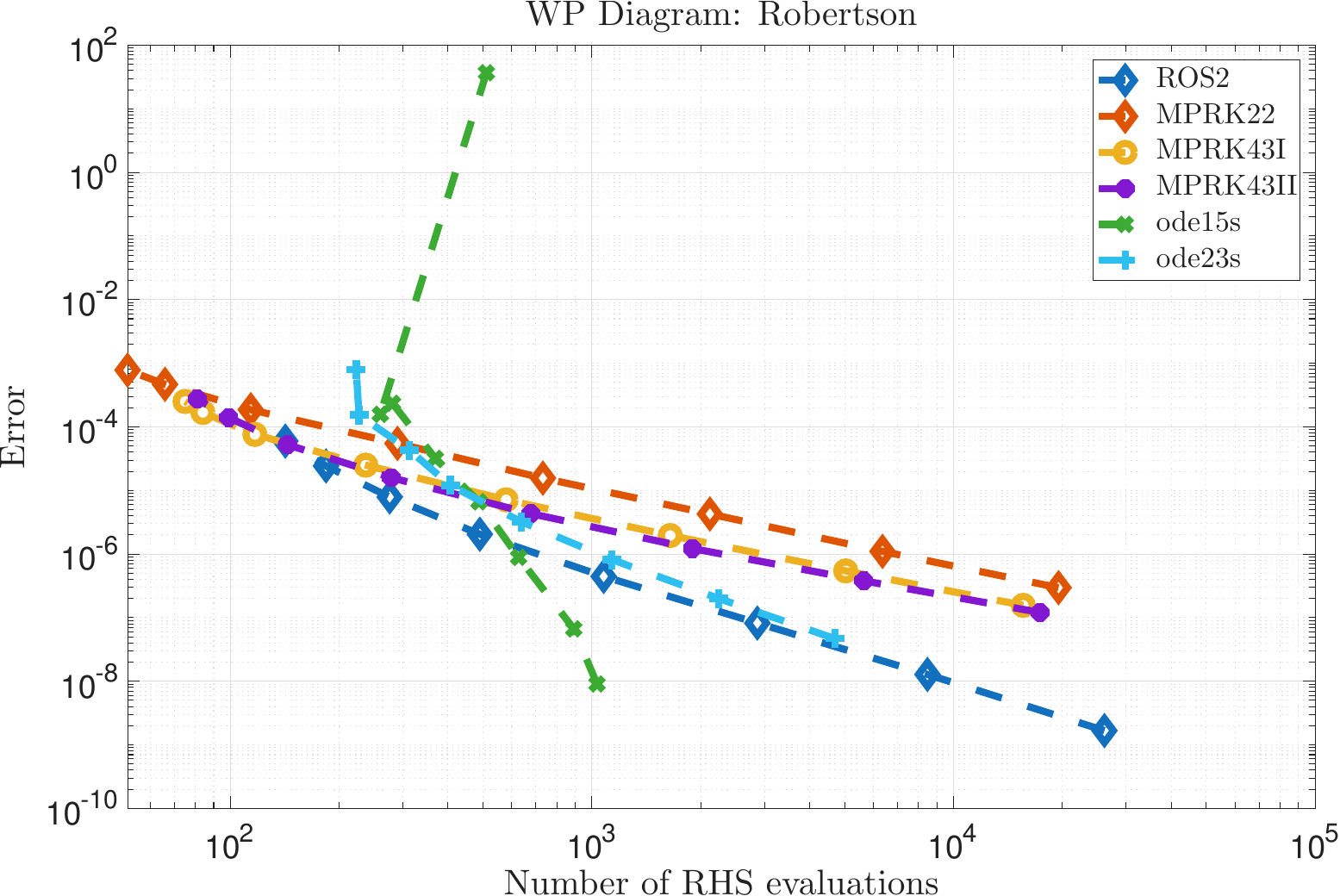}
				\end{subfigure}
				\begin{subfigure}[t]{0.495\textwidth}
					\includegraphics[width=\textwidth]{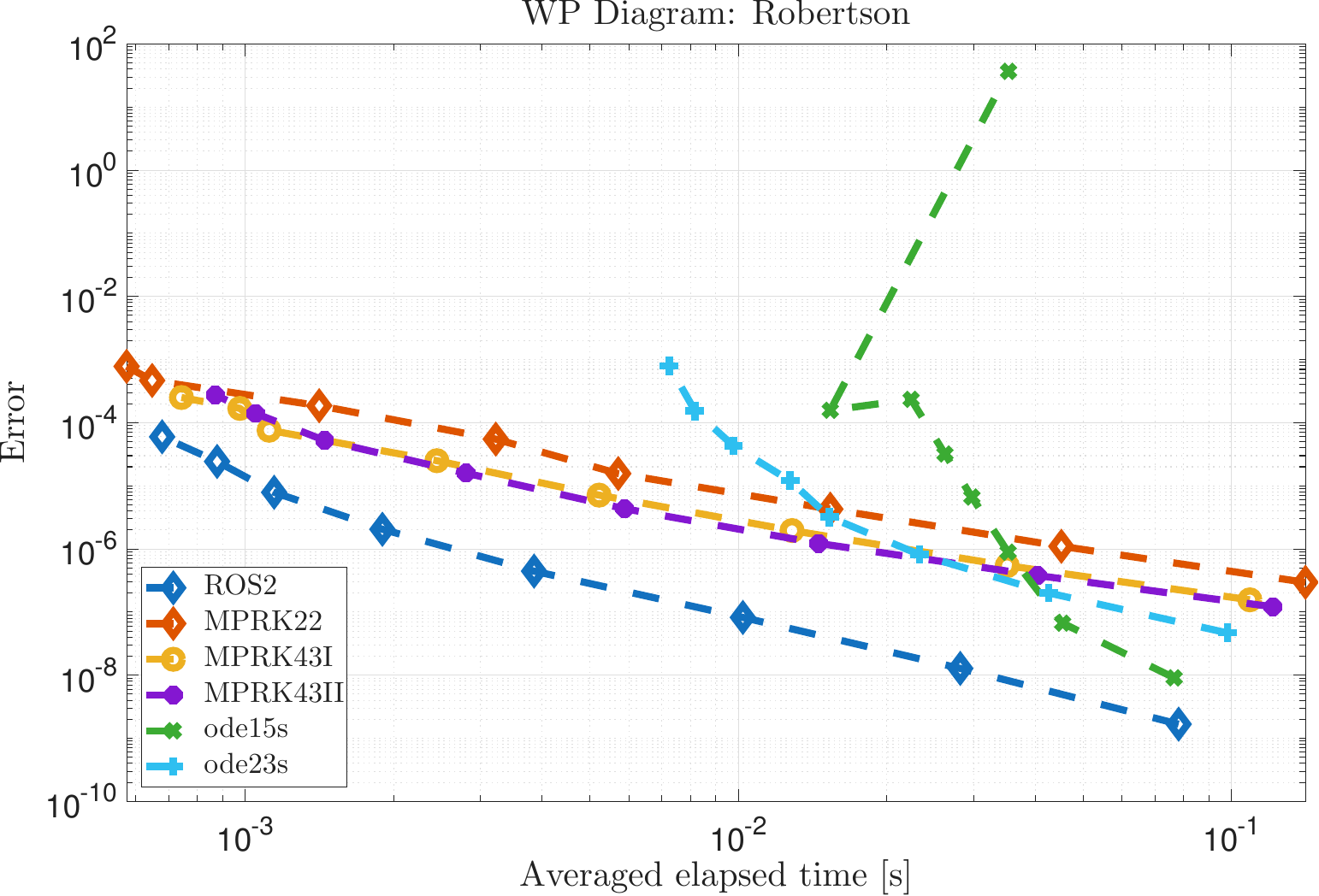}
				\end{subfigure}
												\begin{subfigure}[t]{0.495\textwidth}
					\includegraphics[width=\textwidth]{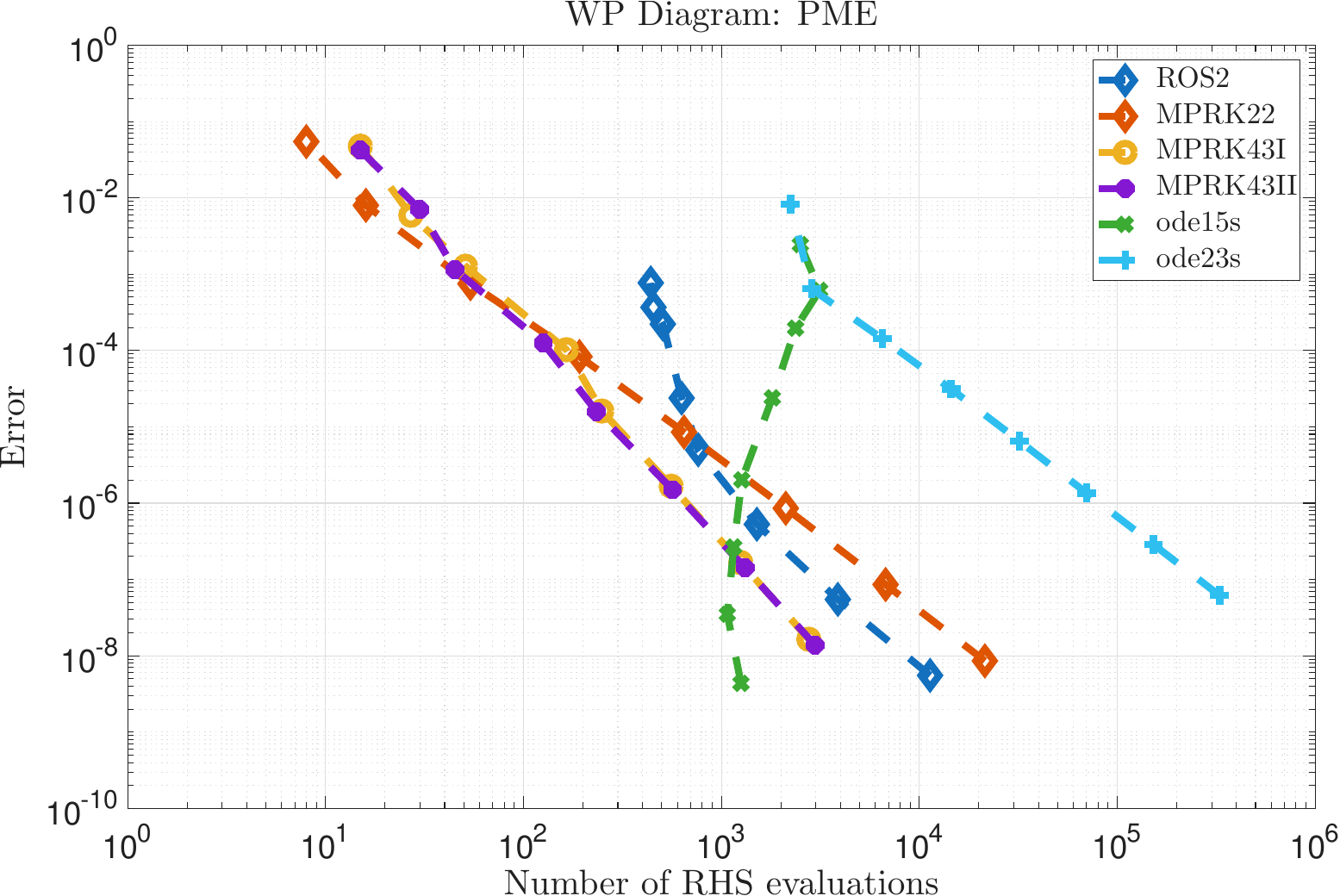}
				\end{subfigure}
				\begin{subfigure}[t]{0.495\textwidth}
					\includegraphics[width=\textwidth]{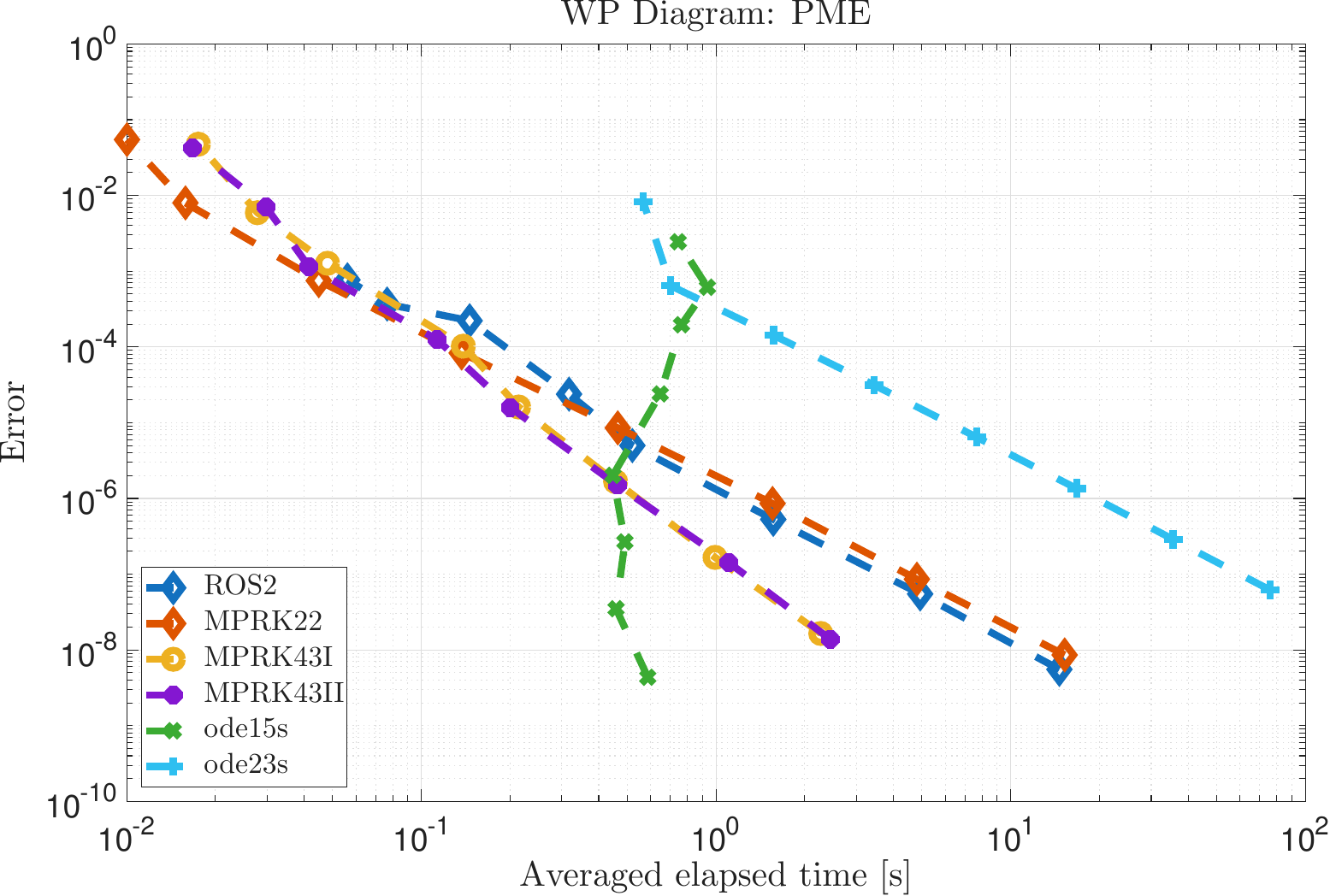}
				\end{subfigure}
					\caption{WP diagrams for the NPZD \eqref{eq:NPZD}, Robertson \eqref{eq:robertson}, and PME \eqref{eq:PME} problem  using the ROS2 and MPRK schemes with the DSP controllers from Table~\ref{tab:optiparas2}, and standard built-in \MATLAB solvers.}\label{Fig:WP_comp}
				\end{figure}

				The results are depicted in Figure~\ref{Fig:WP_comp} for NPZD and Robertson problems. The results for the remaining test cases did not provide further insights and are thus omitted for the sake of a compact presentation.

				In all cases, \code{ode15s} is the best in terms of RHS evaluations for the strictest tolerances. However, this superior performance may be traced back to the variable order of the method and should be further investigated in the future, \eg by comparing it to higher-order MP-type schemes.
				Nevertheless, for coarse and mid-ranged tolerances \code{ode15s} is beaten by MPRK and our Rosenbrock implementation. 



				We also note that in all test cases, our cost function provided a controller which is computationally stable, whereas this cannot be said about the used built-in ODE solvers in \MATLAB R2025b.
				This can be particularly seen for \code{ode15s} applied to the Robertson and PME problem with coarser tolerances, where the tolerance convergence can be observed only for stricter tolerances.

			In general, it may be expected that MPRK schemes are performing better than linear methods for coarser tolerances and test problems whenever the time step restrictions for a given tolerance are dominated by the demand for positive approximations.
			One example is given by the NPZD problem where a method producing negative approximations is in danger to diverge.
			Hence, while MPRK schemes are guaranteed to produce positive approximations for every time step, no matter how big, the built-in solvers, and our ROS2 scheme, can achieve this goal only by means of a time step reduction yielding a larger number of steps, and thus, RHS evaluations and elapsed time. Now, if the given tolerance is getting smaller, the time step size restrictions for positivity are fulfilled by the requirement of producing more accurate approximations, and hence, the benefit of MPRK schemes is lost.

			However, we compare here with the positivity-preserving variants of the built-in solvers.  Thus, we argue that MPRK schemes are worth considering whenever positivity of the numerical solution is interlinked with severe time step constraints. However, we discovered the same dominating behavior for all given test problems, which is why the reason for the better performance for coarser tolerances may lie elsewhere.

\section{Summary and Conclusion}\label{sec:summary}

We have developed an approach to design time step size controllers using
Bayesian optimization for modified Patankar--Runge--Kutta (MPRK) and Rosenbrock methods.
The basic idea is to use a general ansatz of step size controllers from
digital signal processing including PI and PID controllers. Then, we
search for a relevant set of test problems and design an appropiate cost
function used for the Bayesian optimization of the controller parameters.
Finally, the resulting controllers are validated for another set of
additional benchmarks including ODEs and PDEs.

The choice of test problems and the cost function is important to
obtain good results. This approach is comparable in spirit to
\cite{optimizedRK}, where a numerical search based on well-chosen test
problems was used in the context of designing adaptive time integration
methods for compressible computational fluid dynamics. As described therein,
the set of test problems needs to represent different regimes that can
appear in practicar applications, i.e., the problems should have different
characteristics such as stiffness properties. Focusing not only on MPRK methods,
we have chosen several differential equations linked to
chemical reactions and biological systems for this task.
In contrast to \cite{optimizedRK}, we have developed a cost function
to be used in Bayesian optimization to obtain the parameters of the
controllers. The design process of the cost function took several iterations
and took into account desirable properties such as computational stability
and tolerance convergence.
The cost function \eqref{eq:cost_fun} is developed specifically for
schemes such as MPRK methods and with good stability properties that are able
to handle even some stiff problems well. If the same approach was used
for other schemes such as explicit Runge--Kutta schemes, additional properties
such as step size control stability for mildly stiff problems should be
taken into account.
The extension of our approach to cases like this is left for future work.

Finally, we have applied the optimization-based approach to the second-order Rosenbrock scheme that is used for \code{ode23s} and three
MPRK schemes proposed in the literature:
the second-order method MPRK22$(1)$ \cite{KM18}
and the third-order methods MPRK43$(0.5,0.75)$ and MPRK43$(0.563)$
\cite{KM18Order3}. The optimized controllers for these methods are
summarized in Table~\ref{tab:optiparas2}.
A proof-of-concept implementation of the methods with the optimized
controllers are already comparable or even superior to \MATLAB solvers such as \code{ode23s}
for loose tolerances.

Along the way, we have also some new results on MPRK schemes that are
interesting on their own. In particular, we have extended MPRK schemes to
general time-dependent production-destruction-rest systems in
Section~\ref{sec:schemes}, broadening the scope of applications of
this class of methods.

\section*{Acknowledgments}

T.\ Izgin thanks Dr.~Andreas Linß for many fruitful discussions on the optimization problem.

\section*{Declarations}

\bmhead{Funding}
The author T.\ Izgin gratefully acknowledges the financial support by the Deutsche Forschungsgemeinschaft (DFG) through the grant ME 1889/10-1 (DFG project number 466355003).
H.~Ranocha was supported by the Deutsche Forschungsgemeinschaft
(DFG, German Research Foundation, project number 513301895)
and the Daimler und Benz Stiftung (Daimler and Benz foundation,
project number 32-10/22).

\bmhead{Conflict of interest}
The authors declare that they have no conflict of interest.

\bmhead{Availability of code, data, and materials}
The source code is available via the repository \cite{IR2025repository}.

\bmhead{Authors' contributions}
Conceptualization: Thomas Izgin, Hendrik Ranocha;
Data curation: Thomas Izgin;
Formal analysis and investigation: Thomas Izgin, Hendrik Ranocha;
Software: Thomas Izgin;
Visualization: Thomas Izgin;
Writing - original draft preparation: Thomas Izgin, Hendrik Ranocha;
Writing - review and editing: Thomas Izgin, Hendrik Ranocha.

\newpage
\appendix
\section*{Appendix}
\section{Data}\label{sec:data}
\begin{table}[!htpb]
	\begin{tabular}{llllllll}
		$\beta_1$ & $\beta_2$ & $\beta_3$ & $\alpha_2$ & $\kappa_2$ & Bayes & Costs & Time \\
		\hline
		1.5193 & -0.42576 & -0.078535 & -0.29465 & 2 & MPRK22adap\_3000 & 12.5187 & 0.037472 \\
		1.8476 & -0.11068 & -0.2863 & -0.24606 & 2 & MPRK43Iadap\_3000 & 12.5194 & 0.031583 \\
		1.9322 & -0.90687 & 0.21994 & -1.2967 & 3 & MPRK43IIadap\_3000 & 12.5221 & 0.030188 \\
		\textbf{1.4783} & \textbf{-1.027} & \textbf{-0.27743} & \textbf{-0.23112} & \textbf{3} & ROS2\_3000 & \textbf{2.279} & \textbf{15.2579} \\
		0.6 & -0.2 & 0 & 0 & 1 & - & 13.6323 & 0.56658 \\
		0.7 & -0.4 & 0 & 0 & 1 & - & 12.5167 & 0.032705 \\
		0.16667 & -0.33333 & 0 & 0 & 1 & - & 14.387 & 0.037796 \\
		0.16667 & 0.16667 & 0 & 0 & 1 & - & 12.5202 & 0.030455 \\
		1 & 0 & 0 & 0 & 1 & - & 12.5172 & 0.029884 \\
		2 & -1 & 0 & -1 & 1 & - & 12.5188 & 0.032984 \\
		0.5 & 0.5 & 0 & 0.5 & 1 & - & 13.7713 & 0.55764 \\
		0.056 & 0.111 & 0.056 & 0 & 1 & - & 12.5234 & 0.031419 \\
		0.25 & 0.25 & 0.25 & 0 & 1 & - & 12.518 & 0.0305 \\
		\hline
	\end{tabular}
	\caption{Data for ROS2}\label{tab:optiparasROS2}
\end{table}
\begin{table}[!htpb]
	\begin{tabular}{llllllll}
		$\beta_1$ & $\beta_2$ & $\beta_3$ & $\alpha_2$ & $\kappa_2$ & Bayes & Costs & Time \\
		\hline
		\textbf{1.5193} & \textbf{-0.42576} & \textbf{-0.078535} & \textbf{-0.29465} & \textbf{2} & MPRK22adap\_3000 & \textbf{3.6998} & \textbf{127.5217} \\
		1.8476 & -0.11068 & -0.2863 & -0.24606 & 2 & MPRK43Iadap\_3000 & 12.5454 & 0.037791 \\
		1.9322 & -0.90687 & 0.21994 & -1.2967 & 3 & MPRK43IIadap\_3000 & 12.5414 & 0.031324 \\
		1.4783 & -1.027 & -0.27743 & -0.23112 & 3 & ROS2\_3000 & 13.7706 & 1.0785 \\
		0.6 & -0.2 & 0 & 0 & 1 & - & 3.7463 & 121.068 \\
		0.7 & -0.4 & 0 & 0 & 1 & - & 3.7557 & 121.3198 \\
		0.16667 & -0.33333 & 0 & 0 & 1 & - & 14.3934 & 0.068642 \\
		0.16667 & 0.16667 & 0 & 0 & 1 & - & 14.3523 & 1.3355 \\
		1 & 0 & 0 & 0 & 1 & - & 12.5475 & 0.03051 \\
		2 & -1 & 0 & -1 & 1 & - & 3.7188 & 129.5072 \\
		0.5 & 0.5 & 0 & 0.5 & 1 & - & 12.5477 & 0.051973 \\
		0.056 & 0.111 & 0.056 & 0 & 1 & - & 14.352 & 1.3515 \\
		0.25 & 0.25 & 0.25 & 0 & 1 & - & 12.5491 & 0.031522 \\
		\hline
	\end{tabular}
	\caption{Data for MPRK22adap}\label{tab:optiparasMPRK22}
\end{table}
\begin{table}[!htpb]
	\begin{tabular}{llllllll}
		$\beta_1$ & $\beta_2$ & $\beta_3$ & $\alpha_2$ & $\kappa_2$ & Bayes & Costs & Time \\
		\hline
		1.5193 & -0.42576 & -0.078535 & -0.29465 & 2 & MPRK22adap\_3000 & 4.2665 & 24.8296 \\
		\textbf{1.8476} & \textbf{-0.11068} & \textbf{-0.2863} & \textbf{-0.24606} & \textbf{2} & MPRK43Iadap\_3000 & \textbf{4.2631} & \textbf{24.4657} \\
		1.9322 & -0.90687 & 0.21994 & -1.2967 & 3 & MPRK43IIadap\_3000 & 12.5584 & 0.037299 \\
		1.4783 & -1.027 & -0.27743 & -0.23112 & 3 & ROS2\_3000 & 4.4127 & 38.0836 \\
		0.6 & -0.2 & 0 & 0 & 1 & - & 4.3631 & 23.2243 \\
		0.7 & -0.4 & 0 & 0 & 1 & - & 4.36 & 22.9173 \\
		0.16667 & -0.33333 & 0 & 0 & 1 & - & 14.5452 & 0.13138 \\
		0.16667 & 0.16667 & 0 & 0 & 1 & - & 14.6537 & 0.81954 \\
		1 & 0 & 0 & 0 & 1 & - & 4.3014 & 23.6615 \\
		2 & -1 & 0 & -1 & 1 & - & 4.4895 & 80.7768 \\
		0.5 & 0.5 & 0 & 0.5 & 1 & - & 12.5671 & 0.038517 \\
		0.056 & 0.111 & 0.056 & 0 & 1 & - & 14.6664 & 0.84004 \\
		0.25 & 0.25 & 0.25 & 0 & 1 & - & 4.394 & 24.4239 \\
		\hline
	\end{tabular}
	\caption{Data for MPRK43Iadap}\label{tab:optiparasMPRK43I}
\end{table}
\begin{table}[!htpb]
	\begin{tabular}{llllllll}
		$\beta_1$ & $\beta_2$ & $\beta_3$ & $\alpha_2$ & $\kappa_2$ & Bayes & Costs & Time \\
		\hline
		\textbf{1.5193} & \textbf{-0.42576} & \textbf{-0.078535} & \textbf{-0.29465} & \textbf{2} & MPRK22adap\_3000 & \textbf{4.2745} & \textbf{27.2006} \\
		1.8476 & -0.11068 & -0.2863 & -0.24606 & 2 & MPRK43Iadap\_3000 & 4.2801 & 27.5557 \\
		1.9322 & -0.90687 & 0.21994 & -1.2967 & 3 & MPRK43IIadap\_3000 & 4.2785 & 51.7644 \\
		1.4783 & -1.027 & -0.27743 & -0.23112 & 3 & ROS2\_3000 & 4.393 & 45.0676 \\
		0.6 & -0.2 & 0 & 0 & 1 & - & 4.3659 & 25.2486 \\
		0.7 & -0.4 & 0 & 0 & 1 & - & 4.3712 & 25.6635 \\
		0.16667 & -0.33333 & 0 & 0 & 1 & - & 14.5454 & 0.11329 \\
		0.16667 & 0.16667 & 0 & 0 & 1 & - & 14.654 & 0.86466 \\
		1 & 0 & 0 & 0 & 1 & - & 4.3124 & 26.3163 \\
		2 & -1 & 0 & -1 & 1 & - & 4.4916 & 84.3145 \\
		0.5 & 0.5 & 0 & 0.5 & 1 & - & 4.4146 & 26.4771 \\
		0.056 & 0.111 & 0.056 & 0 & 1 & - & 14.6669 & 0.88539 \\
		0.25 & 0.25 & 0.25 & 0 & 1 & - & 16.7426 & 22.2278 \\
		\hline
	\end{tabular}
	\caption{Data for MPRK43IIadap}\label{tab:optiparasMPRK43II}
\end{table}
\section{Training Problems}\label{sec:training}
In this section we present several test problems for deriving optimized DSP parameters by the methodology described in the Section~\ref{sec:methodology}. There, the controller together with the numerical scheme is challenged to learn how to efficiently increase and decrease the step size in order to solve stiff problems.

Before proceeding, it is important to note that MPRK schemes need strictly positive initial data.
Hence, zeros in the initial vector will be replaced by $\code{realmin}$ which is around $10^{-308}$ for 64 bit floating point numbers used in the implementation.

\subsection{Prothero \& Robinson Problem}
\label{sec:PR4}
Introducing a $\mathcal C^1$-map $\b g$, we consider
\begin{equation}\label{eq:PR}
	\b y'(t)=\b \Lambda(\b y(t)-\b g(t))+\b g'(t), \quad \b y(0)=\b y^0\in \R^N_{>0}.
\end{equation}
We want to note that choosing $\b \Lambda\in \R^{N\times N}$ as a Metzler matrix satisfying $\bm 1\in \ker(\b \Lambda^T)$ and $\b g=\b 0$, this problem degenerates to \eqref{eq:Stabtestprob}. The problem \eqref{eq:PR} is related to the \mbox{Prothero \& Robinson} problem introduced in \cite{prothero74}.

If $\b g$ is smooth with $\b g(0)=\b y(0)$, then $\b y=\b g$ is the unique solution to the initial value problem. Therefore, choosing positive initial data and $\b g>\b 0$ guarantees the positivity of the solution of \eqref{eq:PR}.

In what follows, we solve the problem \eqref{eq:PR} for $N=4$ using a matrix $\b \Lambda_{\xi}\in \R^{4\times 4}$, whose spectrum runs over a vertical line in $\C^-$ as $\xi$ passes through the interval $[0,1]$. In particular,
\[\b \Lambda_\xi=\Vec{-1 & 1-\xi  & \xi &0\\ \xi &  -1 & 0 & 1-\xi\\ 1-\xi &  0 & -1 & \xi\\ 0 &  \xi & 1-\xi & -1} \]
satisfies $\sigma(\b \Lambda_\xi)=\{0, -2, (2\xi-1)\ii -1, (1-2\xi)\ii -1 \}$. Now according to \cite{StabMetz}, any eigenvalue $\lambda$ of a Metzler matrix $\b \Lambda\in \R^{4\times 4}$ with $\re(\lambda)=-1$ satisfies $\im(\lambda)\in [-1,1]$, so that we cover all possible eigenvalues when choosing $\xi\in[0,1]$. Indeed, due to symmetry, we only need to consider $\xi\in [0,1/2]$.
For the training we restrict to only one test case using $\xi=0,0.1,\dotsc,0.5$ and use the average cost in order to avoid an artificially large influence of this test problem on the cost function. However, we will use several values of $\xi$ for the validation of our results.

\begin{figure}
	\centering
	\includegraphics[width=.7\textwidth]{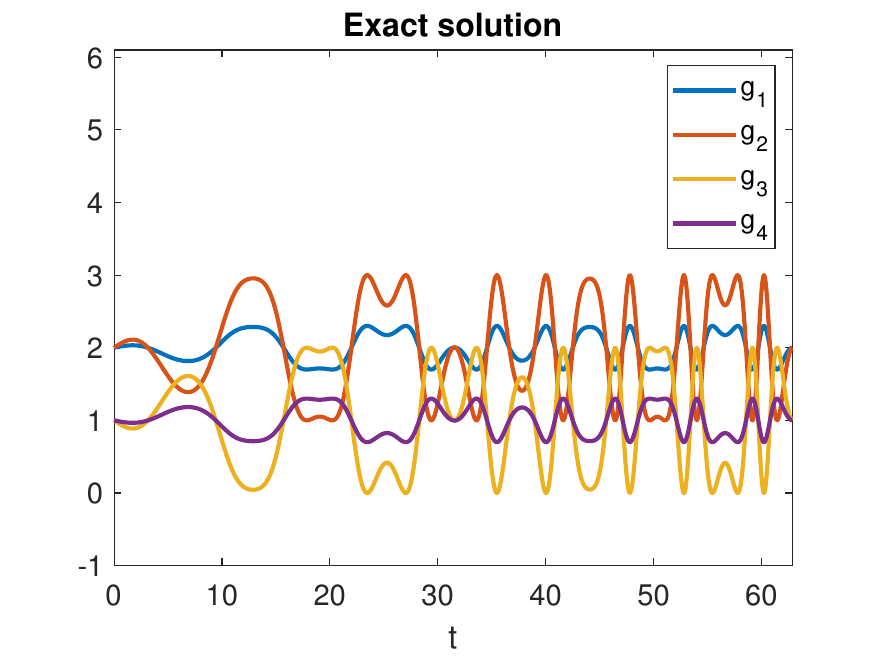}
	\caption{Graph of $\b g$ from \eqref{eq:gPR4} over the time interval $[0,20\pi]$.}\label{fig:graphg}
\end{figure}

Next, we use \begin{equation}\label{eq:gPR4}\b g(t)=\Vec{2+0.3\sin(0.5\cos(0.5t)t)\\ 2+\sin(0.5\cos(0.5t)t)\\1-\sin(0.5\cos(0.5t)t)\\1-0.3\sin(0.5\cos(0.5t)t)}\end{equation} and note that $g_1'(t)=-g_4'(t)$ and $g_2'(t)=-g_3'(t)$. The graph of $\b g$, \ie the solution of \eqref{eq:PR} can be seen in Figure~\ref{fig:graphg}. As one can observe, the time-dependent frequency results in an increasing amount of local maxima and minima of different magnitude.
Finally, the corresponding problem \eqref{eq:PR} will be approximated over the interval $[0,20\pi]$ can be written as a PDS using $p_{ij}=d_{ji}$ and
\begin{equation}
	\begin{aligned}
		p_{12}(\b y,t) &= y_2, &
		p_{13}(\b y,t) &= g_1(t),&
		p_{14}(\b y,t) &= \xi(y_3 + g_2(t)) + \min(0,g'_1(t)),\\
		p_{21}(\b y,t) &= g_2(t),&
		p_{24}(\b y,t) &= y_4,&
		p_{23}(\b y,t) &= \xi(g_4(t)+ y_1) + \min(0,g'_2(t)),\\
		p_{31}(\b y,t) &= y_1,&
		p_{34}(\b y,t) &= g_3(t),&
		p_{32}(\b y,t) &= \xi (g_1(t) + y_4) + \min(0,g'_3(t)),\\
		p_{42}(\b y,t) &= g_4(t),&
		p_{43}(\b y,t) &= y_3,&
		p_{41}(\b y,t) &= \xi (y_2 + g_3(t)) + \min(0,g'_4(t)).
	\end{aligned}
\end{equation}
In the following, we will refer to this problem as PR4 and use the initial time step size $\dt_0=1$.

\subsection{Robertson's Problem}
\begin{figure}[!htbp]
	\centering
	\includegraphics[width=.7\textwidth]{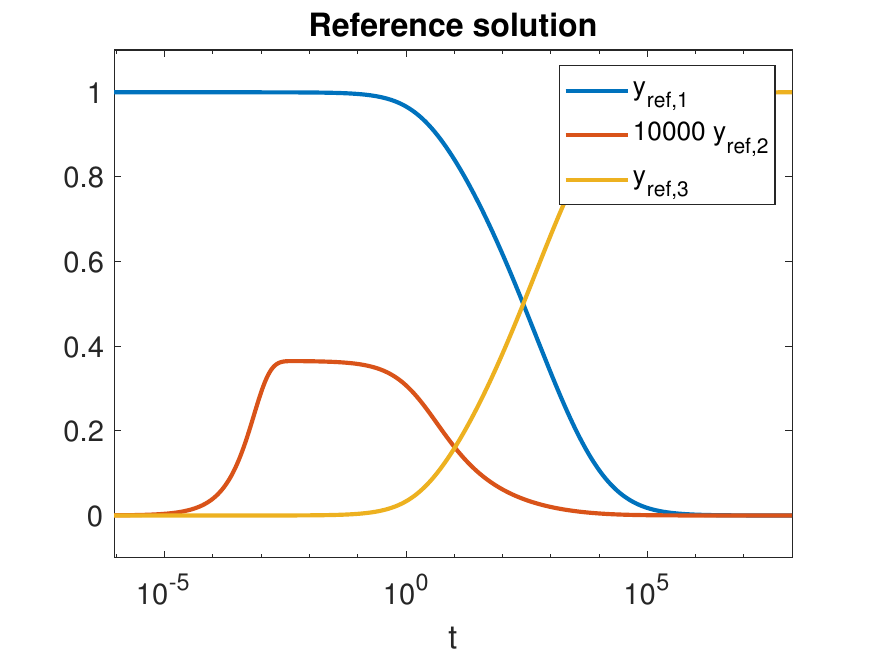}
	\caption{Reference solution of the Robertson problem \eqref{eq:robertson} with logarithmic $t$-axis over the interval $[0,10^{8}]$.}\label{fig:robertsonref}
\end{figure}

The  well-known Robertson problem is stiff \cite[Section~IV.1]{HNWII} and reads
\begin{equation}\label{eq:robertson}
	\begin{aligned}
		y_1'(t)&=10^4y_2(t)y_3(t)-0.04y_1(t),\\
		y_2'(t)&=0.04y_1(t)-10^4y_2(t)y_3(t)-3\cdot 10^7(y_2(t))^2,\\
		y_3'(t)&=3\cdot 10^7(y_2(t))^2.
	\end{aligned}
\end{equation}
We use the initial condition $\b y=(1,0,0)^T$ and initial time step size $\dt_0=10^{-6}$. A reference solution to this problem for $t\in [0,10^{8}]$ is depicted in Figure~\ref{fig:robertsonref}, where $y_2$ is multiplied with $10^4$ for visualization purposes. Moreover, this problem is positive and conservative, \ie $\b n=(1,1,1)^T$ satisfies $\b n^T\b y(t)=\const$.

\subsection{HIRES Problem}
\begin{figure}[!htbp]
	\centering
	\includegraphics[width=0.7\textwidth]{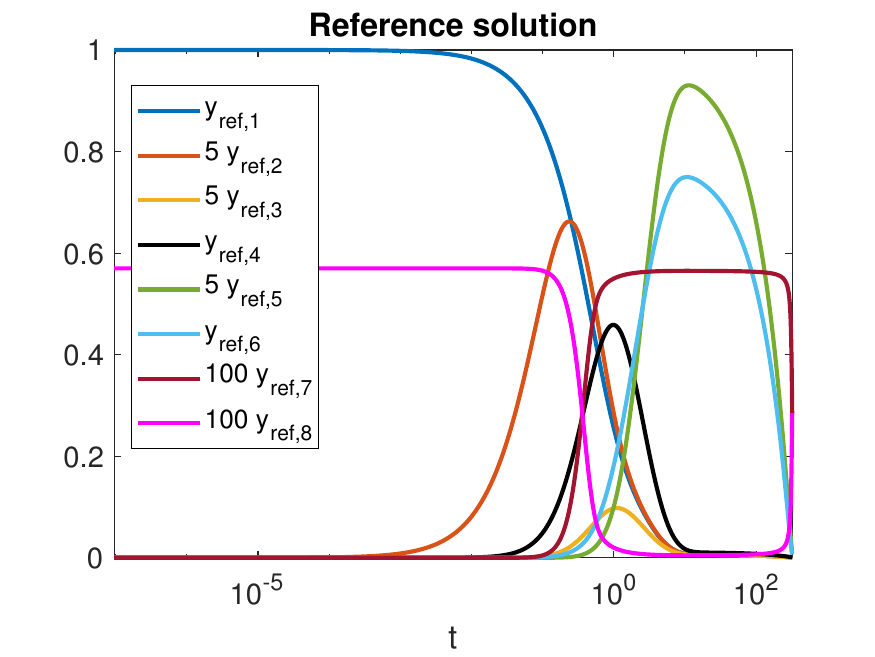}
	\caption{Reference solution of the HIRES problem \eqref{eq:hires} plotted over the time interval $[0,321.8122]$ with logarithmic $t$-axis. Thereby, the solution vector was scaled with $\diag(1,5,5,1,5,1,100,100)$.}\label{fig:hiresref}
\end{figure}

The "High Irradiance RESponse" (HIRES) problem, see \cite{S75,mazzia2003test,HNWII}, takes the form
\begin{equation}\label{eq:hires}
	\begin{aligned}
		y_1'&=-1.71y_1+0.43y_2+8.32y_3+0.0007,\\
		y_2'&=1.71y_1-8.75y_2,\\
		y_3'&= -10.03y_3+0.43y_4+0.035y_5,\\
		y_4'&=8.32y_2+1.71y_3-1.12y_4,\\
		y_5'&=-1.745y_5+0.43y_6+0.43y_7,\\
		y_6'&=-280y_6y_8+0.69y_4+1.71y_5-0.43y_6+0.69y_7,\\
		y_7'&= 280y_6y_8-1.81y_7,\\
		y_8'&=-280y_6y_8+1.81y_7.
	\end{aligned}
\end{equation}
The nonzero entries of the initial vector are $y_1(0)=1$ and $y_8(0)=0.0057$.
The system of ODEs can be written as a PDRS \eqref{eq:PDS_ODE} with $\b r^p(\b y)=(0.0007,0,0,0,0.43y_7,0.69y_7,0,0)^T$, $\b r^d(\b y)=(0,0,0,0,0,280y_6y_8,0,0)^T$ and
\begin{equation*}
	\begin{aligned}
		p_{12}(\b y) &= 0.43y_2, &
		p_{13}(\b y) &= 8.32y_3,&
		p_{21}(\b y) &= 1.71y_1,\\
		p_{34}(\b y) &=0.43y_4,&
		p_{35}(\b y) &=0.035y_5,&
		p_{42}(\b y) &=8.32y_2,\\
		p_{43}(\b y) &=1.71y_3,&
		p_{56}(\b y) &=0.43y_6,&
		p_{64}(\b y) &=0.69y_4 ,\\
		p_{65}(\b y) &=1.71y_5 ,&
		p_{78}(\b y) &=280y_6y_8,&
		p_{87}(\b y) &=1.81y_7
	\end{aligned}
\end{equation*}
being the non-vanishing production terms and $d_{ij}=p_{ji}$. The problem will be approximated on the interval $[0,321.8122]$ using the initial time step size $\dt_0=0.5\cdot 10^{-3}$, and a reference solution can be found in Figure~\ref{fig:hiresref}.
This problem is interesting since it is a non-conservative PDS and the controller has to learn how to handle the different behavior of the solution on the respective time scales.
\subsection{NPZD Problem}
\label{sec:NPZD}
\begin{figure}[!htbp]
	\centering
	\includegraphics[width=0.7\textwidth]{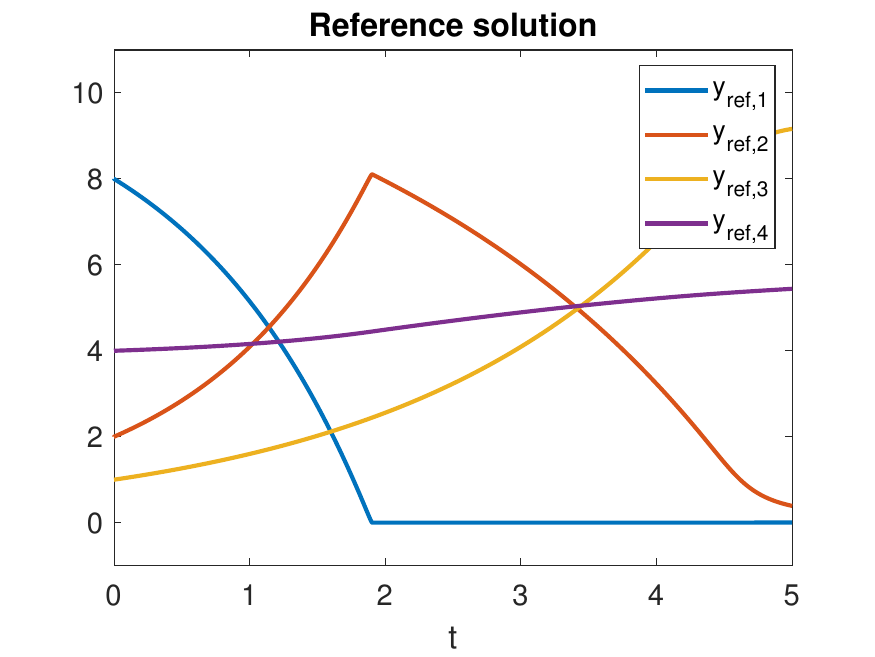}
	\caption{Reference solution of the NPZD problem \eqref{eq:NPZD} plotted over the time interval $[0,5]$ in hours.}\label{fig:npzdref}
\end{figure}

The NPZD problem
\begin{equation}\label{eq:NPZD}
	\begin{aligned}
		y_1'&=0.01y_2+0.01y_3+0.003y_4-\frac{y_1y_2}{0.01+y_1},\\
		y_2'&=\frac{y_1y_2}{0.01+y_1}-0.06y_2-0.5(1-\exp(-1.21y_2^2))y_3,\\
		y_3'&= 0.5(1-\exp(-1.21y_2^2))y_3 -0.03y_3,\\
		y_4'&=0.05y_2 +0.02y_3-0.003y_4
	\end{aligned}
\end{equation}
from \cite{BDM2005} models the interaction of Nutrients, Phytoplankton, Zooplankton and Detritus. The problem \eqref{eq:NPZD} can be written as a PDS using
\begin{equation*}
	\begin{aligned}
		p_{12}(\b y)&=0.01y_2, & p_{13}(\b y)&=0.01y_3, \quad p_{14}(\b y)=0.003y_4, \\
		p_{21}(\b y)&=\frac{y_1y_2}{0.01+y_1}, & p_{32}(\b y)&=0.5(1-\exp(-1.21y_2^2))y_3,   \\
		p_{42}(\b y)&=0.05y_2, & p_{43}(\b y)&=0.02y_3
	\end{aligned}
\end{equation*}
and $d_{ij}=p_{ji}$ while the remaining production and destruction terms are zero. A reference solution to the NPZD problem is depicted in Figure~\ref{fig:npzdref}.
It was demonstrated in \cite{BDM2005} that the occurrence of negative approximations\footnote{There is a small margin of acceptable negative approximations, see \cite{BDM2005} for the details.} leads to the divergence of the method, and hence, severe time step restrictions to schemes not preserving the positivity unconditionally are necessary. For the numerical experiments we use an initial time step size of $\dt_0=1$.

\section{Validation Problems}\label{sec:validation}
Here we introduce several validation problems. In particular, we consider the Prothero and Robinson problem with other parameters and the Brusselator problem.

\subsection{Porous Medium Equation}
The porous medium equation
\begin{equation}\label{eq:PME}
	u_t = (u^m)_{xx} = (a(u) u_x)_x, \quad a(u) = m u^{m-1}
\end{equation}
with a free parameter $m>1$, see for instance \cite{Boscarino23}, admits a non-negative weak solution
\[u^{(m)}(t,x)=t^{-k}\left[\max\left(1-\frac{k(m-1)}{2m}\frac{\lvert x\rvert^2}{t^{2k}},0\right)\right]^{\frac{1}{m-1}}\]
with $k=\frac{1}{m+1}$, the so-called \emph{Barenblatt} solution \cite{Barenblatt52}. For every $t>0$, the solution has a compact support $[-\alpha_m(t),\alpha_m(t)]$ where
\[\alpha_m(t)=\sqrt{\frac{2m}{k(m-1)}}t^k.\]
We follow \cite{Boscarino23} using $u(0,x)=u^{(m)}(1,x)$ as an initial condition. We plot the numerical solution at time $t=2$ on the spatial domain $[-6,6]$ using the boundary conditions $u(t,\pm 6)=0$ for $t>1$.

We use the second-order space discretization from \cite{Mattsson12,ranocha2019mimetic} given by
\[f_i(\b u(t))=\frac{a(u_i(t))+a(u_{i+1}(t))}{2\dx^2} u_{i+1}(t)-\frac{a(u_{i-1}(t))+2a(u_i(t))+a(u_{i+1}(t))}{2\dx^2}u_i(t)+\frac{a(u_{i-1}(t))+a(u_i(t))}{2\dx^2} u_{i-1}(t)\]
for $i=2,\dotsc,N$ and
\[f_j(\b u(t))= \frac{a(u_j(t))}{2\dx^2}u_j(t),\quad \text{ for }\quad j\in \{1,N\}.\]
Next, we consider the convex entropy
\[\eta(\b u)=\frac{\dx^2}2\sum_{i=1}^{N_x}u_i^2, \]
which satisfies
\begin{equation*}
	\begin{aligned}
		\frac{\dd}{\dd t}\eta(\b u(t))\leq 0
	\end{aligned}
\end{equation*}
for the boundary conditions mentioned above, see \cite[Theorem~4.1]{ranocha2019mimetic}.
This system of ODEs may be rewritten as a conservative PDS by setting
\begin{equation*}
	\begin{aligned}
		p_{i,i+1}(\b u)&=\frac{a(u_i)+a(u_{i+1})}{2\dx^2} u_{i+1},&\quad p_{i,i-1}(\b u)&= \frac{a(u_{i-1})+a(u_i)}{2\dx^2} u_{i-1},&\quad i&=2,\dotsc,N,\\
		p_{1,2}(\b u)&= \frac{a(u_2)}{2\dx^2}u_2,&\quad p_{N,N-1}(\b u)&= \frac{a(u_{N-1})}{2\dx^2}u_{N-1},&\quad d_{i,j}&=p_{j,i}.
	\end{aligned}
\end{equation*}

According to \cite{Boscarino23}, the cases $m=3,5$ are particularly interesting as the numerical solution of the proposed third-order IMEX method in \cite[p.~10, eq.~(30)]{Boscarino23} generates negative approximations and which cannot happen with MPRK schemes. We particularly choose $m=3$ for the numerical experiments together with $N=200$ and $\dt_0=0.8\dx$ and solve until $t=2$.
\begin{figure}[!htbp]
	\centering
	\includegraphics[width=0.7\textwidth]{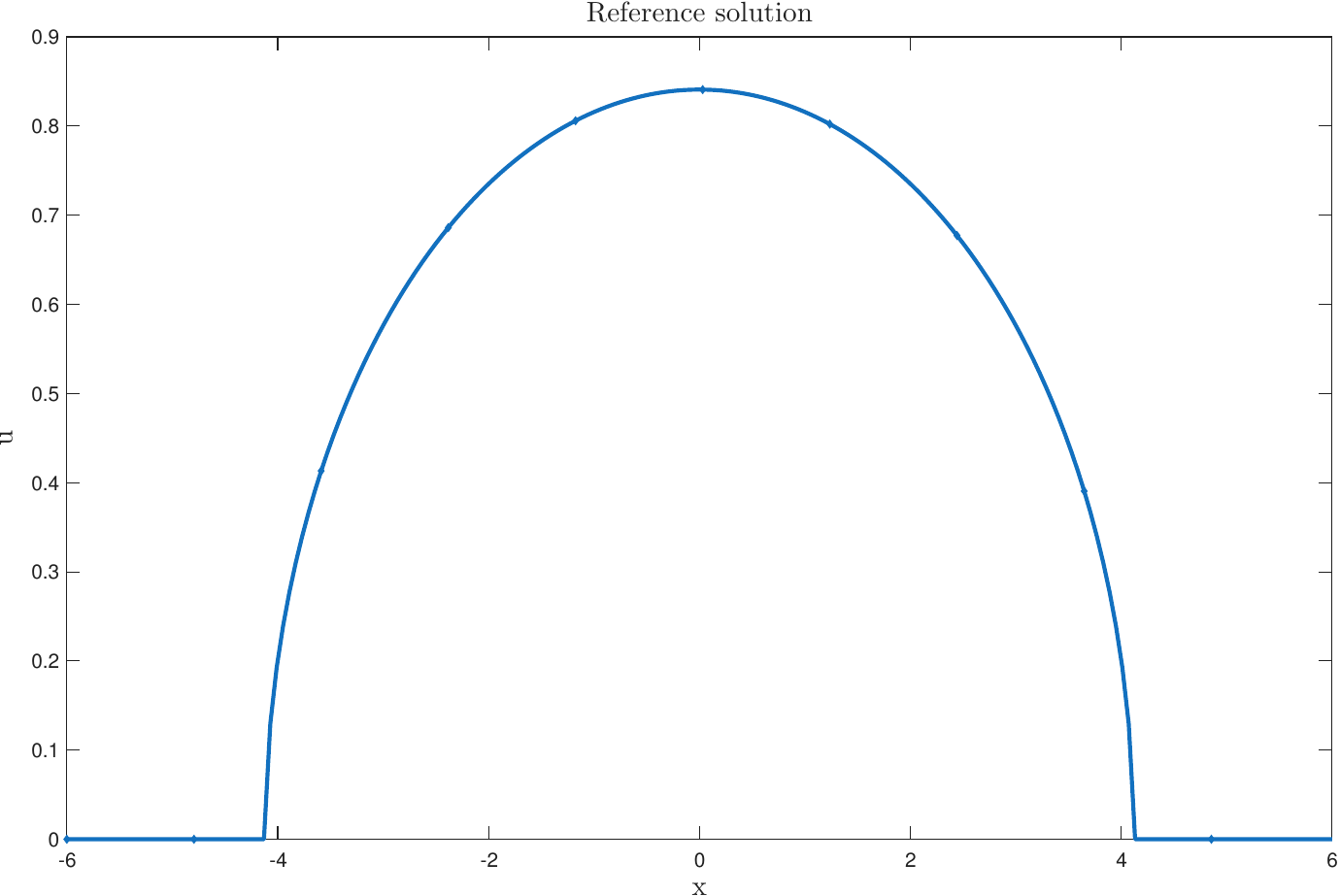}
	\caption{Reference solution of the Porous Medium Equation  \eqref{eq:PME} with $m=3$ plotted over the time interval $[0, 2]$.}\label{fig:PMEref}
\end{figure}
\subsection{Fokker--Planck}
Following \cite{PZ2018,BLR2024}, we consider the Fokker--Planck equation in the form
\begin{equation}\label{eq:FP}
	\partial_tu(w,t) =\partial_w[B[u](w,t)u(w,t)+\partial_w(D(w)u(w,t))],\quad u(w,0)=u^0(w)
\end{equation}
together with no-flux boundary conditions on $I=[-1,1]$. Setting $h_{i+\frac12}\coloneqq h(w_{i+\frac12})$ for any continuous function $h=h(w)$, where $w_i$ denote grid points, $w_{i+\frac12}$ denote cell interfaces, and $\Delta_w$ is the uniform distance between two grid points. In particular, $w_\frac12 = -1$, $w_{N+\frac12}=1$, and $w_i=\frac{\Delta_w}{2}+w_{i-\frac12}$ for $i=1,\dotsc,N$, where $N$ is the number of grid points.
Next, we set $u_i(t)\coloneqq u(w_i,t)$ and apply the Chang--Cooper method to obtain the semidiscretisation
\begin{equation*}
	\frac{\dd u_i(t)}{\dd t} = \frac{F_{i+\frac12}(t)-F_{i-\frac12}(t)}{\Delta_w},\quad i=1,\dotsc,N
\end{equation*}
where 		 $F_\frac12 = F_{N+\frac12} = 0$,
\begin{equation*}
	\begin{aligned}
		F_{i+\frac12}= C_{i+\frac12}((1-\delta_{i+\frac12})u_{i+1}+\delta_{i+\frac12}u_i)+D_{i+\frac12}\frac{u_{i+1}-u_i}{\Delta_w}, \quad i=1,\dotsc,N-1,
	\end{aligned}
\end{equation*}
and
\[C_{i+\frac12} = \frac{\lambda_{i+\frac12}D_{i+\frac12}}{\Delta_w},\quad \delta_{i+\frac12}=\frac{1}{1-\exp(\lambda_{i+\frac12})}+\frac{1}{\lambda_{i+\frac12}},\quad \lambda_{i+\frac12}=\frac{\Delta_w(B[u](w_{i+\frac12})+D'_{i+\frac12})}{D_{i+\frac12}}. \]
According to \cite{BLR2024}, this can be written in form of a conservative PDS, where the nonzero production and destruction terms are given by
\begin{equation*}
	\begin{aligned}
		p_{i,i+1}(\b u) &= d_{i+1,i}(\b u) \coloneqq \frac{\max(0,C_{i+\frac12}((1-\delta_{i+\frac12})u_{i+1}+\delta_{i+\frac12}u_i +\frac{D_{i+\frac12}u_{i+1}}{\Delta_w}}{\Delta_w},\\
				p_{i,i-1}(\b u) &= d_{i-1,i}(\b u) \coloneqq \frac{-\min(0,C_{i-\frac12}((1-\delta_{i-\frac12})u_i+\delta_{i-\frac12}u_{i-1} +\frac{D_{i-\frac12}u_{i-1}}{\Delta_w}}{\Delta_w}
	\end{aligned}
\end{equation*}
for $i=2,\dotsc,N-1$. Within the numerical experiments, we use $\dt_0=0.8\Delta_w$ with $N=200$, and solve until $t=10$. The initial condition is given by
\[u^0(w)=\beta\left(\exp(-30(w+\tfrac12)^2)+\exp(-30(w-\tfrac12)^2) \right),\] where $\beta$ is chosen such that $\int_I u^0\dd w=1$.
\begin{figure}[!htbp]
		\begin{subfigure}[t]{0.45\textwidth}
			\includegraphics[width=1\textwidth]{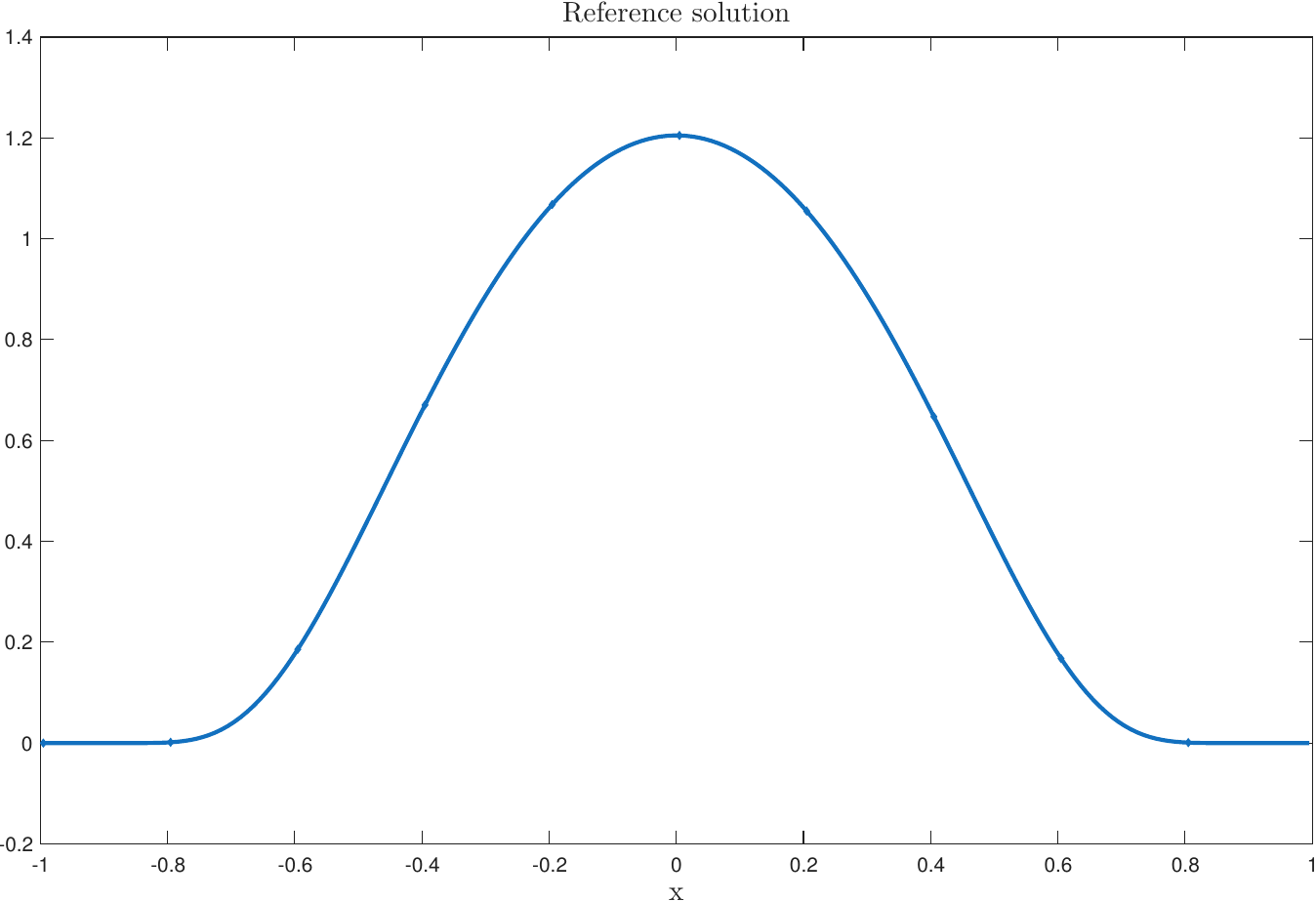}
		\end{subfigure}
		\begin{subfigure}[t]{0.6\textwidth}
			\includegraphics[width=\textwidth]{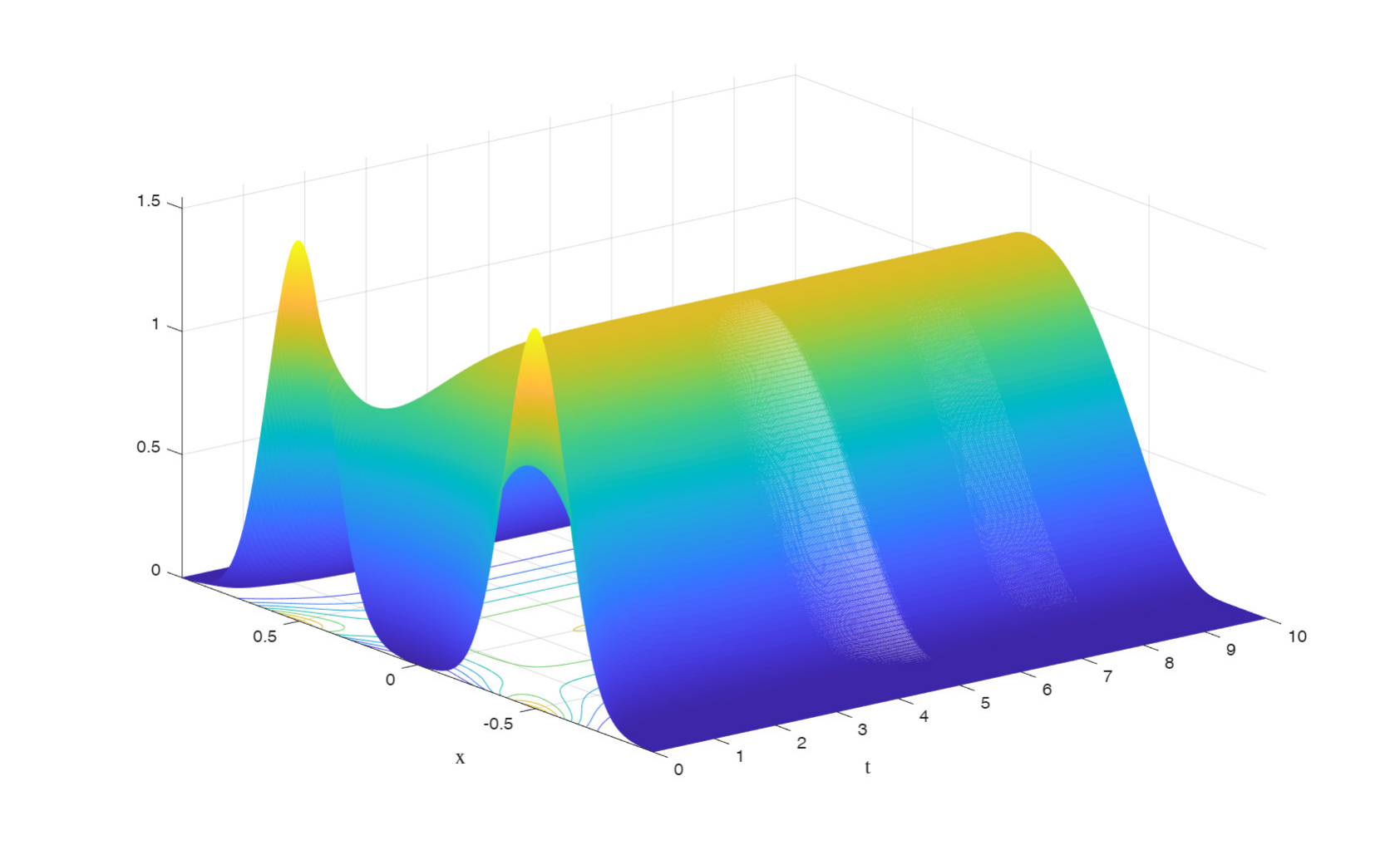}
		\end{subfigure}
	\caption{Reference solution of the Fokker--Planck equation \eqref{eq:FP} plotted over the time interval $[0, 10]$.}\label{fig:FP}
\end{figure}
\subsection{Prothero \& Robinson Problem}
As we have constructed infinitely many problems \eqref{eq:PR}, it is an obvious decision to also take some of the test cases as a validation problem. In particular, we use $\xi=0.25$.
\subsection{Brusselator Problem}
\begin{figure}[!htbp]
	\centering
	\includegraphics[width=0.7\textwidth]{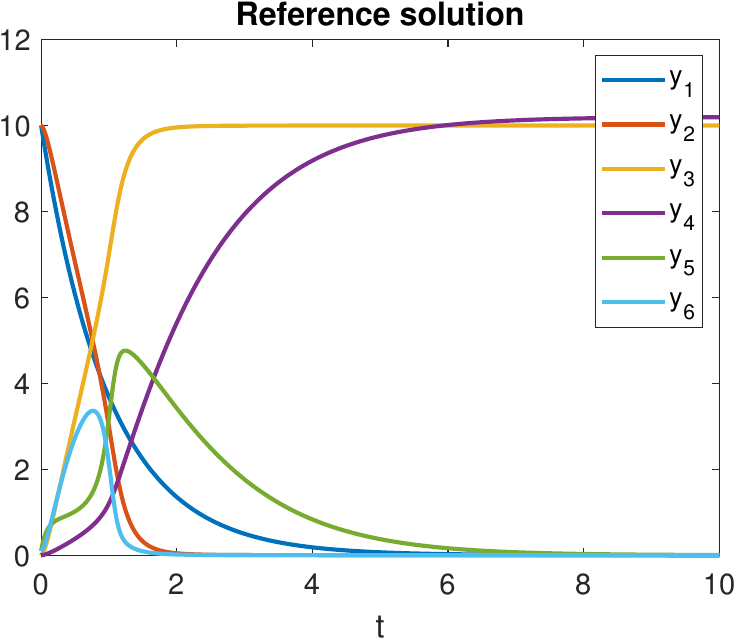}
	\caption{Reference solution of the Brusselator problem \eqref{eq:bruss} plotted over the time interval $[0, 10]$.}\label{fig:brussref}
\end{figure}

As a non-stiff validation problem we consider the Brusselator problem \cite{LefeverNicolis1971,HNW1993} which reads
\begin{equation}\label{eq:bruss}
	\begin{aligned}
		y_1'&=-k_1y_1,\\
		y_2'&=-k_2y_2y_5,\\
		y_3'&= k_2y_2y_5,\\
		y_4'&=k_4y_5,\\
		y_5'&=k_1y_1-k_2y_2y_5+k_3y_5^2y_6-k_4y_5,\\
		y_6'&=k_2y_2y_5-k_3y_5^2y_6,
	\end{aligned}
\end{equation}
where we set $k_i>0$.  In form of PDS the Brusselator problem is determined by
\begin{equation*}
	\begin{aligned}
		p_{32}(\b y) &=k_2 y_2y_5, &
		p_{45}(\b y) &=k_4 y_5,&
		p_{51}(\b y) &= k_1y_1,\\
		p_{56}(\b y) &=k_3y_5^2y_6,&
		p_{65}(\b y) &=k_2y_2y_5
	\end{aligned}
\end{equation*}
and $d_{ij}=p_{ji}$. We set $k_i=1$ and use the initial condition $\b y(0)=(10,10,0,0,0.1,0.1)^T$. The time interval of interest is $[0,10]$, and the reference solution can be seen in Figure~\ref{fig:brussref}. The numerical method start with a time step size of $\dt_0=0.1$.

\bibliography{cas-refs}

\end{document}